\documentclass[DIV=12]{scrartcl}
\usepackage{amsmath}
 \numberwithin{equation}{section}
\usepackage{amssymb}
\usepackage{amsthm}

\theoremstyle{definition}

\theoremstyle{remark}

\usepackage[main=english,german]{babel}
\usepackage[noisbn,noissn,fixlanguage]{babelbib}
 \setbtxfallbacklanguage{english}
\usepackage{cite}
\usepackage[T1]{fontenc}
\usepackage{hyperref}
\usepackage{lmodern}
\usepackage{mathscinet}
\usepackage{mathrsfs}
\usepackage{mathtools}
 \mathtoolsset{mathic=true}
\usepackage{microtype}
\usepackage{pmat}

\usepackage{211arxiv}

\title{The Weyl matrix balls corresponding to the matricial truncated Hamburger moment problem}
\author{Bernd Fritzsche \and Bernd Kirstein \and Susanne Kley \and Conrad M\"adler}
\date{}

\begin{document}
\maketitle

\begin{abstract}
 The main goal of the paper is to parametrize the Weyl matrix balls associated with an arbitrary matricial truncated Hamburger moment problem.
 For the special case of a non-degenerate matricial truncated Hamburger moment problem the corresponding Weyl matrix balls were computed by I.~V.~Kovalishina \cite{MR703593} in the framework of V.~P.~Potapov's method of `Fundamental matrix inequalities'.
\end{abstract}

\begin{description}
 \item[Mathematics Subject Classification (2010)] 44A60, 47A57, 30E05.
 \item[Keywords] Matricial truncated Hamburger moment problem, orthogonal matrix polynomials, parametrization of the solution set via Stieltjes transform, Weyl matrix balls.
\end{description}

\section{Introduction}
 This paper is closely related to the authors' recent investigations \cite{FKKM20} on the matricial truncated Hamburger moment problem.
 In \cite{FKKM20}, the authors obtained a parametrization of the Stieltjes transforms of the solution set in terms of a linear fractional transformation.
 This parametrization works for the general matrix case of the problem under consideration.
 Our main goal in this paper is to use results of \cite{FKKM20} to achieve a description of the set of values of matrices which will be attained by the Stieltjes transforms of the solutions of the moment problem at a prescribed fixed point of the open upper half plane.
 We show that this set of values fills a (closed) matrix ball and present explicit expressions for the center and the semi-radii of this matrix ball, which will also be called the Weyl matrix ball associated with the concrete matricial Hamburger moment problem under consideration.
 Following the classical monograph Akhiezer \cite{MR0184042}, the terminology ``Weyl circles'' or later ``Weyl matrix balls'' was consequently used in the Soviet literature (see, \teg{}\ \cite{MR0425671,MR647177,MR703593,MR752056,MR751390,MR1473266}).
 As explained in \zitaa{MR0184042}{\cch{1}{}} or in \cite{MR2077204}, the history of the scalar case is intimately related to the classical papers \cite{MR1511560,MR1512075,zbMATH02604576,MR1503221}.
 Hellinger \cite{MR1512075} mentioned the analogies of his considerations on the moment problem with  H.~Weyl's paper \cite{MR1511560} on a boundary value problem for ordinary differential equations.
 In 1934, H.~Weyl studies connections between the papers \cite{MR1512075} and \cite{MR1511560}.
 However, he did not consider the moment problem, but the Nevanlinna--Pick problem for holomorphic functions in the open upper half plane, the imaginary part of which is non-negative.
 The main purpose in \cite{MR1503221} is to formulate that problem in the area of ordinary differential equations, which corresponds to the Nevanlinna--Pick interpolation problem.
 What concerns the matrix case, one can observe two lines of investigations.
 The first one is connected with particular types of ordinary differential equations.
 It is a generalization of the topic opened by H.~Weyl in his fundamental paper \cite{MR1511560} on Sturm--Liouville equations on the semi-axis.
 In this framework, we mention the landmark paper \cite{MR0425671} by S.~A.~Orlov.
 The central theme of this paper can be considered as part of the second line of investigations in the matrix case.
 Starting in the 1970's, this line was essentially initiated  by the research of the school of V.~P.~Potapov on matricial generalizations of classical interpolation and moment problems.
 What concerns the matricial version of the non-degenerate Hamburger moment problem, the corresponding matrix balls were explicitly computed by I.~V.~Kovalishina in \cite{MR703593}.
 Her approach is based on V.~P.~Potapov's method of ``Fundamental matrix inequalities'' (FMI method).
 
 This paper is part of the investigations of the first and second authors on matricial versions of classical interpolation and moment problems in the general (possibly degenerate) matrix case.
 In the first period of this research together with A.~Lasarow, they concentrated on interpolation problems for holomorphic matrix functions in the open unit disc, namely the matrix versions of the classical interpolation problems named after Carath\'eodory and Schur.
 The key instrument of this approach is in both problems the so-called central solution.
 The central features of this work can be summarized by the three steps:

 1) Analysis of the structure of the central solution (see \cite{MR2104258,MR2066852})
 
 2) Inspired by the first step construction of a particular polynomial matrix which parametrizes the solution set via a linear fractional transformation of the Schur class (see \cite{MR2222523,MR2493510})
 
 3) Closer analysis of the blocks of the generating matrix polynomial realizing the parametrization of the solution in order to achieve information on the concrete shape of parameters of the associated Weyl matrix ball (see \cite{MR2656833,MR2390673})
 
 This strategy could be also used to handle the corresponding interpolation problem for \(J\)\nobreakdash-Potapov functions.
 This was done in collaboration with U.~Raabe and K.~Sieber (see \cite{MR2535572,MR2647541,MR2746081,MR2890903}).
 Later the authors started their studies on matricial versions of truncated power moment problems on the real axis in the general matrix case.
 There the holomorphic matrix functions are defined in the open upper half plane.
 
 What are the differences to the former investigation in the case of the unit disk and which suggestions can be maintained.
 The first difference now comes from the absence of an analogue of the central solution.
 Remember that the existence of this solution essentially governed what we have done in the unit disk case.
 This forces us to find a convenient way to find a matrix polynomial, called resolvent matrix in the sequel, which produces the parametrization of the solution set (, more precisely, the set of their corresponding Stieltjes transforms).
 In \cite{FKKM20}, we obtained the resolvent matrix on the basis of a Schur type algorithm.
 This procedure provided the four \tqqa{blocks} of the resolvent matrix in a form which is not convenient for obtaining information about the corresponding matrix balls.
 For this reason, in this paper, we will take a closer look at the resolvent matrix.
 In particular, we will recognize that these \tqqa{blocks} can be interpreted as parts of a corresponding quadruple of orthogonal \tqqa{matrix} polynomials.
 This observation enables us now to continue in a way suggested by \cite{MR2222523}.
 More precisely, the blocks of the resolvent matrix suggest to consider now a distinguished rational matrix function, which proves to be a key instrument for the subsequent considerations.
 Roughly speaking, we obtain now a conjecture about the parameters of the Weyl matrix ball we are looking for.
 This conjecture will be verified then.
 
 This paper is organized as follows.
 In \rsec{S-II}, we summarize some basic facts on matricial Hamburger moment problems.
 In particular, we recall the well-known characterizations of the existence of solutions (see \rthmss{SolvHam}{T0708}).
 Following the classical line, in \rsec{S-III}, we reformulate the original moment problem as an interpolation problem for appropriate classes of holomorphic matrix-valued functions in the open upper half plane \(\pip\).
 For this reason, we discuss several integral representations of these functions (see \rthmss{Theo21}{FKMM16_32}).
 The main aim of \rsec{S-IV} is to provide essential ingredients for the parametrization of the set of Stieltjes transforms of all solutions.
 More precisely, we introduce a special class of ordered pairs of meromorphic \tqqa{matrix}-valued functions in \(\pip\), which serves as set of parameters later.
 In \rsec{S1111}, we take a closer look at the structure  of \tHnnde{} (sometimes also called Hankel non-negative definite extendable) sequences \(\seqs{2n}\) of complex \tqqa{matrices} and we introduce the \tsHp{} (see \rdefn{CM2.1.64}).
 We consider the so-called \tnatext{} \(\seqsh{2n+1}\) of \(\seqs{2n}\) (see \rdefn{D.nat-ext}).
 In \zitaa{FKKM20}{\cthm{8.11}{}}, we obtained a parametrization of the set of all Stieltjes transforms of the solution set via  a linear fractional transformation of meromorphic pairs.
 The generating matrix-valued function is a \taaa{2q}{2q}{matrix} polynomial, which was constructed via a Schur type algorithm.
 For the purposes of this paper, it is essential to recognize that the four \tqqa{blocks} are parts of a \tabcd{}.
 This is done in \rsec{S-VII} (see in particular \rthm{CM123}).
 The investigations in \rsec{S1129} are guided by the experiences from \cite{MR2222523}.
 As in \cite{MR2222523}, by a closer look to the resolvent matrix we are led to a sequence of rational matrix-valued functions (see \rdefn{K17.1}), which contain the key information about the Weyl matrix balls we are striving for.
 These functions turn out to be related by our Schur type algorithm (see \rlem{P1056}).
 This observation provides in combination with \zitaa{MR3380267}{\cprop{8.6}{}} that this function belongs to the \tqqa{Herglotz}--Nevanlinna class (see \rprop{P1129}).
 \rsec{Cha13} contains the central result of this paper (see \rthm{T45CD}).
 More precisely, it is a description of the Weyl matrix ball.
 In \rsec{S1031}, we indicate how the matrix ball description, which I.~V.~Kovalishina \cite{MR703593} obtained in the non-degenerate case, can be derived within our approach which works for the most general case.

\section{Notation and preliminaries}\label{S-II}

 First we state some notation.
 Let \(\C\), \(\R\), \(\NO\), and \(\N\) be the set of all complex numbers, the set of all real numbers, the set of all \tnn{} integers, and the set of all positive integers, respectively.
 Further, for every choice of \(\alpha,\beta\in\R\cup\set{-\infty,\infty}\), let \(\mn{\alpha}{\beta}\) be the set of all integers \(k\) such that \(\alpha\leq k\leq\beta\).
 Throughout this paper, if not explicitly mentioned otherwise, then let \(p,q,r\in\N\).
 If \(\mathcal{X}\) is a \tne{} set, then \(\mathcal{X}^\pq\) represents the set of all \tpqa{matrices} each entry of which belongs to \(\mathcal{X}\), and \(\mathcal{X}^p\) is short for \(\mathcal{X}^{p\times 1}\).
 The notation \(\CHq\) is used to denote the set of all \tH{} complex \tqqa{matrices}.
 We write \(\Cggq\) and \(\Cgq\) to designate the set of all \tnnH{} complex \tqqa{matrices} and the set of all positive \tH{} complex \tqqa{matrices}, respectively.
 
 Let \(\OA\) be a measurable space.
 Then each countably additive mapping defined on \(\fA\) with values in \(\Cggq\) is called a \tnnH{} \tqqa{measure} on \(\OA\) and the notation \(\Mggqa{\Omega,\fA}\) stands for the set of all \tnnH{} \tqqa{measures} on \(\OA\).
 Let \(\mu=\matauuuo{\mu_{jk}}{j}{k}{1}{q}\) be a \tnnH{} \tqqa{measure} on \(\OA\).
 Then we use \(\LOC{\mu}\) to denote the set of all Borel measurable functions \(f\colon\Omega\to\C\) for which \(\int_\Omega\abs{f}\dif\nu_{jk}<\infty\) holds true for every choice of \(j\) and \(k\) in \(\mn{1}{q}\), where \(\nu_{jk}\) is the variation of the complex measure \(\mu_{jk}\) (see also \rrem{M.8-1}).
 If \(f\in\LOC{\mu}\), then let \(\int_\Omega f\dif\mu\defeq\matauuuo{\int_\Omega f\dif\mu_{jk}}{j}{k}{1}{q}\) and we also write \(\int_\Omega f(\omega)\mu\rk{\dif\omega}\) for this integral.

 Denote by \(\BsA{\R}\) (\tresp{}\ \(\BsA{\C}\)) the \(\sigma\)\nobreakdash-algebra of all Borel subsets of \(\R\) (or \(\C\), respectively).
 Let \(\Omega\in\BsA{\R}\setminus\set{\emptyset}\).
 Then designate by \(\BsA{\Omega}\) the \(\sigma\)\nobreakdash-algebra of all Borel subsets of \(\Omega\) and by \(\Mggqa{\Omega}\) the set of all \tnnH{} \tqqa{measures} on \(\rk{\Omega,\BsA{\Omega}}\), \tie{}, \(\Mggqa{\Omega}\) is short for \(\Mggqa{\Omega,\BsA{\Omega}}\).
 Let \(\kappa\in\NOinf\).
 Then denote by \(\Mgguoa{\kappa}{q}{\Omega}\) the set of all \(\sigma\in\Mggqa{\Omega}\) such that, for all \(j\in\mn{0}{ \kappa}\), the function \(f_j\colon\Omega\to\C\) defined by \(f_j(\omega)\defeq\omega^j\) belongs to \(\Ls{1}{\Omega}{\BsA{\Omega}}{\sigma}{\C}\).
 If \(\sigma\in\Mgguoa{\kappa}{q}{\Omega}\), then, for all \(j\in\mn{0}{\kappa}\), let \(\mpm{\sigma}{j}\defeq\int_\Omega\omega^j\sigma\rk{\dif\omega}\).
 We will consider the following types of a matricial Hamburger power moment problem:

\begin{problem}[\mprobR{\kappa}{=}]
 Let \(\kappa\in\NOinf\) and let \(\seqska\) be a sequence of complex \tqqa{matrices}.
 Parametrize the set \(\MggqRsg{\kappa}\) of all \(\sigma\in\MgguqR{\kappa}\) fulfilling \(\su{j}=\mpm{\sigma}{j}\) for all \(j\in\mn{0}{\kappa}\).
\end{problem}

\begin{problem}[\mprobR{2n}{\lleq}]
 Let \(n\in\NO\) and let \(\seqs{2n}\) be a sequence of complex \tqqa{matrices}.
 Parametrize the set \(\MggqRskg{2n}\) of all \(\sigma\in\MgguqR{2n}\) for which the matrix \(\su{2n}-\mpm{\sigma}{2n}\) is \tnnH{} and, in the case \(n\geq1\), for which additionally \(\su{j}=\mpm{\sigma}{j}\) is fulfilled for all \(j\in\mn{0}{2n-1}\).
\end{problem}

 For our further consideration, we introduce certain sets of sequences of complex \tqqa{matrices} which are determined by properties of particular \tbHms{} built of them.
 If \(n\in\NO\) and if \(\seqs{2n}\) is a sequence of complex \tqqa{matrices}, then \(\seqs{2n}\) is called \emph{\tHnnd{}} (\emph{\tHpd{}}, respectively) if the \tbHm{}
\beql{N2}
 \Hn
 \defeq\mat{\su{j+k}}_{j,k =0}^n
\eeq
 is \tnnH{} (positive \tH{}, respectively).
 (Note that \tHnnd{} (\tHpd{}, respectively) sequences of complex \tqqa{matrices} are also said to be Hankel non-negative definite (Hankel positive definite, respectively).)
 For all \(n\in\NO\), we will write \(\Hggq{2n}\) (or \(\Hgq{2n}\), respectively) for the set of all sequences \(\seqs{2n}\) of complex \tqqa{matrices} which are \tHnnd{} (\tHpd{}, respectively).
 A well-known solvability criterion for Problem~\mprobR{2n}{\lleq} is the following:

\bthml{SolvHam}
 Let \(n\in\NO\) and let \(\seqs{2n}\) be a sequence of complex \tqqa{matrices}.
 Then \(\MggqRskg{2n}\neq\emptyset\) if and only if \(\seqs{2n}\in\Hggq{2n}\).
\ethm 

 There are various proofs of \rthm{SolvHam}, for example \cite[\cthm{3.2}{}]{MR1624548} 
 and \cite[\cthm{4.16}{795}]{MR2570113}.
 If \(n\in\NO\) and if \(\seqs{2n}\in\Hggq{2n}\) (or \(\seqs{2n}\in\Hgq{2n}\), respectively), then, for each \(m\in\mn{0}{n}\), the sequence \(\seqs{2m}\) obviously belongs to \(\Hggq{2m}\) (or \(\Hgq{2m}\), respectively).
 Thus, let \(\Hggqinf\) (or \(\Hgqinf\), respectively) be the set of all sequences \(\seqsinf\) of complex \tqqa{matrices} such that, for all \(n\in\NO\), the sequence \(\seqs{2n}\) belongs to \(\Hggq{2n}\) (or \(\Hgq{2n}\), respectively).
 For all \(n\in\NO\), let \(\Hggeq{2n}\) (\tresp{}\ \(\Hgeq{2n}\)) be the set of all sequences \(\seqs{2n}\) of complex \tqqa{matrices} for which there exist complex \tqqa{matrices} \(\su{2n+1}\) and \(\su{2n+2}\) such that \(\seqs{2(n+1)}\) belongs to \(\Hggq{2\rk{n+1}}\) (\tresp{}\ \(\Hgq{2\rk{n+1}}\)).
 Furthermore, for all \(n\in\NO\), we will use \(\Hggeq{2n+1}\) (\tresp{}\ \(\Hgeq{2n+1}\)) to denote the set of all sequences \(\seqs{2n+1}\) of complex \tqqa{matrices} for which there exists a complex \tqqa{matrix} \(\su{2n+2}\) such that \(\seqs{2(n+1)}\) belongs to \(\Hggq{2\rk{n+1}}\) (\tresp{}\ \(\Hgq{2\rk{n+1}}\)).
 For each \(m\in\NO\), the elements of the set \(\Hggeq{m}\) are called \emph{\tHnnde{}} (or Hankel non-negative definite extendable) sequences.
 For technical reasons, we set \(\Hggeqinf\defeq\Hggqinf\) and \(\Hgeqinf\defeq\Hgqinf\).
 A well-known solvability criterion for Problem~\mprobR{\kappa}{=} is the following:

\bthml{T0708}
 Let \(\kappa\in\NOinf\) and let \(\seqska\) be a sequence of complex \tqqa{matrices}.
 Then \(\MggqRsg{\kappa}\neq\emptyset\) if and only if \(\seqska\in\Hggeqka\).
\ethm
 
 A proof of \rthm{T0708} is given, \teg{}, \zitaa{MR2805417}{\cthm{6.6}{486}}.

\breml{Schr2.7}
 According to \zitaa{MR2570113}{\crem{2.29}{} and \cprop{2.30}{}}, we have \(\Hgeq{2\tau}=\Hgq{2\tau}\subseteq\Hggeq{2\tau}\subseteq\Hggq{2\tau}\) and \(\Hgq{2\tau}\neq\Hggeq{2\tau}\) for all \(\tau\in\NOinf\).
 Furthermore, \(\Hgeq{\kappa}\subseteq\Hggeq{\kappa}\) for all \(\kappa\in\NOinf\) as well as \(\Hggeq{0}=\Hggq{0}\), whereas \(\Hggeq{2\tau}\neq\Hggq{2\tau}\) for all \(\tau\in\Ninf\).
\erem

 The following result is essential for a parametrization of the set \(\MggqRskg{2n}\):

\bthml{SK7}
 Let \(n\in\NO\) and let \(\seqs{2n}\in\Hggq{2n}\).
 Then there exists a unique sequence \(\seqa{\tilde{s}}{2n}\in\Hggeq{2n}\) such that \(\MggqRakg{\seqa{\tilde{s}}{2n}}=\MggqRskg{2n}\).
\ethm

 The existence of such a sequence \(\seqa{\tilde{s}}{2n}\) was formulated first in \cite[\clem{2.12}{}]{MR1395706}.
 A complete proof of \rthm{SK7} can be found in \cite[\cthm{7.3}{806--808}]{MR2570113}.
 The main goal of \cite{MR3889658} is to present a purely matrix theoretical object which yields, applied to a special case, the explicit construction of the desired sequence \(\seqa{\tilde{s}}{2n}\in\Hggeq{2n}\).
 
 At the end of this introductory section, we give some further notation.
 We will write \(\Iq\) to denote the identity matrix in \(\Cqq\), whereas \(\Opq\) is the zero matrix belonging to \(\Cpq\).
 Sometimes, if the size is clear from the context, we will omit the indices and write \(\EM\) and \(\NM\), \tresp{}
 For each \(A\in\Cpq\), let \(\ran{A}\) be the column space of \(A\), let \(\nul{A}\) be the null space of \(A\), and let \(\rank A\) be the rank of \(A\).
 For each \(A\in\Cqq\), we will use \(\rre A\) and \(\rim A\) to denote the real part and the imaginary part of \(A\), respectively:
 \(\rre A\defeq\frac{1}{2}(A+A^\ad)\) and \(\rim A\defeq\frac{1}{2\iu}(A-A^\ad)\).
 Furthermore, for each \(A\in\Cpq\), let \(\normS{A}\) be the operator norm of \(A\).
 A complex \tpqa{matrix} \(A\) is said to be contractive if \(\normS{A}\leq1\).
 We use \(\Kpq\) in order to designate the set of all contractive complex \tpqa{matrices}.
 If \(A\in\Cqq\), then \(\det A\) denotes the determinant of \(A\). %
 For each \(A\in\Cqp\), let \(A^\mpi\) be the Moore--Penrose inverse of \(A\).
 If \(p_1,p_2,q_1,q_2\in\N\) and if \(A_j\in\Coo{p_j}{q_j}\) for every choice of \(j\in\set{1,2}\), then let \(\diag(A_1,A_2)\defeq\smat{A_1 & \Ouu{p_1}{q_2}\\ \Ouu{p_2}{q_1} & A_2}\).
 Furthermore, within the set \(\CHq\), we use the L\"owner semi-ordering:
 If \(A\) and \(B\) are complex \tH{} \tqqa{matrices}, then we will write \(A\lleq B\) (or \(B\lgeq A\)) to indicate that \(B-A\) is a \tnnH{} matrix.
 
 For all \(x,y\in\Cq\), by \(\ipE{x}{y}\) we denote the (left-hand side) Euclidean inner product of \(x\) and \(y\), \tie{}, we have \(\ipE{x}{y}\defeq y^\ad x\).
 If \(\mM\) is a \tne{} subset of \(\Cq\), then let \(\mM^\oc\) be the set of all vectors in \(\Cq\) which are orthogonal to \(\mM\) (with respect to the Euclidean inner product \(\ipE{.}{.}\)).
 If \(\mU\) is a linear subspace of \(\Cq\), then let \(\OPu{\mU}\) be the orthogonal projection matrix onto \(\mU\).
 
 Let \(\dom\) be a a \tne{} open subset of the complex plane.
 If \(f\) is a complex-valued function meromorphic in \(\dom\), then we use \(\holpt{f}\) to denote the set of all points \(w\) at which \(f\) is holomorphic and we use the notation
\beql{NM}
 \nst{f}
 \defeq\setaca{w\in\holpt{f}}{f\rk{w}=0}.
\eeq
  A \tpqa{matrix}-valued function \(F=\mat{f_{jk}}_{\substack{j=1,\dotsc,p\\k=1,\dotsc,q}}\) is said to be meromorphic in \(\dom\) if \(f_{jk}\) is meromorphic in \(\dom\) for each \(j\in\mn{1}{p}\) and each \(k\in\mn{1}{q}\).
  In this case, we signify \(\holpt{F}\defeq\bigcap_{j=1}^p\bigcap_{k=1}^q\holpt{f_{jk}}\).

\section{Particular classes of holomorphic matrix functions}\label{S-III}
 We will reformulate the matricial moment problems under consideration as equivalent interpolation problems for particular classes of holomorphic matrix-valued functions.
 For this reason, we introduce in this section the corresponding function classes and summarize some of their essential properties needed in the sequel.
 Most of this material is taken from \cite{MR2988005,MR3380267}.
 Let \(\pip\defeq\setaca{z\in\C}{\rim z\in(0,\infty)}\) be the upper half-plane of the complex plane.
 The class \(\RFq\) of all  \emph{\tqqa{Herglotz--Nevanlinna} functions in \(\pip\)} consists of all matrix-valued functions \(F\colon\pip\to\Cqq\) which are holomorphic in \(\pip\) and which satisfy \(\rim F\rk{z}\in\Cggq\) for all \(z\in\pip\).
 Detailed observations about matrix-valued Herglotz--Nevanlinna functions can be found in \cite{MR1784638,MR2988005}.
 Especially, the functions belonging to \(\RFq\) admit a well-known integral representation.
 Before we formulate this matricial generalization of a famous result due to R.~Nevanlinna, we observe that, for every choice of \(\nu\in\MggqR\) and \(z\in\CR\), the function \(f_z\colon\R\to\C\) given by \(f_z(x) \defeq\rk{1 + xz}/\rk{x-z}\) belongs to \(\LRC{\nu}\).

\bthmnl{Nevanlinna}{Theo21}
\benui
 \item For each \(F\in\RFq\), there exists a unique triple \((\alpha, \beta, \nu)\in\CHq  \times \Cggq  \times \MggqR\) such that 
\begin{align}\label{FKMM16_31_1}
 F\rk{w}&=\alpha +\beta w +\int_{\R}\frac{1 + xw}{x - w}\nu \rk{\dif x}&\text{for each }w&\in\pip.
\end{align}
 \item If \(\alpha\in\CHq\), if \(\beta\in\Cggq\), and if \(\nu\in\MggqR\), then \(F\colon\pip\to\Cqq\) defined by \eqref{FKMM16_31_1} belongs to \(\RFq\).
\eenui
\ethm

 For each \(F\in\RFq\), the unique triple \((\alpha, \beta, \nu)\in\CHq  \times \Cggq  \times \MggqR\) for which the representation \eqref{FKMM16_31_1} holds true is called the \emph{\tNpo{\(F\)}} and we also write \((\alpha_F, \beta_F, \nu_F)\) instead of \((\alpha, \beta, \nu)\).
 
\blemnl{See \zitaa{MR1784638}{\cthm{5.4(iv)}{}}}{fkm12bP3.3}%
 If \(F\in\RFq\), then \(\beta_F = \lim_{y\to\infty} \rk{\iu y}^\inv F\rk{\iu y}\).
\elem
 
 For our consideration, the class \(\RFOq\) given by
\[
 \RFOq
 \defeq\setaca*{ F\in\RFq }{ \sup_{y\in[1,\infty)} y\normS*{F\rk{\iu y}}<\infty},
\]
 plays a key role.
 The functions belonging to \(\RFOq\) admit a special integral representation as well.
 Before we formulate this result, let us observe that, for every choice of \(\sigma\in\MggqR\) and \(z\in\C\setminus\R\), one can easily see that the function \(h_{z}\colon\R\to\C\) given by \(h_{z}(x)\defeq\rk{x-z}^\inv\) describes a bounded and continuous function which, in particular, belongs to \(\LRC{\sigma}\).
 Now we formulate the well-known matricial generalization of another classical result due to R.~Nevanlinna \cite{zbMATH02604576}:

\bthml{FKMM16_32}
\benui
 \il{FKMM16_32.a} For each \(F\in\RFOq\), there is a unique \(\sigma\in\MggqR\) such that 
\begin{align}\label{FKMM16_N1140D}
 F\rk{w}&=\int_\R\frac{1}{x -w}\sigma\rk{\dif x}&\text{for each }w&\in\pip.
\end{align}
 \il{FKMM16_32.b} If \(\sigma\in\MggqR\), then \(F\colon\pip\to\Cqq\) defined by \eqref{FKMM16_N1140D} belongs to \(\RFOq\).
\eenui
\ethm

 \rthm{FKMM16_32} can be proved by using its well-known scalar version in the case \(q=1\) as well as \rrem{M.8-1}, \zitaa{MR2988005}{\clem{B.3}{1788}}, and the fact that, for each \(F\in\RFOq\) and each \(u\in\Cq\), the function \(u^\ad Fu\) belongs to \(\RFOu{1}\).
 If \(F\in\RFOq\), then the unique \(\sigma\in\MggqR\) for which \eqref{FKMM16_N1140D} holds true is the so-called \emph{\tRSm{} (or matricial spectral measure) of \(F\)} and we also write \(\smF\) instead of \(\sigma\).
 If \(\sigma\in\MggqR\) is given, then \(F\colon\pip\to\Cqq\) defined by \eqref{FKMM16_N1140D} is said to be the \emph{\tRSto{\(\sigma\)}}.

\breml{R16}
 In view of \rthm{FKMM16_32}, now one can reformulate Problems~\mprobR{\kappa}{=} and~\mprobR{2n}{\lleq} in the language of \tRSt{s}:
\begin{problem}[\rprobR{\kappa}{=}]
 Let \(\kappa\in\NOinf\) and let \(\seqska\) be a sequence of complex \tqqa{matrices}.
 Parametrize the set \(\RFOqsg{\kappa}\) of all matrix-valued functions \(F\in\RFOq\) the \tRSm{} of which belongs to \(\MggqRsg{\kappa}\).
\end{problem}
\begin{problem}[\rprobR{2n}{\lleq}]
 Let \(n\in\NO\) and let \(\seqs{2n}\) be a sequence of complex \tqqa{matrices}.
 Parametrize the set \(\RFOqskg{2n}\) of all matrix-valued functions \(F\in\RFOq\) the \tRSm{} of which belongs to \(\MggqRskg{2n}\).
\end{problem}
\erem

 Note that parametrizations of the solution sets of Problems~\rprobR{\kappa}{=} and~\rprobR{2n}{\lleq} in the general case are given in \cite{MR1624548,MR1740433,MR1395706,MR3380267,Thi06,FKKM20}.

\bleml{CS314}
 If \(F\in\RFOq\), then, for each \(z\in\pip\), the equations \(\ran{F\rk{z}}=\ran{\sigma_F\rk{\R}}\) and \(\nul{F\rk{z}}=\nul{\sigma_F\rk{\R}}\) 
 hold true.
\elem

 There is a proof of \rlem{CS314}, \teg{}, in \cite[\clem{8.2}{1785} and \cprop{8.9}{1786}]{MR2988005}.
 Let
\[
 \RFoq{-2}
 \defeq\setaca*{F\in\RFq}{\lim_{y\to\infty}\rk*{\frac{1}{y}\normS*{F\rk{\iu y}}}=0}.
\]
 
\bexanl{\tcf{}~\zitaa{MR3380267}{\cexa{3.3}{190}}}{E0520}
 Let \(E\in\Cqq\) be such that \(\rim E\in\Cggq\).
 Then \(F\colon\pip\to\Cqq\) defined by \(F(z)\defeq E\) belongs to \(\RFoq{-2}\).
\eexa

 For all \(A\in\Cpq\), let \(\Pqevena{A}\defeq\setaca{F\in\RFoq{-2}}{\nul{A}\subseteq\nul{\alpha_F}\cap\Nul{\nu_F\rk{\R}}}\).

\bleml{L1654}
 Let \(E\in\Cqq\) be such that \(\rim E\in\Cggq\) and let \(A\in\Cpq\) be such that \(\nul{A}\subseteq\nul{E}\).
 Then \(F\colon\pip\to\Cqq\) defined by \(F(z)\defeq E\) belongs to \(\Pqevena{A}\).
\elem
\bproof
 From \rexa{E0520} we see \(F\in\RFoq{-2}\).
 In particular, \(F\in\RFq\).
 Therefore, we can apply \zitaa{MR2988005}{\cprop{3.7}{1772}} to get \(\nul{F\rk{z}}=\nul{\alpha_F}\cap\nul{\beta_F}\cap\nul{\nu_F\rk{\R}}\) for all \(z\in\pip\), where \(\rk{\alpha_F,\beta_F,\nu_F}\) denotes the \tNpo{\(F\)}.
 Consequently, \(\nul{A}\subseteq\nul{E}=\nul{\alpha_F}\cap\nul{\beta_F}\cap\nul{\nu_F\rk{\R}}\subseteq\nul{\alpha_F}\cap\nul{\nu_F\rk{\R}}\) follows.
\eproof

\section{Some facts on Nevanlinna pairs}\label{S-IV}
 In this section, we state some results on certain pairs of matrix-valued functions meromorphic in \(\pip\).
 These pairs take over the role of the free parameters within the parametrization of the set of solutions to the matricial power moment problems.
 Before we recall the definition of this well-known class of so-called Nevanlinna pairs, we observe the following well-known fact:

\breml{AB53N}
 The matrix \(\Jimq\) given by 
\beql{JQ}
 \Jimq 
 \defeq
\Mat{
 \Oqq&-\iu\Iq \\ 
 \iu\Iq&\Oqq
} 
\eeq
 is a \taaa{2q}{2q}{signature} matrix, \tie{}, \(\Jimq^\ad=\Jimq\) and \(\Jimq^2=\Iu{2q}\) hold true.
 Moreover, \(\tmat{A \\ B}^\ad \rk{-\Jimq}\tmat{A \\ B}=2\rim\rk{B^\ad A}\) 
 for all \(A,B\in\Cqq\).
 In particular, the case \(B=\Iq\) is of interest.
\erem

 Let \(\dom\) be a a \tne{} open subset of the complex plane.
 As usual, we will call a subset \(\mD\) of \(\dom\) a discrete subset of \(\dom\) if \(\mD\) does not have an accumulation point in \(\dom\).

\bdefnl{def-nev-paar}
 Let \(\phi\) and \(\psi\) be \tqqa{matrix}-valued functions meromorphic in \(\pip\).
 The pair \(\copa{\phi}{\psi}\) is called \emph{\tqqa{Nevanlinna} pair in \(\pip\)} if there is a discrete subset \(\mD\) of \(\pip\) such that the following three conditions are fulfilled:
 \baeqi{0}
  \il{def-nev-paar.i} \(\phi\) and \(\psi\) are holomorphic in \(\pip \setminus  \mD\). 
  \il{def-nev-paar.ii} \(\rank\smat{\phi \rk{w}\\\psi\rk{w}} =q\) for each \(w\in\pip \setminus  \mD\).
  \il{def-nev-paar.iii} \(\smat{\phi \rk{w}\\\psi\rk{w}}^\ad \rk{-\Jimq} \smat{\phi \rk{w}\\\psi\rk{w}}\in\Cggq\) for each \(w\in\pip \setminus  \mD\).
 \eaeqi
 We denote the set of all \tqqa{Nevanlinna} pairs in \(\pip\) by \(\PRF{q}\).
\edefn
 
\breml{IG.b}
 \rrem{AB53N} shows that condition~\ref{def-nev-paar.iii} of \rdefn{def-nev-paar} is equivalent to:
\begin{aseqi}{2}
 \il{IG.iii'} \(\rim\rk{\ek{\psi\rk{w}}^\ad \phi\rk{w}}\in\Cggq\) for all \(w\in\pip\setminus\mD\).
\end{aseqi}
\erem

\breml{Ma09_1.11}
 Let \(\copa{\phi}{\psi}\in\PRF{q}\).
 For each \tqqa{matrix}-valued function \(g\) meromorphic in \(\pip\) such that the function \(\det g\) does not vanish identically, one can easily see that the pair \(\copa{\phi g}{\psi g}\) belongs to \(\PRF{q}\) as well.
 Pairs  \(\copa{\phi_1}{\psi_1},\copa{\phi_2}{\psi_2}\) belonging to \(\PRFq\) are said to be \emph{equivalent} if there exist a \tqqa{matrix}-valued function \(g\) meromorphic in \(\pip\) and a discrete subset \(\mD\) of \(\pip\) such that \(\phi_{1}\), \(\psi_{1}\), \(\phi_{2}\), \(\psi_{2}\), and \(g\) are holomorphic in \(\pip\setminus\mD\) and that \(\det g\rk{w} \neq0\) as well as \(\phi_2\rk{w}=\phi_1\rk{w}g\rk{w}\) and \(\psi_2\rk{w}=\psi_1\rk{w}g\rk{w}\) hold true for each \(w\in\pip\setminus\mD\).
 Indeed, it is readily checked that this relation defines an equivalence relation on \(\PRF{q}\).
 For each \(\copa{\phi}{\psi}\in\PRF{q}\), let  \(\cpcl{\phi}{\psi}\)  denote the equivalence class generated by \(\copa{\phi}{\psi}\).
 Furthermore, if \(\mM\) is a \tne{} subset of \(\PRF{q}\), then let \(\cpsetcl{\mM}\defeq\setaca{\cpcl{\phi}{\psi}}{\copa{\phi}{\psi}\in\mM}\).
\erem

 Now we want to study special subclasses of the class \(\PRF{q}\).
 For each linear subspace \(\cU\) of \(\Cq\), let \(\OPu{\cU}\) be the orthogonal projection matrix onto \(\cU\) (see \rremss{PM}{R.P}).

\bnotal{CbezA5}
 Let \(M\in\Cqp\).
 We denote by \(\PRFqa{M}\) the set of all pairs \(\copa{\phi}{\psi}\in\PRF{q}\) such that \(\OPu{\ran{M}}\phi=\phi\) is fulfilled.
\enota

 The construction of \rnota{CbezA5} will later be used to treat Problem~\mprobR{2n}{\lleq} by choosing \(M=\hp{2n}\), where \(\hp{2n}\) is given by \rdefn{CM2.1.64} below.

\breml{CbemA8}%
 Let \(M\in\Cqp \) and let \(\copa{\phi}{\psi}\in\PRFqa{M}\).
 In view of \rrem{Ma09_1.11} and \rnota{CbezA5}, then \(\copa{\eta}{\theta}\in\PRFqa{M}\) holds true for all \(\copa{\eta}{\theta}\in\cpcl{\phi}{\psi}\).
\erem

 The following lemma plays an important role in the proof of \rprop{L1204}, which is an essential step to a main result of this paper.

\bleml{CfolA15}
 Let \(M\in\CHq\), let \(\copa{\phi}{\psi}\in\PRFqa{M}\), and let \(P\defeq\OPu{\ran{M}}\) and \(Q\defeq\OPu{\nul{M}}\).
 Then there exists a pair \(\copa{\cpfP }{\cpfQ}\in\cpcl{\phi}{\psi}\) such that
\begin{align}\label{ST37N}
 P\cpfP &=\cpfP,&
 \cpfP P&=\cpfP,&
 P\cpfQ &=\cpfQ -Q,&
 &\text{and}&
 \cpfQ P&=\cpfQ -Q.
\end{align}
\elem
\bproof
 Let \(r\defeq\rank M\).
 First we consider the case \(r\geq 1\).
 Let \(u_1,u_2,\dotsc,u_r\) be an orthonormal basis of \(\ran{M}\) and let \(U\defeq\mat{u_1,u_2,\dotsc,u_r}\).
 Then, from \zitaa{FKKM20}{\cprop{4.12(b)}{14}} we can infer that there exists a pair \(\copa{\tilde\phi}{\tilde\psi}\in\PRF{r}\) such that \(\copa{U\tilde\phi U^\ad}{U\tilde\psi U^\ad+\OPu{\ran{M}^\oc}}\in\PRFq\) and \(\cpcl{U\tilde\phi U^\ad}{U\tilde\psi U^\ad+\OPu{\ran{M}^\oc}}=\cpcl{\phi}{\psi}\), \tie{}\ \(\copa{U\tilde\phi U^\ad}{U\tilde\psi U^\ad+\OPu{\ran{M}^\oc}}\in\cpcl{\phi}{\psi}\) hold true.
 Thus, setting \(\cpfP\defeq U\tilde\phi U^\ad\) and \(\cpfQ\defeq U\tilde\psi U^\ad+\OPu{\ran{M}^\oc}\), we have \(\copa{\cpfP}{\cpfQ}\in\cpcl{\phi}{\psi}\).
 \rrem{R.P} shows \(P^2=P\) and \(P^\ad=P\).
 Clearly, \(PU=U\).
 Hence, we get \(U^\ad=\rk{PU}^\ad=U^\ad P\).
 Consequently, \(P\cpfP=PU\tilde\phi U^\ad=\cpfP\) and \(\cpfP P=U\tilde\phi U^\ad P=\cpfP\).
 From the assumption \(M\in\CHq\) and \rrem{tsa2} we get \(\nul{M}=\nul{M^\ad}=\ran{M}^\oc\), implying \(Q=\OPu{\ran{M}^\oc}\).
 Hence, \(\cpfQ=U\tilde\psi U^\ad+Q\). 
 Using \rrem{A.R.0<P<1}, we obtain furthermore \(Q=\Iq-P\).
 In view of \(P^2=P\), then \(QP=\Oqq\) and \(PQ=\Oqq\) follow.
 Consequently, \(P\cpfQ=P\rk{U\tilde\psi U^\ad+Q}=PU\tilde\psi U^\ad+PQ=U\tilde\psi U^\ad=\cpfQ-Q\) and \(\cpfQ P=\rk{U\tilde\psi  U^\ad+Q}P=U\tilde\psi U^\ad P+QP=U\tilde\psi U^\ad=\cpfQ-Q\).
 Thus, in the case \(r\geq1\), the proof is complete.
 If \(r=0\), then the assertion can be easily checked using \zitaa{FKKM20}{\cprop{4.12(a)}{14}}.
 We omit the details.
\eproof

\section{\hHp{s}}\label{S1111}
 Let \(\kappa\in\NOinf\) and let \(\seqska\) be a sequence of complex \tpqa{matrices}.
 For every choice of integers \(\ell \) and \(m\) fulfilling \(0\leq \ell \leq m \leq\kappa\), let
\begin{align}\label{yz}
 \yuu{\ell }{m}&\defeq\Mat{\su{\ell }\\\vdots\\\su{m}}&
&\text{and}&
 \zuu{\ell }{m}&\defeq\mat{\su{\ell },\dotsc,\su{m}}.
\end{align}
 Let
\begin{align}\label{Trip} 
 \Trip{0}&\defeq\Opq&
&\text{and}&
 \Trip{n}\defeq\zuu{n}{2n-1}\Hu{n-1}^\mpi\yuu{n}{2n-1}
\end{align}
 for each \(n\in\N\) such that \(2n-1\leq\kappa\).
 For all \(n\in\NO\) fulfilling \(2n+1\leq\kappa\), we also introduce the \tbHm{} \(\Kn\defeq\matauuuo{\su{j+k+1}}{j}{k}{0}{n}\). 
 For every choice of \(n\in\N\) fulfilling \(2n-1\leq\kappa\), we set
\[
 \Sigma_n
 \defeq\zuu{n}{2n-1}\Hu{n-1}^\mpi\Ku{n-1}\Hu{n-1}^\mpi\yuu{n}{2n-1}.
\]
 For each \(n\in\N\) fulfilling \(2n\leq\kappa\), let
\begin{align*}
 M_n&\defeq\zuu{n}{2n-1}\Hu{n-1}^\mpi\yuu{n+1}{2n}&
&\text{and}&
 N_n&\defeq\zuu{n+1}{2n}\Hu{n-1}^\mpi\yuu{n}{2n-1}.
\end{align*}
 Let 
\begin{align}\label{Lam}
 \Lam{0}&\defeq\Opq&
&\text{and}&
 \Lam{n}&\defeq M_n+N_n -\Sigma_n
\end{align}
 for all \(n\in\N\) fulfilling \(2n\leq\kappa\).
 Now we turn our attention to sequences of complex \tqqa{matrices} which are introduced in \rsec{S-II} and which are defined by certain properties of the \tbHms{} built from the given sequence.

\breml{R1624}
 Let \(\kappa\in\NOinf\) and let \(\seqska\in\Hggeqka\).
 Then \(\su{j}^\ad=\su{j}\) for each \(j\in\mn{0}{\kappa}\).
\erem

\breml{Schr2.5}
 Let \(\kappa\in\NOinf\) and let  \(\seqs{\kappa}\) be a sequence of complex \tqqa{matrices}.
 It is easy to see that \(\seqs{\kappa}\in\Hggeq{\kappa}\) is valid if and only if \(\seqs{m}\in\Hggeq{m}\) holds true for each \(m\in\mn{0}{\kappa}\).
\erem

 Now we recall useful parameters which will play a key role.
 
\bdefnnl{\zitaa{MR3014199}{\cdefn{2.3}{}}, \zitaa{FKKM20}{\cdefn{5.5}{}}}{CM2.1.64}
 Let \(\kappa\in\NOinf\) and let \(\seqska\) be a sequence of complex \tpqa{matrices}.
 For each \(k\in\NO\) fulfilling \(2k\leq\kappa\), let \(\hp{2k}\defeq\su{2k}-\Trip{k}\), where \(\Trip{k}\) is given by \eqref{Trip}, and, for each \(k\in\NO\) fulfilling \(2k+1\leq\kappa\), let \(\hp{2k+1}\defeq\su{2k+1}-\Lam{k}\), where \(\Lam{k}\) is given by \eqref{Lam}.
 Then \(\sHp{\kappa}\) is called the \emph{\tsHp{} (or sequence of canonical Hankel parameters) of \(\seqska\)}.
\edefn

 In view of \eqref{Trip} and \eqref{Lam}, we have in particular
\begin{align}\label{Hp.01}
 \hp{0}&=\su{0}&
&\text{and}&
 \hp{1}&=\su{1}.
\end{align}

\breml{CM2.1.65}
 Let \(\kappa\in\NOinf\), let \(\seqska\) be a sequence of complex \tpqa{matrices}, and let \(\sHp{\kappa}\) be the \tsHpo{\(\seqska\)}.
 Then, for all \(m\in\mn{0}{\kappa}\), the sequence \(\sHp{m}\) coincides with the \tsHpo{\(\seqs{m}\)}.
\erem

\bdefnl{D.nat-ext}
 Let \(n\in\NO\) and let \(\seqs{2n}\) be a sequence of complex \tpqa{matrices}.
 Let \(\hat{s}_j\defeq\su{j}\) for each \(j\in\mn{0}{2n}\) and let \(\hat{s}_{2n+1}\defeq\Lam{n}\), where \(\Lam{n}\) is given by \eqref{Lam}.
 Then we call the sequence \(\seqsh{2n+1}\) the \emph{\tnatexto{\(\seqs{2n}\)}}.
\edefn

 In view of the construction of \(\seqsh{2n+1}\), we can immediately see from \eqref{Lam} that this sequence is completely determined by the original sequence \(\seqs{2n}\).
 This observation will play an important role in the rest of the paper.
 It should be mentioned that the notion introduced in \rdefn{D.nat-ext} is in some sense analogous to the central extension of a matricial Carath\'eodory or Schur sequence.
 We get the following remark, which is essential for our further considerations.

\breml{R26-SK1}
 Let \(\kappa\in\NOinf\), let \(\seqska\) be a sequence of complex \tpqa{matrices}, and let \(\sHp{\kappa}\) be the \tsHpo{\(\seqska\)}.
 Let \(n\in\NO\) be such that \(2n\leq\kappa\) and let \(\shHp{2n+1}\) be the \tsHp{} of the \tnatext{} \(\seqsh{2n+1}\) of \(\seqs{2n}\).
 Then \rrem{CM2.1.65} and \rdefnss{CM2.1.64}{D.nat-ext} show that \(\hhp{j}=\hp{j}\) for all \(j\in\mn{0}{2n}\) and \(\hhp{2n+1}=\Opq\).
\erem

 In the sequel, we will repeatedly use the sequence \(\seqsh{2n+1}\) given in \rdefn{D.nat-ext}.

\breml{R8Z}
 Let \(\kappa\in\NOinf\), let \(\seqska\) be a sequence of complex \tpqa{matrices}, and let \(\sHp{\kappa}\) be the \tsHpo{\(\seqska\)}.
 Then it is readily checked that \(\seq{\hp{j}^\ad}{j}{0}{\kappa}\) is the \tsHpo{\(\seqa{s^\ad}{\kappa}\)}.
\erem

\bpropnl{\tcf{}~\zitaa{MR2805417}{\cprop{2.10(b)}{454} and \cprop{2.15(b)}{457}}}{DFKMT2.30N}%
 Let \(\kappa\in\NOinf\) and let \(\seqska\) be a sequence of complex \tqqa{matrices} with \tsHp{} \(\sHp{\kappa}\).
 Then \(\seqska\) belongs to \(\Hggeq{\kappa}\) if and only if the following three conditions are fulfilled:
\bAeqi{0}
 \il{DFKMT2.30N.I} \(\hp{2k}\in\Cggq\) for all \(k\in\NO\) such that \(2k\leq\kappa\).
 \il{DFKMT2.30N.II} If \(\kappa\geq1\), then \(\hp{2k-1}^\ad=\hp{2k-1}\) as well as \(\ran{\hp{2k-1}}\subseteq\ran{\hp{2k-2}}\) hold true for all \(k\in\N\) fulfilling \(2k-1\leq\kappa\).
 \il{DFKMT2.30N.III} If \(\kappa\geq2\), then \(\ran{\hp{2k}}\subseteq\ran{\hp{2k-2}}\) for all \(k\in\N\) such that \(2k\leq\kappa\).
\eAeqi
\eprop

\bpropnl{\tcf{}~\zitaa{MR2805417}{\cprop{2.10(d)}{454} and \cprop{2.15(c)}{457}}}{P1644}
 Let \(\kappa\in\NOinf\) and let \(\seqs{2\kappa}\) be a sequence of complex \tqqa{matrices} with \tsHp{} \(\sHp{2\kappa}\).
 Then \(\seqs{2\kappa}\) belongs to \(\Hgq{2\kappa}\) if and only if the following two conditions are fulfilled:
\bAeqi{0}
 \il{P1644.I} \(\hp{2k}\in\Cgq\) for all \(k\in\mn{0}{\kappa}\).
 \il{P1644.II} If \(\kappa\geq1\), then \(\hp{2k-1}^\ad=\hp{2k-1}\) for all \(k\in\mn{1}{\kappa}\).
\eAeqi
\eprop

\bpropl{L1237-A}
 Let \(\kappa\in\NOinf\) and let \(\seqska\in\Hggeqka\).
 Let \(n\in\NO\) be such that \(2n\leq\kappa\).
 Then the \tnatext{} \(\seqsh{2n+1}\) of \(\seqs{2n}\) belongs to \(\Hggeq{2n+1}\).
\eprop
\bproof
 Since \rrem{Schr2.5} yields \(\seqs{2n}\in\Hggeq{2n}\), the assertion follows from \rrem{R1624}, \zitaa{MR2805417}{\crem{2.1}{}}, and \zitaa{MR2570113}{\cprop{2.18}{770}}.
\eproof

 Now we want to recall the definition of first \tScht{} of sequences from \(\Cpq\).
 This requires a little preparation.
 Let \(\kappa\in\NOinf\) and let \(\seqska\) be a sequence of complex \tpqa{matrices}.
 Then the sequence \(\seqa{s^\rez}{\kappa}\) of complex \tqpa{matrices} given by \(\su{0}^\rez\defeq\su{0}^\mpi\) and, for all \(k\in\mn{1}{\kappa}\), recursively by
\[%
 \su{k}^\rez
 \defeq-\su{0}^\mpi\sum_{j=0}^{k-1}\su{k-j}\su{j}^\rez,
\]
 is called the \emph{\trso{\(\seqska\)}}.
 A detailed discussion of \trs{s} is given in \cite{MR3014197}.
 As our main application of \trs{s} we explain the elementary step of the Schur type algorithm under consideration.
 Let \(\kappa\in\minf{2}\cup\set{\infi}\) and let \(\seqska\) be a sequence of complex \tpqa{matrices} with \trs{} \(\seqa{s^\rez}{\kappa}\).
 Then the sequence \(\seqa{s^\St{1}}{\kappa-2}\) defined, for all \(j\in\mn{0}{\kappa-2}\), by
\[%
 \su{j}^\St{1}
 \defeq-\su{0}\su{j+2}^\rez\su{0}
\]
 is said to be the \emph{first \tSchto{\(\seqska\)}}.
 As considered in \zitaa{MR3014199}{\cdefn{9.1}{167}} already, the repeated application of the first \tScht{} in a natural way generates a corresponding algorithm for (finite or infinite) sequences of complex \tpqa{matrices}:

\breml{K-TRA}
 Let \(\kappa\in\NOinf\) and let \(\seqska\) be a sequence of complex \tpqa{matrices}.
 Then the sequence \(\seqa{s^\St{0}}{\kappa}\) given by \(\su{j}^\St{0}\defeq\su{j}\) for all \(j\in\mn{0}{\kappa}\), is called the \emph{\ath{0} \tSchto{\(\seqska\)}}.
 If \(\kappa\geq2\), then the \ath{k} \tScht{} is defined recursively:
 For all \(k\in\N\) fulfilling \(2k\leq\kappa\), the first \tScht{} \(\seqa{s^\St{k}}{\kappa-2k}\) of \(\seqa{s^\St{k-1}}{\kappa-2(k-1)}\) is called the \emph{\ath{k} \tSchto{\(\seqska\)}}.
\erem

\section{Special matrix polynomials}\label{S-VII}
 In this section, we discuss special matrix polynomials, which have been used, \teg{}, in \cite[\cform{(4.13)}{226}]{MR1624548} and \cite[\capp{C}{}]{MR3380267} already.
 Such matrix polynomials will be applied for the description of the solution set \(\RFOqskg{2n}\) of Problem~\rprobR{2n}{\lleq}, which was recognized as an equivalent reformulation of the Hamburger moment problem \mprobR{2n}{\lleq}.
 More precisely, these matrix polynomials act as generating matrix-valued functions of the linear fractional transformations which establish the description of the set \(\RFOqskg{2n}\).

\bnotal{BezWV}
 Let \(A,B\in\Cpq\).
 Then let \(V_{A,B}\colon\C\to\Coo{\rk{p+q}}{\rk{p+q}}\) be defined by
\[
 V_{A,B}\rk{z}
 \defeq
\begin{pmat}[{|}]
 \Opp  & -A\cr\- 
 A^\mpi  & z\Iq-A^\mpi B\cr
\end{pmat}.
\]
 We set \(V_A\defeq V_{A, \Opq}\) (see also \cite[\cch{4}{219}]{MR1624548}).
\enota

\bnotanl{see {\cite[\cpage{267}]{MR3380267}}}{SN-Tra}
 Let \(\kappa\in\NOinf\) and let \(\seqska\) be a sequence of complex \tpqa{matrices}.
 For all \(k\in\NO\) such that \(2k\leq\kappa\), let \(\seqa{s^\St{k}}{\kappa-2k}\) be the \ath{k} \tSchto{\(\seqska\)}.
 In view of \rnota{BezWV}, for all \(n\in\NO\) fulfilling \(2n\leq\kappa\), let
\beql{SN-Tra.0}%
 \fV^{\rk{\seqs{2n}}}
 \defeq
\begin{cases}
 V_{\su{0}^\St{0}}\tincase{n =0}\\ 
 V_{\su{0}^\St{0},\su{1}^\St{0}}V_{\su{0}^\St{1},\su{1}^\St{1}}\dotsm V_{\su{0}^\St{n-1},\su{1}^\St{n-1}}V_{\su{0}^\St{n}}\tincase{n\geq1}
\end{cases} 
\eeq
 and, for all \(n\in\NO\) fulfilling \(2n+1\leq\kappa\), let
\beql{SN-Tra.1}%
 \fV^{\rk{\seqs{2n+1}}}
 \defeq V_{\su{0}^\St{0},\su{1}^\St{0}}V_{\su{0}^\St{1},\su{1}^\St{1}}\dotsm V_{\su{0}^\St{n},\su{1}^\St{n}}.
\eeq
\enota

 The matrix-valued functions introduced in \rnota{SN-Tra} can be expressed by the \tsHp{} if the given sequence \(\seqska\) belongs to \(\Hggeqka\).
 This is caused by the following theorem:

\bthmnl{\zitaa{MR3014199}{\cthm{9.15}{178}}}{HS}
 Let \(\kappa\in\NOinf\) and let \(\seqska\in\Hggeq{\kappa}\) with \tsHp{} \(\sHp{\kappa}\).
 For each \(k\in\NO\) such that \(2k\leq\kappa\), let \(\seqa{s^\St{k}}{\kappa-2k}\) be the \ath{k} \tSchto{\(\seqska\)}.
 Then \(\hp{2k}=\su{0}^\St{k}\) for all \(k\in\NO\) such that \(2k\leq\kappa\) and \(\hp{2k+1}=\su{1}^\St{k}\) for all \(k\in\NO\) fulfilling \(2k+1\leq\kappa\).
\ethm

 Given \(n\in\N\) and arbitrary rectangular complex matrices \(A_1,A_2,\dotsc,A_n\), we write \(\col\seq{A_j}{j}{1}{n}\) (\tresp{}, \(\row\seq{A_j}{j}{1}{n}\)) for the block column (\tresp{}, block row) built from the matrices \(A_1,A_2,\dotsc,A_n\) if their numbers of columns (\tresp{}, rows) are all equal.

\bnotal{B.N.deg}
 Let \(P\) be a complex \tpqa{matrix} polynomial.
 For each \(n\in\NO\), let \(
 \cvuo{n}{P}
 \defeq\col\seq{A_j}{j}{0}{n}
\), where \((A_j)_{j=0}^\infi\) is the uniquely determined sequence of complex \tpqa{matrices}, such that \(P(w)=\sum_{j=0}^\infi w^jA_j\) holds true for all \(w\in\C\).
 Denote by \(\deg P\defeq\sup\setaca{j\in\NO}{A_j\neq\Opq}\) the \emph{degree of \(P\)}.
 If \(k\defeq\deg P\geq0\), then the matrix \(A_k\) is called the \emph{leading coefficient matrix of \(P\)}.
\enota

 In particular, we have \(\deg P=-\infty\), if \(P(z)=\Opq\) for all \(z\in\C\).

\breml{B.R.P=euY}
 If \(P\) is a complex \tqqa{matrix} polynomial, then \(P=\eu{n}\cvuo{n}{P}\) for all \(n\in\NO\) with \(n\geq \deg P\), where \(\eu{n}\colon\C\to\Coo{q}{(n+1)q}\) is defined by
\beql{NMB}
 \eua{n}{z}
 \defeq\mat{z^0\Iq,z^1\Iq,z^2\Iq,\dotsc,z^n\Iq}.
\eeq
\erem

\bdefnl{143.D1419}
 Let \(\kappa\in\NOinf\) and let \(\seqs{2\kappa}\) be a sequence of complex \tqqa{matrices}.
 A sequence \(\seq{P_k}{k}{0}{\kappa}\) of complex \tqqa{matrix} polynomials is called \emph{\tmosb{\(\seq{\su{j}}{j}{0}{2\kappa}\)}}, if it satisfies the following conditions:
\bAeqi{0}
 \il{143.D1419.I} For each \(k\in\mn{0}{\kappa}\), the matrix polynomial \(P_k\) has degree \(k\) and leading coefficient matrix \(\Iq\).
 \il{143.D1419.II} \(\ek{\cvuo{n}{P_j}}^\ad\Hu{n}\ek{\cvuo{n}{P_k}}=\Oqq\) for all \(j,k\in\mn{0}{\kappa}\) with \(j\neq k\), where \(n\defeq\max\set{j,k}\) and \(\Hu{n}\) is given by \eqref{N2}.
\eAeqi
\edefn
  
 Using \rrem{B.R.P=euY}, we can conclude from~\zitaa{MR2805417}{\cpropss{5.8(a1)}{479}{5.9(a)}{481}}: 
 
\bpropnl{\zitaa{MR4181333}{\cprop{D.5}{488}}}{B.P.oMP-lGs}
 Let \(\kappa\in\NOinf\), let \(\seqs{2\kappa}\in\Hggq{2\kappa}\), and let \(\seq{P_k}{k}{0}{\kappa}\) be a sequence of complex \tqqa{matrix} polynomials, satisfying condition~\ref{143.D1419.I} of \rdefn{143.D1419}.
 Then \(\seq{P_k}{k}{0}{\kappa}\) is a \tmosb{\(\seqs{2\kappa}\)} if and only if \(\Hu{k-1}X_{k}=\yuu{k}{2k-1}\) for all \(k\in\mn{1}{\kappa}\), where \(X_{k}\) is taken from the \tbr{} \(\cvuo{k}{P_k}=\tmat{-X_{k}\\ \Iq}\) and where \(\Hu{k-1}\) and \(\yuu{k}{2k-1}\) are given by \eqref{N2} and \eqref{yz}, respectively.
\eprop

 The following considerations of this section are aimed at having closer look at the four \tqqa{blocks} of the \taaa{2q}{2q}{matrix}-valued function \(\fV^{\rk{\seqs{2n}}}\).
 We will see that these are particular \tqqa{matrix} polynomials which satisfy recurrence formulas with the \tsHp{} as coefficients.
 For each \(\kappa\in\NOinf\), let
\beql{K5.11}
 \ev{\kappa}
 \defeq\sup\setaca{k\in\NO}{2k-1\leq\kappa}.
\eeq

\bdefnl{K19}
 Let \(\kappa\in\NOinf\), let \(\seqska\) be a sequence of complex \tqqa{matrices}, and let \(\sHp{\kappa}\) be the \tsHpo{\(\seqska\)}.
 Let \(\pa{0}, \pb{0}, \pc{0}, \pd{0}\colon\C\to\Cqq\) be defined by
\begin{align}\label{K19.0}
 \pa{0}\rk{z}&\defeq\Oqq,&
 \pb{0}\rk{z}&\defeq\Iq, &
 \pc{0}\rk{z}&\defeq\Oqq,&
&\text{and}&
 \pd{0}\rk{z}&\defeq\Iq.
\end{align}
 If \(\kappa\geq 1\), then let \(\pa{1}, \pb{1}, \pc{1}, \pd{1}\colon\C\to\Cqq\) be given via
\begin{align}\label{K19.1}%
 \pa{1}\rk{z}&\defeq\hp{0}, &
 \pb{1}\rk{z}&\defeq z\Iq-\hp{0}^\mpi \hp{1}, &
 \pc{1}\rk{z}&\defeq\hp{0}, &
&\text{and} &
 \pd{1}\rk{z}&\defeq z\Iq-\hp{1}\hp{0}^\mpi.
\end{align}
 If \(\kappa\geq 2\), then, for all \(k\in\mn{2}{\infi}\) fulfilling \(2k-1\leq\kappa\), let \(\pa{k}, \pb{k}, \pc{k}, \pd{k}\colon\C\to\Cqq\) be defined recursively by
\begin{align}
 \paa{k}{z}&\defeq\pa{k-1}\rk{z}(z\Iq-\hp{2k-2}^\mpi \hp{2k-1})-\pa{k-2}\rk{z}\hp{2k-4}^\mpi \hp{2k-2},\label{K19.pa}\\%\label{K8.3}
 \pba{k}{z}&\defeq\pb{k-1}\rk{z}(z\Iq-\hp{2k-2}^\mpi \hp{2k-1})-\pb{k-2}\rk{z}\hp{2k-4}^\mpi \hp{2k-2},\label{K19.pb}\\%\label{K8.4}
 \pca{k}{z}&\defeq\rk{z\Iq-\hp{2k-1}\hp{2k-2}^\mpi}\pc{k-1}\rk{z}-\hp{2k-2}\hp{2k-4}^\mpi \pc{k-2}\rk{z},\label{K19.pc}%
\intertext{and}
 \pda{k}{z}&\defeq\rk{z\Iq-\hp{2k-1}\hp{2k-2}^\mpi}\pd{k-1}\rk{z}-\hp{2k-2}\hp{2k-4}^\mpi \pd{k-2}\rk{z}.\label{K19.pd}%
\end{align}
 Regarding \eqref{K5.11}, we call the quadruple \(\abcd{\ev{\kappa}}\) the \emph{\tabcd{}}, abbreviating \sabcd{}, \emph{associated with \(\seqska\)}. 
\edefn

\breml{SK32R}
 Under the assumptions of \rdefn{K19}, for all \(k\in\NO\) such that \(2k-1\leq\kappa\), the matrix-valued functions \(\pa{k}\), \(\pb{k}\), \(\pc{k}\), and \(\pd{k}\) indeed are matrix polynomials, where the matrix polynomials \(\pb{k}\) and \(\pd{k}\) both have degree \(k\) and the same leading coefficient matrix \(\Iq\) (see also \zitaa{MR2805417}{\cthm{5.5}{475}}).
\erem

\breml{BK8.1}
 Let \(\kappa\in\NOinf\), let \(\seqska\) be a sequence of complex \tqqa{matrices}, and let \(\abcd{\ev{\kappa}}\) be the \sabcdo{\(\seqska\)}.
 In view of \rrem{CM2.1.65} and \rdefn{K19}, for all \(m\in\mn{0}{\kappa}\), then \(\abcd{\ev{m}}\) equals the \sabcdo{\(\seqs{m}\)}. 
\erem

\bpropnl{\tcf{}~\zitaa{MR2805417}{\cthm{5.5(a)}{475}}}{143.T1336}
 Let \(\kappa\in\NOinf\), let \(\seqs{2\kappa}\in\Hggq{2\kappa}\), and let \(\abcd{\kappa}\) be the \sabcdo{\(\seqs{2\kappa}\)}.
 Then \(\seq{\pb{k}}{k}{0}{\kappa}\) is a \tmosb{\(\seqs{2\kappa}\)}.
\eprop

 Let \(\kappa\in\NOinf\) and let \(\seqska\) be a sequence of complex \tpqa{matrices}.
 For all \(m\in\mn{0}{\kappa}\), then let
\beql{M.N.S}
 \SUu{m}
 \defeq
 \Mat{
   s_0    & s_1      & s_2    & \hdots    & s_m   \\
   \NM             & s_0      & s_1    & \hdots    & s_{m-1}   \\
   \NM             & \NM               & s_0    & \hdots    & s_{m-2}  \\
   \vdots        & \vdots          & \vdots        & \ddots    & \vdots  \\
   \NM             & \NM               & \NM             & \hdots    & s_0
 }.
\eeq

\bnotal{B.N.sec}
 Let \(\kappa\in\NOinf\), let \(\seqs{\kappa}\) be a sequence of complex \tqqa{matrices}, and let \(P\) be a complex \tqqa{matrix} polynomial with degree \(k\defeq\deg P\) satisfying \(k\leq\kappa+1\).
 Then let \(P^\secra{s}\colon\C\to\Cqq\) be defined by \(P^\secra{s}(z)=\Oqq\) if \(k\leq0\) and by
\(
 P^\secra{s}(z)
 \defeq\eua{k-1}{z}\mat{\Ouu{kq}{q},\SUu{k-1}}\cvuo{k}{P}
\)
 if \(k\geq1\).
\enota

\breml{143.R1532}
 Let \(\kappa\in\NOinf\), let \(\seqska\) be a sequence of complex \tqqa{matrices}, and let \(P\) and \(Q\) be two complex \tqqa{matrix} polynomials, each having degree at most \(\kappa+1\).
 Then \((P+Q)^\secra{s}=P^\secra{s}+Q^\secra{s}\).
 Furthermore, \((PA)^\secra{s}=P^\secra{s}A\) for all \(A\in\Cqq\).
\erem

\blemnl{\zitaa{MR4181333}{\clem{E.4}{489}}}{143.L1644}
 Let \(\kappa\in\Ninf\), let \(\seqska\) be a sequence of complex \tqqa{matrices}, let \(k\in\N\) with \(2k-1\leq\kappa\), and let \(P\) be a complex \tqqa{matrix} polynomial with degree \(k\) and leading coefficient matrix \(\Iq\), satisfying \(\Hu{k-1}X=\yuu{k}{2k-1}\), where the matrix \(X\) is taken from the \tbr{} \(\cvuo{k}{P}=\tmat{-X\\ \Iq}\).
 Let the matrix polynomial \(Q\) be defined by \(Q(w)\defeq wP(w)\).
 For all \(z\in\C\), then \(Q^\secra{s}(z)=zP^\secra{s}(z)\).
\elem

 Using \rnota{B.N.sec} we have:

\bpropl{P1203}
 Let \(\kappa\in\NOinf\), let \(\seqska\in\Hggeq{\kappa}\), and let \(\abcd{\ev{\kappa}}\) be the \sabcdo{\(\seqska\)}.
 For all \(k\in\NO\) with \(2k-1\leq\kappa\), then \(\pa{k}=\pb{k}^\secra{s}\).
\eprop
\bproof
 For all \(\ell\in\NO\) with \(2\ell-1\leq\kappa\), from \rrem{SK32R} we see that \(\pb{\ell}\) is a matrix polynomial with \(\deg\pb{\ell}=\ell\leq\kappa+1\) and leading coefficient matrix \(\Iq\).
 Because of \(\deg\pb{0}=0\) and \rnota{B.N.sec}, we have \(\pb{0}^\secra{s}(z)=\Oqq\) for all \(z\in\C\).
 In view of \eqref{K19.0}, hence \(\pb{0}^\secra{s}=\pa{0}\).
 In the case \(\kappa=0\), the proof is complete.
 Now suppose \(\kappa\geq1\).
 Using \(\deg\pb{1}=1\), \rnotass{B.N.sec}{B.N.deg}, \eqref{NMB}, \eqref{M.N.S}, and \eqref{Hp.01}, we obtain
\[
 \pb{1}^\secra{s}(z)
 =\eua{0}{z}\mat{\Oqq,\SUu{0}}\cvuo{1}{\pb{1}}
 =\Iq\cdot\mat{\Oqq,\su{0}}\matp{\uk}{\Iq}
 =\su{0}
 =\hp{0}
\]
 for all \(z\in\C\).
 In view of \eqref{K19.1}, hence \(\pb{1}^\secra{s}=\pa{1}\).
 In the case \(1\leq\kappa\leq2\) the proof is complete.
 Now suppose \(\kappa\geq3\).
 Then there exists an integer \(\ell \in\minf{2}\) with \(2\ell -1\leq\kappa\) such that \(\pa{\ell -1}=\pb{\ell -1}^\secra{s}\) and \(\pa{\ell -2}=\pb{\ell -2}^\secra{s}\) hold true.
 From \(\seqska\in\Hggeq{\kappa}\) we can infer \(\seqs{2\ell -2}\in\Hggq{2\ell -2}\).
 Regarding \rrem{BK8.1}, we can thus apply \rprop{143.T1336} to see that \(\seq{\pb{k}}{k}{0}{\ell -1}\) is a \tmosb{\(\seqs{2\ell -2}\)}.
 In view of \rprop{B.P.oMP-lGs}, then \rlem{143.L1644} yields \(B_{\ell -1}^\secra{s}(z)=z\pb{\ell -1}^\secra{s}(z)\) for all \(z\in\C\), where \(B_{\ell -1}\colon\C\to\Cqq\) is defined by \(B_{\ell -1}(z)\defeq z\pba{\ell -1}{z}\).
 Because of \(\pa{\ell -1}=\pb{\ell -1}^\secra{s}\), we hence obtain \(B_{\ell -1}^\secra{s}=A_{\ell -1}\), where \(A_{\ell -1}\colon\C\to\Cqq\) is defined by \(A_{\ell -1}(z)\defeq z\paa{\ell -1}{z}\).
 According to \eqref{K19.pb} and \eqref{K19.pa}, we have furthermore \(\pb{\ell }=B_{\ell -1}-\pb{\ell -1}\hp{2\ell -2}^\mpi\hp{2\ell -1}-\pb{\ell -2}\hp{2\ell -4}^\mpi\hp{2\ell -2}\) and \(\pa{\ell }=A_{\ell -1}-\pa{\ell -1}\hp{2\ell -2}^\mpi\hp{2\ell -1}-\pa{\ell -2}\hp{2\ell -4}^\mpi\hp{2\ell -2}\). 
 Regarding \(\deg B_{\ell -1}=\ell \leq\kappa+1\), the application of \rrem{143.R1532} yields then
 \[\begin{split}
  \pb{\ell }^\secra{s}
  &=\rk{B_{\ell -1}-\pb{\ell -1}\hp{2\ell -2}^\mpi\hp{2\ell -1}-\pb{\ell -2}\hp{2\ell -4}^\mpi\hp{2\ell -2}}^\secra{s}\\
  &=B_{\ell -1}^\secra{s}-\pb{\ell -1}^\secra{s}\hp{2\ell -2}^\mpi\hp{2\ell -1}-\pb{\ell -2}^\secra{s}\hp{2\ell -4}^\mpi\hp{2\ell -2}\\
  &=A_{\ell -1}-\pa{\ell -1}\hp{2\ell -2}^\mpi\hp{2\ell -1}-\pa{\ell -2}\hp{2\ell -4}^\mpi\hp{2\ell -2}
  =\pa{\ell }.
 \end{split}\]
  Thus, \(\pa{k}=\pb{k}^\secra{s}\) is inductively proved for all \(k\in\NO\) with \(2k-1\leq\kappa\).
\eproof

 Now we write the recurrence formulas stated in \rdefn{K19} in an alternative form. %

\breml{MD51N}
 Let \(\kappa\in\Ninf\) and let \(\seqska\) be a sequence of complex \tqqa{matrices} with \tsHp{} \(\sHp{\kappa}\) and \sabcd{} \(\abcd{\ev{\kappa}}\).
 For every choice of \(k\in\N\) such that \(2k+1\leq\kappa\) and \(z\in\C\), then 
\begin{align*}
 \Mat{
 \paa{k}{z}&\pa{k+1}\rk{z}\\
 \pba{k}{z}&\pb{k+1}\rk{z}
 }
 &=
 \Mat{
 \pa{k-1}\rk{z}&\paa{k}{z}\\
 \pb{k-1}\rk{z}&\pba{k}{z}
 }
 \begin{pmat}[{|}]
 \Oqq&-\hp{2k-2}^\mpi \hp{2k}\cr\-
 \Iq &z\Iq -\hp{2k}^\mpi\hp{2k+1}\cr
 \end{pmat}
\intertext{and}
 \Mat{
 \pca{k}{z}&\pda{k}{z}\\
 \pc{k+1}\rk{z}&\pd{k+1}\rk{z}
 }
 &=
 \begin{pmat}[{|}]
  \Oqq&\Iq \cr\-
 -\hp{2k}\hp{2k-2}^\mpi&z\Iq -\hp{2k+1}\hp{2k}^\mpi\cr
 \end{pmat}
 \Mat{
 \pc{k-1}\rk{z}&\pd{k-1}\rk{z}\\
 \pca{k}{z}&\pda{k}{z}
 }.
\end{align*}
\erem

\breml{BK8.2}
 Let \(\kappa\in\NOinf\), let \(\seqska\) be a sequence of \tH{} complex \tqqa{matrices}, and let \(\abcd{\ev{\kappa}}\) be the \sabcdo{\(\seqska\)}.
 In view of \rremss{R8Z}{A.R.A++*}, then \(\pca{k}{z}=[\paa{k}{\ko{z}}]^\ad\) and \(\pda{k}{z}=[\pba{k}{\ko{z}}]^\ad\) hold true for every choice of \(z\in\C\) and \(k\in\NO\) fulfilling \(2k-1\leq\kappa\).
\erem
 
 Now we are going to consider a quadruple of \tqqa{matrix} polynomials which, as we will see, connects the \sabcdo{\(\seqs{2n}\)} introduced in \rdefn{K19} with the particular sequence \(\seqsh{2n+1}\) introduced in \rdefn{D.nat-ext}.

\bnotal{N-abcdO}
 Let \(\kappa\in\NOinf\) and let \(\seqska\) be a sequence of complex \tqqa{matrices} with \tsHp{} \(\sHp{\kappa}\) and \sabcd{} \(\abcd{\ev{\kappa}}\), where \(\ev{\kappa}\) is given in \eqref{K5.11}.
 Let \(\pao{1},\pbo{1},\pco{1},\pdo{1}\colon\C\to\Cqq\) be defined by
\begin{align}\label{abcdO-1}
  \pao{1}\rk{z}&\defeq\hp{0},&
  \pbo{1}\rk{z}&\defeq z\Iq,&
  \pco{1}\rk{z}&\defeq\hp{0},&
  \pdo{1}\rk{z}&\defeq z\Iq.
\end{align}
 For all \(k\in\mn{2}{\infi}\) fulfilling \(2k-2\leq\kappa\), let \(\pao{k},\pbo{k},\pco{k},\pdo{k}\colon\C\to\Cqq\) be defined by
\begin{align*}
  \pao{k}\rk{z}&\defeq z\paa{k-1}{z}-\pa{k-2}\rk{z}\hp{2k-4}^\mpi \hp{2k-2},&
  \pbo{k}\rk{z}&\defeq z\pba{k-1}{z}-\pb{k-2}\rk{z}\hp{2k-4}^\mpi \hp{2k-2},\\
  \pco{k}\rk{z}&\defeq z\pca{k-1}{z}-\hp{2k-2}\hp{2k-4}^\mpi \pc{k-2}\rk{z},&
  \pdo{k}\rk{z}&\defeq z\pda{k-1}{z}-\hp{2k-2}\hp{2k-4}^\mpi \pd{k-2}\rk{z}.
\end{align*}
\enota

 In view of \eqref{K19.0} and \eqref{K19.1}, we have in particular
\begin{align*}
 \pao{2}\rk{z}&=z\hp{0},&
&&
 \pbo{2}\rk{z}
 &=z^2\Iq-z\hp{0}^\mpi\hp{1}-\hp{0}^\mpi\hp{2},\\
 \pco{2}\rk{z}&=z\hp{0},&
&\text{and}&
 \pdo{2}\rk{z}
 &=z^2\Iq-z\hp{1}\hp{0}^\mpi-\hp{2}\hp{0}^\mpi
\end{align*}
 for all \(z\in\C\).
 Regarding \rdefn{K19} and \rnota{N-abcdO}, we see:

\breml{R1358}
 Let \(\kappa\in\Ninf\) and let \(\seqska\) be a sequence of complex \tqqa{matrices} with \tsHp{} \(\sHp{\kappa}\) and \sabcd{} \(\abcd{\ev{\kappa}}\).
 Then
\begin{align*}
 \paa{k}{z}&=\paoa{k}{z}-\paa{k-1}{z}\hp{2k-2}^\mpi\hp{2k-1},&
 \pba{k}{z}&=\pboa{k}{z}-\pba{k-1}{z}\hp{2k-2}^\mpi\hp{2k-1},\\
 \pca{k}{z}&=\pcoa{k}{z}-\hp{2k-1}\hp{2k-2}^\mpi\pca{k-1}{z},&
 \pda{k}{z}&=\pdoa{k}{z}-\hp{2k-1}\hp{2k-2}^\mpi\pda{k-1}{z}
\end{align*}
 for all \(k\in\N\) with \(2k-1\leq\kappa\) and all \(z\in\C\).
\erem

 The following result shows that the \tqqa{matrix} polynomials introduced in \rnota{N-abcdO} occur in the \sabcdo{\(\seqsh{2n+1}\)}.

\bleml{L1237-B}
 Let \(\kappa\in\NOinf\) and let \(\seqska\) be a sequence of complex \tqqa{matrices}.
 Let \(n\in\NO\) be such that \(2n\leq\kappa\) and let \(\seqsh{2n+1}\) be the \tnatexto{\(\seqs{2n}\)}.
 Let \(\abcd{\ev{\kappa}}\) be the \sabcdo{\(\seqska\)} and let \(\habcd{n+1}\) be the \sabcdo{\(\seqsh{2n+1}\)}.
 Then
\begin{align*}%
 \pah{k}&=\pa{k}, &
 \pbh{k}&=\pb{k}, &
 \pch{k}&=\pc{k}, &
 &\text{and} &
 \pdh{k}&=\pd{k}
\end{align*}
 for each \(k\in\mn{0}{n}\).
 Furthermore,
\begin{align*}%
 \pah{n+1}&=\pao{n+1}, &
 \pbh{n+1}&=\pbo{n+1}, &
 \pch{n+1}&=\pco{n+1}, &
 &\text{and} &
 \pdh{n+1}&=\pdo{n+1}.
\end{align*}
\elem
\bproof
 Regarding \rrem{R26-SK1} and \rnota{N-abcdO}, this can be seen from the (recursive) construction in \rdefn{K19}.
\eproof

 In particular, we see from \rlem{L1237-B} that the quadruple \([\pao{n+1},\pbo{n+1},\pco{n+1},\pdo{n+1}]\) is completely determined by the sequence \(\seqs{2n}\).
 In the remaining considerations of this section, we always consider sequences of \tqqa{matrices} which are \tHnnde{}.
 This is done against to the background of \rthm{SK7}.
 
\bleml{BK8.6}
 Let \(\kappa\in\NOinf\), let \(\seqska\in\Hggeq{\kappa}\), and let \(\abcd{\ev{\kappa}}\) be the \sabcdo{\(\seqska\)}.
 Then \(\det\pba{k}{z}\neq 0\) as well as \(\det\pda{k}{z}\neq 0\) hold true for every choice of \(k\in\NO\) fulfilling \(2k-1\leq\kappa\) and for each \(z\in\C\setminus\R\).
\elem
\bproof
 The case \(\kappa=0\) is trivial.
 Consider the case that \(\kappa=2\tau\) with some \(\tau\in\Ninf\).
 Since \eqref{K5.11} yields \(\ev{\kappa}=\tau\) and \rrem{Schr2.7} provides \(\Hggeq{2\tau}\subseteq\Hggq{2\tau}\), then the assertion follows immediately from \zitaa{MR2805417}{\cthm{5.5}{475}}.
 Consider the case that \(\kappa=2n-1\) with some \(n\in\N\).
 Since we have supposed that \(\seqska\) belongs to \(\Hggeq{\kappa}\), there exists a matrix \(\su{2n}\in\Cqq\) such that \(\seqs{2n}\) belongs to \(\Hggq{2n}\).
 Let \(\dabcd{n}\) be the \sabcdo{\(\seqs{2n}\)}.
 Regarding \eqref{K5.11}, the application of \rrem{BK8.1} with \(m=2n-1\) to the sequence \(\seqs{2n}\) then yields that \(\dabcd{n}\) is the \sabcdo{\(\seqs{2n-1}\)}, \tie{}, with \(\seqska\).
 Thus, we receive that \(\pb{k}^\diamond=\pb{k}\) and \(\pd{k}^\diamond=\pd{k}\) for all \(k\in\mn{0}{n}\).
 Consequently, regarding \(\ev{\kappa}=n\), the application of \zitaa{MR2805417}{\cthm{5.5}{475}} to the sequence \(\seqs{2n}\) completes the proof.   
\eproof

\bleml{K15-1}
 Let \(\kappa\in\NOinf\) and let \(\seqska\in\Hggeq{\kappa}\).
 For every choice of \(n\in\NO\) fulfilling \(2n\leq\kappa\) and \(z\in\C\setminus\R\), then \(\det\pboa{n+1}{z}\neq0\) and \(\det\pdoa{n+1}{z}\neq0\).
\elem
\bproof
 Let \(n\in\NO\) be such that \(2n\leq\kappa\).
 \rprop{L1237-A} shows that the \tnatext{} \(\seqsh{2n+1}\) of \(\seqs{2n}\) belongs to \(\Hggeq{2n+1}\).
 Let \(\habcd{n+1}\) be the \sabcdo{\(\seqsh{2n+1}\)}.
 Using \rlem{L1237-B}, then we perceive that \(\pbh{n+1}=\pbo{n+1}\) and \(\pdh{n+1}=\pdo{n+1}\) hold true.
 Since the application of \rlem{BK8.6} to \(\seqsh{2n+1}\) yields \(\det\pbh{n+1}\rk{z}\neq 0\) and \(\det\pdh{n+1}\rk{z}\neq 0\) for all \(z\in\C\setminus\R\), the proof is complete.
\eproof

\breml{R1019}
 Let \(\kappa\in\NOinf\) and let \(\seqska\in\Hggeq{\kappa}\) with \tsHp{} \(\sHp{\kappa}\).
 From \rprop{DFKMT2.30N} we infer
\begin{align}\label{R1019.0}%
 \hp{j}^\ad&=\hp{j}&\text{for all }j&\in\mn{0}{\kappa}
\end{align}
 and \(\ran{\hp{2k-1}}\subseteq\ran{\hp{2k-2}}\) for all \(k\in\N\) fulfilling \(2k-1\leq\kappa\) as well as \(\ran{\hp{2k}}\subseteq\ran{\hp{2k-2}}\) for all \(k\in\N\) fulfilling \(2k\leq\kappa\).
 Using \eqref{R1019.0} and \rrem{tsa2}, we can conclude \(\nul{\hp{2k-2}}\subseteq\nul{\hp{2k-1}}\) for all \(k\in\N\) fulfilling \(2k-1\leq\kappa\) and \(\nul{\hp{2k-2}}\subseteq\nul{\hp{2k}}\) for all \(k\in\N\) such that \(2k\leq\kappa\).
 Hence, the application of \rrem{tsa12} yields
\begin{align}
 \hp{2k-2}\hp{2k-2}^\mpi\hp{2k-1}&=\hp{2k-1},&
 \hp{2k-1}\hp{2k-2}^\mpi\hp{2k-2}&=\hp{2k-1}&
 \text{for all }k&\in\N\text{ with }2k-1\leq\kappa\label{R1019.1}%
\intertext{and}
 \hp{2k-2}\hp{2k-2}^\mpi\hp{2k}&=\hp{2k},&
 \hp{2k}\hp{2k-2}^\mpi\hp{2k-2}&=\hp{2k}&
 \text{ for all }k&\in\N\text{ with }2k\leq\kappa.\label{R1019.2}%
\end{align}
\erem

 The following considerations yield interesting interrelations between the sequences forming a \tabcd{}.

\bleml{L0908}
 Let \(\kappa\in\Ninf\) and let \(\seqska\in\Hggeq{\kappa}\) with \tsHp{} \(\sHp{\kappa}\) and  \sabcd{} \(\abcd{\ev{\kappa}}\).
 For all \(n\in\NO\) such that  \(2n+1\leq\kappa\), then
\beql{L0908.R4A}
 \Mat{\pc{n}&\pd{n}\\\pc{n+1}&\pd{n+1}}
 \Mat{\Oqq&-\Iq \\\Iq &\Oqq}
 \Mat{\pa{n}&\pa{n+1}\\\pb{n}&\pb{n+1}}
 =
 \begin{pmat}[{|}]
  \pd{n}\pa{n}-\pc{n}\pb{n}&\pd{n}\pa{n+1}-\pc{n}\pb{n+1}\cr\-
  \pd{n+1}\pa{n}-\pc{n+1}\pb{n}&\pd{n+1}\pa{n+1}-\pc{n+1}\pb{n+1}\cr
 \end{pmat}
\eeq
 and
\beql{L0908.R4c}
 \Mat{\pc{n}&\pd{n}\\\pc{n+1}&\pd{n+1}}
 \Mat{\Oqq&-\Iq \\\Iq &\Oqq}
 \Mat{\pa{n}&\pa{n+1}\\\pb{n}&\pb{n+1}}
 =
 \Mat{
  \Oqq&\hp{2n}\\
  -\hp{2n}&\Oqq
  }.
\eeq
\elem
\bproof
 The identity \eqref{L0908.R4A} follows for all  \(n\in\NO\) with \(2n+1\leq\kappa\) by straightforward computation.
 According to \rrem{R1019}, we have \eqref{R1019.1} and \eqref{R1019.2}.
 Let \(\varepsilon\colon\C\to\C\) be defined by \(\varepsilon\rk{z}\defeq z\).
 From \eqref{L0908.R4A}, \rdefn{K19}, and \eqref{R1019.1} we get
\begin{multline}\label{L0908.R4BaaD}
 \Mat{\pc{0}&\pd{0}\\\pc{1}&\pd{1}}
 \Mat{\Oqq&-\Iq \\\Iq &\Oqq}
 \Mat{\pa{0}&\pa{1}\\\pb{0}&\pb{1}}
 =
 \begin{pmat}[{|}]
   \pd{0}\pa{0}-\pc{0}\pb{0}&\pd{0}\pa{1}-\pc{0}\pb{1}\cr\-
   \pd{1}\pa{0}-\pc{1}\pb{0}&\pd{1}\pa{1}-\pc{1}\pb{1}\cr
 \end{pmat}\\
 =
 \begin{pmat}[{|}]
   \Oqq&\hp{0}\cr\-
   -\hp{0}&\rk{\varepsilon \Iq -\hp{1}\hp{0}^\mpi}\hp{0}-\hp{0}\rk{\varepsilon \Iq -\hp{0}^\mpi\hp{1}}\cr
 \end{pmat}
 =\Mat{
  \Oqq&\hp{0}\\
  -\hp{0}&\Oqq
  }.
\end{multline}
 If \(\kappa\leq2\), then, regarding  \eqref{L0908.R4BaaD}, the proof is finished.
 
 Now suppose \(\kappa\geq3\).
 Regarding \eqref{L0908.R4BaaD}, the following statement holds true:
\bAeqi{0}
 \il{L0908.I} There is an \(m\in\N\) fulfilling \(2m+1\leq\kappa\) such that \eqref{L0908.R4c} holds true for all \(n\in\mn{0}{m-1}\).
\eAeqi
 Taking into account \rrem{MD51N}, \ref{L0908.I}, and \eqref{R1019.2}, we conclude
\beql{L0908.R4Baz}\begin{split}
  &\Mat{\pc{m}&\pd{m}\\\pc{m+1}&\pd{m+1}}
  \Mat{\Oqq&-\Iq \\\Iq &\Oqq}
  \Mat{\pa{m}&\pa{m+1}\\\pb{m}&\pb{m+1}}\\
  &=
  \begin{pmat}[{|}]
   \Oqq&\Iq \cr\-
   -\hp{2m}\hp{2m-2}^\mpi&\varepsilon \Iq -\hp{2m+1} \hp{2m}^\mpi\cr
  \end{pmat}\\
  &\qquad\times
  \Mat{\pc{m-1}&\pd{m-1}\\\pc{m}&\pd{m}}
  \Mat{\Oqq&-\Iq \\\Iq &\Oqq}
  \Mat{\pa{m-1}&\pa{m}\\\pb{m-1}&\pb{m}}
  \begin{pmat}[{|}]
   \Oqq&-\hp{2m-2}^\mpi \hp{2m}\cr\-
   \Iq &\varepsilon \Iq -\hp{2m}^\mpi \hp{2m+1}\cr
  \end{pmat}\\
  &=
  \begin{pmat}[{|}]
   \Oqq&\Iq \cr\-
   -\hp{2m}\hp{2m-2}^\mpi&\varepsilon \Iq -\hp{2m+1} \hp{2m}^\mpi\cr
  \end{pmat}
  \Mat{\Oqq&\hp{2m-2}\\-\hp{2m-2}&\Oqq}
  \begin{pmat}[{|}]
   \Oqq&-\hp{2m-2}^\mpi \hp{2m}\cr\-
   \Iq &\varepsilon \Iq -\hp{2m}^\mpi \hp{2m+1}\cr
  \end{pmat}\\
  &=\Mat{
  \Oqq&\hp{2m-2}\hp{2m-2}^\mpi \hp{2m}\\
  -\hp{2m}\hp{2m-2}^\mpi\hp{2m-2}&R_{m+1}}
  =\Mat{\Oqq& \hp{2m}\\- \hp{2m}&R_{m+1}},
\end{split}\eeq
 where
\[
 R_{m+1}
 \defeq-\hp{2m}\hp{2m-2}^\mpi\hp{2m-2}\rk{\varepsilon \Iq -\hp{2m}^\mpi \hp{2m+1}}+\rk{\varepsilon \Iq -\hp{2m+1} \hp{2m}^\mpi}\hp{2m-2}\hp{2m-2}^\mpi \hp{2m}.
\]
 Using \eqref{R1019.2} and \eqref{R1019.1}, one can easily see that
\beql{L0908.R4Bax}\begin{split}
 R_{m+1}
 &=-\hp{2m}\rk{\varepsilon \Iq -\hp{2m}^\mpi \hp{2m+1}}+\rk{\varepsilon \Iq -\hp{2m+1} \hp{2m}^\mpi}\hp{2m}\\
 &=-\varepsilon\hp{2m}+\hp{2m+1}+\varepsilon\hp{2m}-\hp{2m+1}
 =\Oqq.
\end{split}\eeq
 By virtue of \eqref{L0908.R4Baz}, \eqref{L0908.R4Bax}, and \ref{L0908.I}, we get that \eqref{L0908.R4c} is fulfilled for all  \(n\in\mn{0}{m}\).
 Therefore, \eqref{L0908.R4c} is proved by mathematical induction.
\eproof

\bleml{P8.14}
 Let \(\kappa\in\NOinf\) and let \(\seqska\in\Hggeq{\kappa}\) with \tsHp{} \(\sHp{\kappa}\) and \sabcd{} \(\abcd{\ev{\kappa}}\).
 Then \(\pda{n}{z}\paa{n}{z}=\pca{n}{z}\pba{n}{z}\) is valid for each \(n\in\NO\) fulfilling \(2n-1\leq\kappa\) and all \(z\in\C\).
 Furthermore, if \(\kappa\geq 1\), then \(\pda{n-1}{z}\paa{n}{z}-\pc{n-1}\rk{z}\pba{n}{z}=\hp{2n-2}\) and \(\pda{n}{z}\pa{n-1}\rk{z}-\pca{n}{z}\pba{n-1}{z}=-\hp{2n-2}\)
 hold true for all \(n\in\N\) fulfilling \(2n-1\leq\kappa\) and all \(z\in\C\).   
\elem
\bproof
 This is a consequence of \eqref{K19.0} and \rlem{L0908}.
\eproof

 Against the background of later application in \rprop{CH428} below, we consider now a sequence \(\seqs{2n+1}\in\Hggeq{2n+1}\).
 We are going to determine the \tqqa{blocks} of \(\fV^{\rk{\seqs{2n+1}}}\) in terms of the \tqqa{matrix} polynomials \(\pa{n},\pa{n+1},\pb{n},\pb{n+1}\) introduced in \rdefn{K19}.

\bpropl{C2A3}
 Let \(n\in\NO\) and let \(\seqs{2n+1}\in\Hggeq{2n+1}\) with \tsHp{} \(\sHp{2n+1}\) and \sabcd{} \(\abcd{n+1}\).
 Then the matrix-valued function \(\fV^{\rk{\seqs{2n+1}}}\colon\C\to\Coo{2q}{2q}\) given by \rnota{SN-Tra} admits, for each \(z\in\C\), the representation 
\[
 \fV^{\rk{\seqs{2n+1}}}\rk{z}
 =
 \Mat{
  -\paa{n}{z}\hp{2n}^\mpi  & -\paa{n+1}{z} \\ 
  \pba{n}{z}\hp{2n}^\mpi  & \pba{n+1}{z}
 }.
\]
\eprop
\bproof
 Because of assumption \(\seqs{2n+1}\in\Hggeq{2n+1}\) and \rthm{HS}, we have \(\hp{2k}=\su{0}^\St{k}\) and \(\hp{2k+1}=\su{1}^\St{k}\) for all \(k\in\mn{0}{n}\).
 In view of \eqref{SN-Tra.1}, \rnota{BezWV}, \eqref{K19.0}, and \eqref{K19.1}, then
\[
 \fV^{(\seqs{1})}\rk{z}
 =V_{\su{0}^\St{0},\su{1}^\St{0}}\rk{z}
 =V_{\hp{0},\hp{1}}\rk{z}
 =
 \begin{pmat}[{|}]
  \Oqq & -\hp{0} \cr\-
  \hp{0}^\mpi  & z\Iq-\hp{0}^\mpi \hp{1}\cr
  \end{pmat}
  =
  \Mat{
  -\pa{0}\rk{z}\hp{0}^\mpi  & -\pa{1}\rk{z} \cr\-
  \pb{0}\rk{z}\hp{0}^\mpi  & \pb{1}\rk{z}\cr
 }
\]
 follows for all \(z\in\C\).
 In the case \(n=0\), the proof is finished.
 
 Now suppose \(n\geq1\).
 Proceeding inductively, we assume that there is an \(m\in\mn{1}{n}\) such that
\[
 \fV^{(\seqs{2k-1})}\rk{z}
 =
 \Mat{
  -\pa{k-1}\rk{z}\hp{2(k-1)}^\mpi  & -\paa{k}{z} \\ 
  \pb{k-1}\rk{z}\hp{2(k-1)}^\mpi  & \pba{k}{z}
 }
\]
 is valid for every choice of \(k\in\mn{1}{m}\) and \(z\in\C\).
 Consequently, from \eqref{SN-Tra.1}, \rnota{BezWV}, \eqref{K19.pa}, and \eqref{K19.pb} we get then
\[\begin{split}
  \fV^{(\seqs{2m+1})}\rk{z}
  &=\fV^{(\seqs{2m-1})}\rk{z}V_{\su{0}^\St{m},\su{1}^\St{m}}\rk{z}
  =\fV^{(\seqs{2m-1})}\rk{z}V_{\hp{2m},\hp{2m+1}}\rk{z}\\
  &=
  \Mat{
  -\pa{m-1}\rk{z}\hp{2(m-1)}^\mpi  & -\pa{m}\rk{z}\\
  \pb{m-1}\rk{z}\hp{2(m-1)}^\mpi  & \pb{m}\rk{z}
  }
  \begin{pmat}[{|}]
  \Oqq & -\hp{2m} \cr\-
  \hp{2m}^\mpi  & z\Iq-\hp{2m}^\mpi \hp{2m+1}\cr
  \end{pmat}\\
  &=
  \begin{pmat}[{|}]
  -\pa{m}\rk{z}\hp{2m}^\mpi  & \pa{m-1}\rk{z}\hp{2(m-1)}^\mpi \hp{2m}-\pa{m}\rk{z}\rk{z\Iq-\hp{2m}^\mpi \hp{2m+1}}\cr\- 
  \pb{m}\rk{z}\hp{2m}^\mpi  &-\pb{m-1}\rk{z}\hp{2(m-1)}^\mpi \hp{2m}+\pb{m}\rk{z}\rk{z\Iq-\hp{2m}^\mpi \hp{2m+1}} \cr
  \end{pmat}\\
  &=
  \Mat{
  -\pa{m}\rk{z}\hp{2m}^\mpi  & -\pa{m+1}\rk{z} \\ 
  \pb{m}\rk{z}\hp{2m}^\mpi  & \pb{m+1}\rk{z}
  } 
\end{split}\]
 for all \(z\in\C\).
 The assertion is proved inductively.
\eproof

\bpropl{CH428}
 Let \(n\in\NO\) and let \(\seqs{2n}\in\Hggeq{2n}\) with \tsHp{} \(\sHp{2n}\) and \sabcd{} \(\abcd{n}\).
 Then the matrix-valued function \(\fV^{\rk{\seqs{2n}}}\colon\C\to\Coo{2q}{2q}\) given by \rnota{SN-Tra} admits the representation 
\beql{N101N}
 \fV^{\rk{\seqs{2n}}}\rk{z}
 =
 \Mat{
  -\paa{n}{z}\hp{2n}^\mpi  & -\paoa{n+1}{z} \\ 
  \pba{n}{z}\hp{2n}^\mpi  & \pboa{n+1}{z}
 }
\eeq
 for each \(z\in\C\), %
 where \(\pao{n+1}\) and \(\pbo{n+1}\) are defined in \rnota{N-abcdO}.
\eprop
\bproof
 \rprop{L1237-A} shows that the \tnatext{} \(\seqsh{2n+1}\) of \(\seqs{2n}\) belongs to \(\Hggeq{2n+1}\).
 Let \(\habcd{n+1}\) be the \sabcdo{\(\seqsh{2n+1}\)}.
 \rlem{L1237-B} yields then \(\pah{n}=\pa{n}\) and \(\pbh{n}=\pb{n}\) as well as \(\pah{n+1}=\pao{n+1}\) and \(\pbh{n+1}=\pbo{n+1}\).
 Let \(\shHp{2n+1}\) be the \tsHpo{\(\seqsh{2n+1}\)}.
 From \rrem{R26-SK1} then we recognize \(\hhp{2n}=\hp{2n}\) and \(\hhp{2n+1}=\Oqq\).
 Because of \rthm{HS}, thus \(\hat{s}_1^\St{n}=\Oqq\).
 According to \rnota{BezWV}, hence \(V_{\hat{s}_{0}^\St{n},\hat{s}_{1}^\St{n}}=V_{\hat{s}_{0}^\St{n}}\).
 Taking into account \eqref{SN-Tra.1} and \eqref{SN-Tra.0}, then we conclude \(\fV^{\rk{\seqsh{2n+1}}}=\fV^{\rk{\seqsh{2n}}}\).
 In view of \rdefn{D.nat-ext}, thus \(\fV^{\rk{\seqsh{2n+1}}}=\fV^{\rk{\seqs{2n}}}\) follows.
 Consequently, the application of \rprop{C2A3} to \(\seqsh{2n+1}\) completes the proof.
\eproof

 Now we are able to rewrite \zitaa{FKKM20}{\cthm{8.11}{43}} in terms of the \tabcd{}.
 
\bthml{CM123}
 Let \(n\in\NO\) and let \(\seqs{2n}\in\Hggeq{2n}\) with \tsHp{} \(\sHp{2n}\) and \sabcd{} \(\abcd{n}\).
 Further, let \(\pars{n}\), \(\pbrs{n}\), \(\paors{n+1}\), and \(\pbors{n+1}\) be the restrictions of \(\pa{n}\), \(\pb{n}\), \(\pao{n+1}\), and \(\pbo{n+1}\) onto \(\pip\), respectively, where \(\pao{n+1}\) and \(\pbo{n+1}\) are defined in \rnota{N-abcdO}.
 Then:
\benui
 \il{CM123.a} For each pair \(\copa{\phi}{\psi}\in\PRFqa{\hp{2n}}\), the function \(\det\rk{\pbrs{n}\hp{2n}^\mpi\phi+\pbors{n+1}\psi}\) does not vanish identically and the matrix-valued function \(F\defeq-\rk{\pars{n}\hp{2n}^\mpi\phi+\paors{n+1}\psi}\rk{\pbrs{n}\hp{2n}^\mpi\phi+\pbors{n+1}\psi}^\inv\) belongs to \(\RFOqskg{2n}\).
 \il{CM123.b} For each \(F\in\RFOqskg{2n}\), there exists a pair \(\copa{\phi}{\psi}\in\PRFqa{\hp{2n}}\) such that both \(\phi\) and \(\psi\) are holomorphic in \(\pip\), that %
\beql{L1204.1}%
 \det\rk*{\pba{n}{z}\hp{2n}^\mpi\phi\rk{z}+\pboa{n+1}{z}\psi\rk{z}}
 \neq0
\eeq
 holds true for all \(z\in\pip\), and that \(F\) admits the representation %
\beql{L1204.2}%
 F\rk{z}
 =-\ek*{\paa{n}{z}\hp{2n}^\mpi\phi\rk{z}+\paoa{n+1}{z}\psi\rk{z}}\ek*{\pba{n}{z}\hp{2n}^\mpi\phi\rk{z}+\pboa{n+1}{z}\psi\rk{z}}^\inv 
\eeq 
 for all \(z\in\pip\).
 \il{CM123.c} For each \(\ell\in\set{1,2}\) let \(\copa{\phi_\ell}{\psi_\ell}\in\PRFqa{\hp{2n}}\) and let \(F_\ell\defeq-\rk{\pars{n}\hp{2n}^\mpi\phi_\ell+\paors{n+1}\psi_\ell}\rk{\pbrs{n}\hp{2n}^\mpi\phi_\ell+\pbors{n+1}\psi_\ell}^\inv\).
 Then \(F_1=F_2\) if and only if \(\cpcl{\phi_1}{\psi_1}=\cpcl{\phi_2}{\psi_2}\).
\eenui
\ethm
\bproof
 Using \rthm{HS} and the notation therein, we receive \(\hp{2n}=\su{0}^\St{n}\).
 From \rprop{CH428} we obtain the \tqqa{\tbr{}} \eqref{N101N} for all \(z\in\C\).
 Applying \zitaa{FKKM20}{\cthm{8.11}{43}} and comparing the \tqqa{\tbr{}} therein with \eqref{N101N} completes the proof.
\eproof

 As already observed above, parametrizations of the set \(\RFOqskg{2n}\) with independent parameters are given in \cite{MR1624548,MR1395706,Thi06,FKKM20}.
 Now we get a refinement of a result (see \zitaa{MR3380267}{\cthm{12.1}{272}}), which concerns the moment problem~\mprobR{2n}{=}.
 
\bthml{T1153}
 Let \(n\in\NO\) and let \(\seqs{2n}\in\Hggeq{2n}\) with \tsHp{} \(\sHp{2n}\) and \sabcd{} \(\abcd{n}\).
 Further, let \(\pars{n}\), \(\pbrs{n}\), \(\paors{n+1}\), and \(\pbors{n+1}\) be the restrictions of \(\pa{n}\), \(\pb{n}\), \(\pao{n+1}\), and \(\pbo{n+1}\) onto \(\pip\), respectively, where \(\pao{n+1}\) and \(\pbo{n+1}\) are defined in \rnota{N-abcdO}.
 Then:
\benui
 \il{T1153.a} For each \(G\in\Pqevena{\hp{2n}}\), the function \(\det\rk{\pbrs{n}\hp{2n}^\mpi G+\pbors{n+1}}\) does not vanish identically and the matrix-valued function \(F\defeq-\rk{\pars{n}\hp{2n}^\mpi G+\paors{n+1}}\rk{\pbrs{n}\hp{2n}^\mpi G+\pbors{n+1}}^\inv\) belongs to \(\RFOqsg{2n}\).
 \il{T1153.b} For each \(F\in\RFOqsg{2n}\), there exists a unique \(G\in\Pqevena{\hp{2n}}\) such that \(\det\rk{\pba{n}{z}\hp{2n}^\mpi G\rk{z}+\pboa{n+1}{z}}
 \neq 0\) and \(F\rk{z}=-\ek{\paa{n}{z}\hp{2n}^\mpi G\rk{z}+\paoa{n+1}{z}}\ek{\pba{n}{z}\hp{2n}^\mpi G\rk{z}+\pboa{n+1}{z}}^\inv\) for all \(z\in\pip\).
\eenui
\ethm
\bproof
 Using \rthm{HS} and the notation therein, we receive \(\hp{2n}=\su{0}^\St{n}\).
 From \rprop{CH428} we obtain the \tqqa{\tbr{}} \eqref{N101N} for all \(z\in\C\).
 Applying \zitaa{MR3380267}{\cprop{11.22(a)}{270} and \cthm{12.1}{272}} completes the proof. 
\eproof

\section{Particular rational matrix-valued functions}\label{S1129}
 In \rsec{Cha13}, we will focus on the values of the respective functions belonging to the solution set \(\RFOqskg{2n}\) of Problem~\rprobR{2n}{\lleq}, \tie{}, we want to describe the set 
\begin{align}\label{N128N}
 \setaca*{F\rk{w}}{F\in\RFOqskg{2n}}
\end{align}
 for each \(\seqs{2n}\in\Hggeq{2n}\) and each \(w\in\pip\).
 In order to realize this aim, we are essentially guided by \rprop{CH428} and \rthm{CM123}.
 These statements contain essential information about the  description of the set \(\RFOqskg{2n}\).
 The concrete shape of the \tqqa{blocks} of the \taaa{2q}{2q}{matrix} polynomial \(\fV^{\rk{\seqs{2n}}}\) given by \eqref{N101N} suggest to consider a particular sequence of rational \tqqa{matrix}-valued functions.
 It will turn out that these functions play a key role in our further considerations.
 We start with a little preparation.

\breml{TC}
 Let \(\kappa\in\NOinf\), let \(\seqska\in\Hggeqka\) with \sabcd{} \(\abcd{\ev{\kappa}}\), and let \(n\in\NO\) be such that \(2n-1\leq\kappa\).
 According to \rrem{SK32R}, the functions \(\det\pb{n}\) and \(\det\pd{n}\) are non-trivial polynomials for which, in view of \eqref{NM} and \rlem{BK8.6}, the respective sets of zeros \(\nst{\det\pb{n}}\) and \(\nst{\det\pd{n}}\) are finite and, in particular, discrete subsets of \(\C\) such that \(\nst{\det\pb{n}}\cup\nst{\det\pd{n}}\subseteq\R\).
 Consequently, \(\pb{n}^\inv\) and \(\pd{n}^\inv\) are matrix-valued functions meromorphic in \(\C\), which fulfill \(\C\setminus\nst{\det\pb{n}}\subseteq\holpt{\pb{n}^\inv }\), \(\C\setminus\nst{\det\pd{n}}\subseteq\holpt{\pd{n}^\inv }\), and, in particular, \(\CR\subseteq\holpt{\pb{n}^\inv }\cap\holpt{\pd{n}^\inv }\). 
\erem

 The following sequence of rational matrix-valued functions will play a key role of our following considerations.
 We note that a procedure is used here which is applied in a similar way in \cite{MR2656833} for the matricial Carath\'eodory problem. 
 
\bdefnl{K17.1}
 Let \(\kappa\in\NOinf\) and let \(\seqska\in\Hggeqka\) with \tsHp{} \(\sHp{\kappa}\) and \sabcd{} \(\abcd{\ev{\kappa}}\).
 Let \(\XF{-1}\colon\C\to\Cqq\) be defined by \(\XF{-1}\rk{z}\defeq\Oqq\).
 For all \(n\in\NO\) such that \(2n\leq\kappa\), let \(\XF{2n}\defeq\hp{2n}\pb{n}^\inv\pbo{n+1}\).
 For all \(n\in\NO\) fulfilling \(2n+1\leq\kappa\), let \(\XF{2n+1}\defeq\hp{2n}\pb{n}^\inv\pb{n+1}\).
 Then \(\seqX{\kappa}\) is called the \emph{\tsXFo{\(\seqska\)}}.
\edefn

 In view of \eqref{K19.0}, \eqref{abcdO-1}, \eqref{Hp.01}, \eqref{K19.1}, and \eqref{R1019.1}, we have
\begin{align}\label{XF.01}
 \XFa{-1}{z}&=\Oqq,&
 \XFa{0}{z}&=z\hp{0}=z\su{0},&
&\text{and}&
 \XFa{1}{z}&=z\hp{0}-\hp{1}=z\su{0}-\su{1}
\end{align}
 for all \(z\in\C\).
 The \tXF{s} can be considered in a wider sense as analogue of the rational \tqqa{matrix} valued function \(\Theta_n\) defined in \zitaa{MR2222523}{\cpage{294}}
 
\breml{XF-tru}
 Let \(\kappa\in\NOinf\), let \(\seqska\in\Hggeqka\) with \tsXF{} \(\seqX{\kappa}\) and let \(m\in\mn{0}{\kappa}\).
 Regarding \rdefn{K17.1}, \rnota{N-abcdO}, and \eqref{K5.11}, in view of \rremsss{Schr2.5}{CM2.1.65}{BK8.1} then \(\seqs{m}\) belongs to \(\Hggeq{m}\) and \(\seqX{m}\) equals the \tsXFo{\(\seqs{m}\)}. 
\erem
 
 We introduce a notation which is operated throughout the attaching considerations.
 For each \(m\in\NO\), let
\begin{align}\label{SK5-11}
 \ef{m}
 &\defeq\max\setaca{k\in\NO}{2k\leq m}&
&\text{and}&
 \eff{m}
 &\defeq2\ef{m},
\end{align}
 \tie{}, if \(m=2n\) or \(m=2n+1\) for some \(n\in\NO\), then \(\ef{m}=n\) and \(\eff{m}=2n\).

\breml{R0735}
 Let \(\kappa\in\NOinf\) and let \(\seqska\in\Hggeqka\) with \tsXF{} \(\seqX{\kappa}\).
 In view of \rdefn{K17.1} and \rrem{TC}, then \(\CR\subseteq\C\setminus\nst{\det\pb{\ef{m}}}\subseteq\holpt{\XF{m}}\) for all \(m\in\mn{0}{\kappa}\).
\erem

\breml{R52SK}
 Let the assumptions of \rdefn{K17.1} be fulfilled.
 In view of \rrem{TC}, then:
\benui
 \il{R52SK.b} Let \(n\in\NO\) be such that \(2n\leq\kappa\).
 For all \(z\in\CR\), then \(\det\pba{n}{z}\neq0\) and \(\XFa{2n}{z}=\hp{2n}\ek{\pba{n}{z}}^\inv\pboa{n+1}{z}\) hold true.
 \il{R52SK.a} Let \(n\in\NO\) be such that \(2n+1\leq\kappa\).
 For all \(z\in\CR\), then \(\det\pba{n}{z}\neq0\) and
\beql{Z18}
 \XFa{2n+1}{z}
 =\hp{2n}\ek*{\pba{n}{z}}^\inv\pba{n+1}{z}.
\eeq
 
\eenui
\erem

 Now we are going to prove that the \tsXF{}, which is introduced in \rdefn{K17.1}, can also be expressed in terms of the sequence \(\seq{\pd{k}}{k}{0}{\ev{\kappa}}\) instead of the sequence \(\seq{\pb{k}}{k}{0}{\ev{\kappa}}\).
 This is a consequence of the following result.
 We use the notation given by \eqref{NM}.

\bpropl{KN}
 Let \(\kappa\in\NOinf\) and let \(\seqska\in\Hggeq{\kappa}\) with \tsHp{} \(\sHp{\kappa}\) and \sabcd{} \(\abcd{\ev{\kappa}}\).
 Then:
\benui
 \il{KN.a} If \(k\in\NO\) is such that \(2k-1\leq\kappa\), then \(\nst{\det\pb{k}}\cup\nst{\det\pd{k}}\subseteq\R\).
 \il{KN.b} If \(\kappa\geq 1\), for all \(n\in\NO\) fulfilling \(2n+1\leq\kappa\), then
\begin{multline}\label{KN.1}%
 \hp{2n}\ek*{\pba{n}{z}}^\inv \pba{n+1}{z}
 =\pda{n+1}{z}\ek*{\pda{n}{z}}^\inv \hp{2n}\\
 \text{for all }z\in\C\setminus\rk*{\nst{\det\pb{n}}\cup\nst{\det\pd{n}}}.
\end{multline}
 \il{KN.c} If \(\kappa\geq 1\), for all \(n\in\NO\) fulfilling \(2n+1\leq\kappa\), then
\begin{multline}\label{KN.2}%
 \hp{2n}^\mpi \hp{2n}\ek*{\pba{n+1}{z}}^\inv \pba{n}{z}\hp{2n}^\mpi
 =\hp{2n}^\mpi \pda{n}{z}\ek*{\pda{n+1}{z}}^\inv \hp{2n}\hp{2n}^\mpi\\
 \text{for all }z\in\C\setminus\rk*{\nst{\det\pb{n+1}}\cup\nst{\det\pd{n+1}}}.
\end{multline}
 \eenui
\eprop
\bproof
 \eqref{KN.a} Use \rlem{BK8.6}.
 
 \eqref{KN.b}, \eqref{KN.c} Suppose \(\kappa\geq 1\).
 According to \rrem{R1019}, we have \eqref{R1019.1} and \eqref{R1019.2}. 
 The subsequent part of our proof proceeds inductively.
 From \eqref{K19.0}, \eqref{K19.1}, and \eqref{R1019.1} we discern that \(\nst{\det\pb{0}}=\emptyset\), that \(\nst{\det\pd{0}}=\emptyset\), and that 
\beql{CN3}\begin{split}
 \hp{0}\ek*{\pb{0}\rk{z}}^\inv \pb{1}\rk{z}
 &=\hp{0}\Iq^\inv \rk{z\Iq-\hp{0}^\mpi \hp{1}}
 =z\hp{0}-\hp{0}\hp{0}^\mpi \hp{1}
 =z\hp{0}-\hp{1}\\
 &=z\hp{0}-\hp{1}\hp{0}^\mpi \hp{0}
 = \rk{z\Iq-\hp{1}\hp{0}^\mpi}\Iq^\inv \hp{0}
 =\pd{1}\rk{z}\ek*{\pd{0}\rk{z}}^\inv \hp{0}
\end{split}\eeq
 holds true for all \(z\in\C\).
 Consequently, \eqref{KN.1} is proved for \(n=0\).
 For each \(z\in\C\setminus\rk{\nst{\det\pb{1}}\cup\nst{\det\pd{1}}}\), from \eqref{CN3} we conclude
\[
  \hp{0}\rk*{\ek*{\pb{1}\rk{z}}^\inv \pb{0}\rk{z}}^\inv
  =\hp{0}\ek*{\pb{0}\rk{z}}^\inv \pb{1}\rk{z}
  =\pd{1}\rk{z}\ek*{\pd{0}\rk{z}}^\inv \hp{0}
  =\rk*{\pd{0}\rk{z}\ek*{\pd{1}\rk{z}}^\inv}^\inv \hp{0}
\]
 and, hence, 
\[
 \pd{0}\rk{z}\ek*{\pd{1}\rk{z}}^\inv \hp{0}
 =\hp{0}\ek*{\pb{1}\rk{z}}^\inv \pb{0}\rk{z},
\]
 which implies 
\[
 \hp{0}^\mpi \hp{0}\ek*{\pb{1}\rk{z}}^\inv \pb{0}\rk{z}\hp{0}^\mpi
 =\hp{0}^\mpi \pd{0}\rk{z}\ek*{\pd{1}\rk{z}}^\inv \hp{0}\hp{0}^\mpi.
\]
 Thus, \eqref{KN.2} is checked for \(n=0\).
 Consequently, there is an \(m\in\NO\) fulfilling \(2m+1\leq\kappa\) such that \eqref{KN.1} and \eqref{KN.2} are valid for all \(n\in\mn{0}{m}\).
 If \(2m+1\in\set{\kappa, \kappa-1}\), then \eqref{KN.b} and \eqref{KN.c} are proved.
 Assume that \(2m+1\leq\kappa -2\).
 We consider an arbitrary \(z\in\C\setminus\rk{\nst{\det\pb{m+1}}\cup\nst{\det\pd{m+1}}}\).
 Using \eqref{K19.pb}, \eqref{R1019.1}, and \eqref{R1019.2}, we obtain
\beql{CN3-1}\begin{split}
  &\hp{2m+2}\ek*{\pb{m+1}\rk{z}}^\inv \pb{m+2}\rk{z}\\
  &= \hp{2m+2}\ek*{\pb{m+1}\rk{z}}^\inv \ek*{\pb{m+1}\rk{z}\rk{z\Iq-\hp{2m+2}^\mpi \hp{2m+3}}-\pb{m}\rk{z}\hp{2m}^\mpi \hp{2m+2}}\\
  &= z\hp{2m+2}-\hp{2m+2}\hp{2m+2}^\mpi \hp{2m+3}-\hp{2m+2}\ek*{\pb{m+1}\rk{z}}^\inv \pb{m}\rk{z}\hp{2m}^\mpi \hp{2m+2}\\
  &= z\hp{2m+2}-\hp{2m+3}-\hp{2m+2}\hp{2m}^\mpi \hp{2m}\ek*{\pb{m+1}\rk{z}}^\inv \pb{m}\rk{z}\hp{2m}^\mpi \hp{2m+2}.
\end{split}\eeq
 Similarly, \eqref{K19.pd}, \eqref{R1019.1}, and \eqref{R1019.2} yield
\beql{CN3-2}%
  \pd{m+2}\rk{z}\ek*{\pd{m+1}\rk{z}}^\inv \hp{2m+2}
  = z\hp{2m+2}-\hp{2m+3}-\hp{2m+2}\hp{2m}^\mpi \pd{m}\rk{z}\ek*{\pd{m+1}\rk{z}}^\inv \hp{2m}\hp{2m}^\mpi \hp{2m+2}.
\eeq%
 Since \eqref{KN.2} is assumed to be valid for \(n=m\), the comparison of \eqref{CN3-1} and \eqref{CN3-2} provides that \eqref{KN.1} holds true for \(n=m+1\).
 If \(z\) belongs to \(\C\setminus(\nst{\det\pb{m+1}}\cup\nst{\det\pd{m+1}}\cup\nst{\det\pb{m+2}}\cup\nst{\det\pd{m+2}})\), then \eqref{KN.1}, used with \(n=m+1\), implies
\[
  \hp{2m+2}\rk*{\ek*{\pb{m+2}\rk{z}}^\inv \pb{m+1}\rk{z}}^\inv
  =\rk*{\pd{m+1}\rk{z}\ek*{\pd{m+2}\rk{z}}^\inv}^\inv \hp{2m+2}
\]
 and, consequently,
\[%
  \hp{2m+2}\ek*{\pb{m+2}\rk{z}}^\inv \pb{m+1}\rk{z}
  = \pd{m+1}\rk{z}\ek*{\pd{m+2}\rk{z}}^\inv \hp{2m+2}.
\]
 Continuity arguments yield that this identity is valid for all \(z\in\C\setminus(\nst{\det\pb{m+2}}\cup\nst{\det\pd{m+2}})\).
 Thus, multiplication from the left and from the right by \(\hp{2m+2}^\mpi\) shows that \eqref{KN.2} holds true for \(n=m+1\).
 Consequently, \rpartss{KN.b}{KN.c} are proved by mathematical induction.
\eproof
 
 Now we state the announced alternative description of the \tsXF{}.
 
\bpropl{K17-5}
 Let \(\kappa\in\NOinf\) and let \(\seqska\in\Hggeqka\) with \tsHp{} \(\sHp{\kappa}\) and \sabcd{} \(\abcd{\ev{\kappa}}\).
 Further, let \(\seqX{\kappa}\) be the \tsXFo{\(\seqska\)}.
 Then:
\benui
 \il{K17-5.b} Let \(n\in\NO\) be such that \(2n\leq\kappa\).
 For all \(z\in\CR\), then \(\det\pda{n}{z}\neq0\) and
\beql{K17-5.A}
 \XFa{2n}{z}
 =\pdoa{n+1}{z}\ek*{\pda{n}{z}}^\inv\hp{2n}.
\eeq
 In particular, \(\XF{2n}=\pdo{n+1}\pd{n}^\inv\hp{2n}\).
 \il{K17-5.a} Let \(n\in\NO\) be such that \(2n+1\leq\kappa\).
 For all \(z\in\CR\), then \(\det\pda{n}{z}\neq0\) and
\beql{GR2}
  \XFa{2n+1}{z}
  =\pda{n+1}{z}\ek*{\pda{n}{z}}^\inv\hp{2n}.
\eeq
 In particular, \(\XF{2n+1}=\pd{n+1}\pd{n}^\inv\hp{2n}\).
\eenui
\eprop
\bproof
 \eqref{K17-5.b} Let \(z\in\CR\).
 In the case \(n=0\), equation \eqref{K17-5.A} follows immediately from \eqref{XF.01}, \eqref{K19.0}, and \eqref{abcdO-1}.
 Now suppose \(n\geq1\).
 According to \rpartss{KN.a}{KN.c} of \rprop{KN}, we realize that \(\det\pba{n}{z}\neq0\) and \(\det\pda{n}{z}\neq0\) hold true, and that \(\hp{2n-2}^\mpi \hp{2n-2}\ek{\pba{n}{z}}^\inv\pba{n-1}{z}\hp{2n-2}^\mpi=\hp{2n-2}^\mpi\pda{n-1}{z}\ek{\pda{n}{z}}^\inv \hp{2n-2}\hp{2n-2}^\mpi\) is valid.
 Since \rrem{R1019} shows \(\hp{2n}\hp{2n-2}^\mpi \hp{2n-2}=\hp{2n}\) and \(\hp{2n-2}\hp{2n-2}^\mpi \hp{2n}=\hp{2n}\), multiplication form the left and from the right by \(\hp{2n}\) yields then \(\hp{2n}\ek{\pba{n}{z}}^\inv\pba{n-1}{z}\hp{2n-2}^\mpi\hp{2n}=\hp{2n}\hp{2n-2}^\mpi\pda{n-1}{z}\ek{\pda{n}{z}}^\inv\hp{2n}\).
 Taking additionally into account \rremp{R52SK}{R52SK.b} and \rnota{N-abcdO}, we thus conclude
\[\begin{split}
 \XFa{2n}{z}
 &=\hp{2n}\ek*{\pba{n}{z}}^\inv\pboa{n+1}{z}
 =\hp{2n}\ek{\pba{n}{z}}^\inv\ek*{z\pba{n}{z}-\pba{n-1}{z}\hp{2n-2}^\mpi\hp{2n}}\\
 &=z\hp{2n}-\hp{2n}\ek*{\pba{n}{z}}^\inv\pba{n-1}{z}\hp{2n-2}^\mpi\hp{2n}
 =z\hp{2n}-\hp{2n}\hp{2n-2}^\mpi\pda{n-1}{z}\ek*{\pda{n}{z}}^\inv\hp{2n}\\
 &=\ek*{z\pda{n}{z}-\hp{2n}\hp{2n-2}^\mpi\pda{n-1}{z}}\ek{\pda{n}{z}}^\inv\hp{2n}
 =\pdoa{n+1}{z}\ek*{\pda{n}{z}}^\inv\hp{2n}.
\end{split}\]
 In view of \rrem{TC} and continuity arguments, we perceive \(\XF{2n}=\pdo{n+1}\pd{n}^\inv\hp{2n}\).
 
 \eqref{K17-5.a} Let \(z\in\CR\).
 According to \rpartss{KN.a}{KN.b} of \rprop{KN}, we realize that \(\det\pba{n}{z}\neq0\) and \(\det\pda{n}{z}\neq0\) hold true, and that \(\hp{2n}\ek{\pba{n}{z}}^\inv\pba{n+1}{z}=\pda{n+1}{z}\ek{\pda{n}{z}}^\inv \hp{2n}\) is valid.
 In view of \rremp{R52SK}{R52SK.a}, then \eqref{GR2} follows.
 \rrem{TC} and continuity arguments yield \(\XF{2n+1}=\pd{n+1}\pd{n}^\inv\hp{2n}\).
\eproof

 We proceed by providing essential as well as technical characteristics of the \tsXFo{given} data \(\seqska\), which prove beneficial when treating the description of the set given in \eqref{N128N} generally.

\bleml{K17-2}
 Let \(\kappa\in\Ninf\) and let \(\seqska\in\Hggeqka\) with \tsHp{} \(\sHp{\kappa}\) and \tsXF{} \(\seqX{\kappa}\).
 For every choice of \(n\in\NO\) fulfilling \(2n+1\leq\kappa\), then \(\XF{2n}-\XF{2n+1}=\hp{2n+1}\) and, for all \(z\in\CR\), in particular,
\beql{N47SK}
 \XFa{2n}{z}-\XFa{2n+1}{z}
 =\hp{2n+1}.
\eeq
\elem
\bproof
 Let \(n\in\NO\) be such that \(2n+1\leq\kappa\).
 \rrem{R1019} yields \(\hp{2n}\hp{2n}^\mpi\hp{2n+1}=\hp{2n+1}\).
 In view of \rremss{R52SK}{R1358}, for all \(z\in\CR\), we then acquire
\[\begin{split}%
 \XFa{2n}{z}
 &=\hp{2n}\ek{\pba{n}{z}}^\inv\pboa{n+1}{z}
 =\hp{2n}\ek*{\pba{n}{z}}^\inv\ek*{\pba{n+1}{z}+\pba{n}{z}\hp{2n}^\mpi\hp{2n+1}}\\
 &=\hp{2n}\ek{\pba{n}{z}}^\inv\pba{n+1}{z}+\hp{2n}\hp{2n}^\mpi\hp{2n+1}
 =\XFa{2n+1}{z}+\hp{2n+1}.
\end{split}\]
 Regarding \rdefn{K17.1} and \rrem{TC}, continuity arguments complete the proof.
\eproof

 The following result contains an interesting connection between a sequence \(\seqs{2n}\in\Hggeq{2n}\) and its \tnatext{} \(\seqsh{2n+1}\), which describes the interplay between sequences of \tXF{s}.

\bleml{L1520}
 Let \(\kappa\in\NOinf\), let \(\seqska\in\Hggeqka\) with \tsXF{} \(\seqX{\kappa}\), and let \(n\in\NO\) with \(2n\leq\kappa\).
 Then the \tnatext{} \(\seqsh{2n+1}\) of \(\seqs{2n}\) belongs to \(\Hggeq{2n+1}\) and the \tsXF{} \(\seqXh{2n+1}\) associated with \(\seqsh{2n+1}\) fulfills \(\XFh{2n+1}=\XF{2n}\) and, for all \(z\in\CR\), in particular, \(\XFh{2n+1}\rk{z}=\XFa{2n}{z}\).
\elem
\bproof
 \rprop{L1237-A} yields \(\seqsh{2n+1}\in\Hggeq{2n+1}\).
 Let \(\shHp{2n+1}\) be the \tsHpo{\(\seqsh{2n+1}\)}.
 Then, the application of \rlem{K17-2} to the sequence \(\seqsh{2n+1}\) provides \(\XFh{2n}-\XFh{2n+1}=\hhp{2n+1}\).
 According to \rrem{R26-SK1}, we have \(\hhp{2n+1}=\Oqq\).
 From \rrem{XF-tru} we can conclude \(\XFh{2n}=\XF{2n}\).
 Consequently, we obtain \(\XFh{2n+1}=\XFh{2n}-\hhp{2n+1}=\XFh{2n}=\XF{2n}\).
 Regarding \rremss{TC}{R0735}, the proof is complete.
\eproof

\bleml{L0625}%
 Let \(\kappa\in\NOinf\), let \(\seqska\in\Hggeqka\) with \tsXF{} \(\seqX{\kappa}\), and let \(m\in\mn{-1}{\kappa}\).
 Then \(\ek{\XFa{m}{z}}^\ad=\XFa{m}{\ko{z}}\) for all \(z\in\holpt{\XF{m}}\).
 In particular, \(\ek{\XFa{m}{x}}^\ad=\XFa{m}{x}\) for all \(x\in\R\cap\holpt{\XF{m}}\).
\elem
\bproof
 \rrem{R1019} yields \eqref{R1019.0}.
 By virtue of \eqref{XF.01}, we see then \(\XF{-1}\rk{\ko{z}}=\Oqq=\ek{\XF{-1}\rk{z}}^\ad\) and \(\XF{0}\rk{\ko{z}}
 =\ko{z}\hp{0}
 =\ko{z}\hp{0}^\ad
 =\ek{\XF{0}\rk{z}}^\ad\) for all \(z\in\C\).
 Thus, in the case \(m\leq0\), the proof is complete.
 Now assume \(m\geq1\).
 First we consider the case that \(m=2n\) is fulfilled with some \(n\in\N\).
 According to \eqref{R1019.0}, we have \(\hp{2n}^\ad=\hp{2n}\) and \(\hp{2n-2}^\ad=\hp{2n-2}\).
 Consequently, \rrem{A.R.A++*} shows that \(\rk{\hp{2n-2}^\mpi}^\ad=\hp{2n-2}^\mpi\) is valid.
 Because of \rrem{R1624}, we can infer from \rrem{BK8.2} furthermore \(\pda{n}{\ko{z}}=\ek{\pba{n}{z}}^\ad\) and \(\pda{n-1}{\ko{z}}=\ek{\pba{n-1}{z}}^\ad\) for all \(z\in\C\).
 For all \(z\in\C\setminus\nst{\det\pb{n}}\), then \(\det\pda{n}{\ko{z}}\neq0\) and, using \rpropp{K17-5}{K17-5.b}, \rnota{N-abcdO}, and \rdefn{K17.1}, we receive moreover
\begin{multline*}
 \XF{2n}\rk{\ko{z}}
 =\pdoa{n+1}{\ko z}\ek*{\pda{n}{\ko z}}^\inv\hp{2n}
 =\ek*{\ko{z}\pda{n}{\ko{z}}-\hp{2n}\hp{2n-2}^\mpi\pd{n-1}\rk{\ko{z}}}\ek*{\pda{n}{\ko{z}}}^\inv\hp{2n}\\
 =\rk*{\ko{z}\ek*{\pba{n}{z}}^\ad-\hp{2n}^\ad (\hp{2n-2}^\mpi)^\ad\ek*{\pba{n-1}{z}}^\ad}\ek*{\pba{n}{z}}^\invad\hp{2n}^\ad
 =\ek*{\pboa{n+1}{z}}^\ad\ek*{\pba{n}{z}}^\invad\hp{2n}^\ad
 =\ek*{\XFa{2n}{z}}^\ad.
\end{multline*}
 In the case that \(m=2n+1\) is fulfilled with some \(n\in\NO\), 
 using \rpropp{K17-5}{K17-5.a} instead of \rpropp{K17-5}{K17-5.b}, we obtain similarly
\[
 \XF{2n+1}\rk{\ko{z}}
 =\pda{n+1}{\ko{z}}\ek*{\pda{n}{\ko{z}}}^\inv\hp{2n}
 =\ek*{\pba{n+1}{z}}^\ad\ek*{\pba{n}{z}}^\invad\hp{2n}^\ad
 =\ek*{\XFa{2n+1}{z}}^\ad
\]
 for all \(z\in\C\setminus\nst{\det\pb{n}}\).
 Taking into account \eqref{SK5-11}, we have thus shown \(\XF{m}\rk{\ko{z}}=\ek{\XFa{m}{z}}^\ad\) for all \(z\in\C\setminus\nst{\det\pb{\ef{m}}}\).
 By continuity, then \(\ek{\XFa{m}{z}}^\ad=\XFa{m}{\ko{z}}\) follows for all \(z\in\holpt{\XF{m}}\).
\eproof

\bleml{K17-6}
 Let \(\kappa\in\NOinf\), let \(\seqska\in\Hggeqka\), and let \(\sHp{\kappa}\) be the \tsHpo{\(\seqska\)}.
 Further, let \(\seqX{\kappa}\) be the \tsXFo{\(\seqska\)}.
 For all \(m\in\mn{0}{\kappa}\) and all \(z\in\CR\), then \(\ran{\XFa{m}{z}}=\ran{\hp{\eff{m}}}\) and \(\nul{\XFa{m}{z}}=\nul{\hp{\eff{m}}}\).
\elem
\bproof
 Let \(z\in\CR\).
 \rlem{BK8.6} yields \(\det\pba{k}{z}\neq0\) and \(\det\pda{k}{z}\neq0\) for all \(k\in\NO\) with \(2k-1\leq\kappa\).
 Furthermore, \rlem{K15-1} yields \(\det\pboa{n+1}{z}\neq0\) and \(\det\pdoa{n+1}{z}\neq0\) for all \(n\in\NO\) with \(2n\leq\kappa\).
 Taking into account \rrem{R52SK} and \rprop{K17-5}, we can conclude from \rrem{R1241} then \(\ran{\XFa{2n+1}{z}}=\ran{\hp{2n}}\) and \(\nul{\XFa{2n+1}{z}}=\nul{\hp{2n}}\) for all \(n\in\NO\) with \(2n+1\leq\kappa\) as well as \(\ran{\XFa{2n}{z}}=\ran{\hp{2n}}\) and \(\nul{\XFa{2n}{z}}=\nul{\hp{2n}}\) for all \(n\in\NO\) with \(2n\leq\kappa\).
 In view of \eqref{SK5-11}, the proof is complete.
\eproof

 Now we consider Moore--Penrose inverses of the matrix-valued functions belonging to the \tsXFo{the} given data \(\seqska\).

\bleml{L0946}
 Let \(\kappa\in\Ninf\) and let \(\seqska\in\Hggeqka\) with \tsHp{} \(\sHp{\kappa}\), \sabcd{} \(\abcd{\ev{\kappa}}\), and \tsXF{} \(\seqX{\kappa}\).
 Then \(\det\pba{n+1}{z}\neq0\) and \(\ek{\XFa{2n+1}{z}}^\mpi=\hp{2n}^\mpi\hp{2n}\ek{\pba{n+1}{z}}^\inv\pba{n}{z}\hp{2n}^\mpi\) as well as \(\det\pda{n+1}{z}\neq0\) and \(\ek{\XFa{2n+1}{z}}^\mpi=\hp{2n}^\mpi\pda{n}{z}\ek{\pda{n+1}{z}}^\inv\hp{2n}\hp{2n}^\mpi\) are valid for all \(n\in\NO\) such that \(2n+1\leq\kappa\) and all \(z\in\CR\).
\elem
\bproof
 Let \(n\in\NO\) be such that \(2n+1\leq\kappa\) and let \(z\in\CR\).
 In view of \rlem{BK8.6}, all the matrices \(B_{n+1}\defeq\pba{n+1}{z}\), \(B_{n}\defeq\pba{n}{z}\), \(D_{n+1}\defeq\pda{n+1}{z}\), and \(D_{n}\defeq\pda{n}{z}\) are invertible.
 We set \(L\defeq D_{n}D_{n+1}^\inv\), \(R\defeq B_{n+1}^\inv B_{n}\), \(M\defeq\hp{2n}\), \(N\defeq LMR^\inv\), and \(X\defeq MR^\inv\).
 Then the matrices \(L\) and \(R\) are invertible and \(R^\inv=B_{n}^\inv B_{n+1}\).
 Consequently, \(N=D_{n}D_{n+1}^\inv MB_{n}^\inv B_{n+1}\) and \(X=MB_{n}^\inv B_{n+1}\).
 \rpropp{KN}{KN.b} provides \(MB_{n}^\inv B_{n+1}=D_{n+1}D_{n}^\inv M\).
 Hence, 
\[
 M
 =D_{n}D_{n+1}^\inv D_{n+1}D_{n}^\inv M
 =D_{n}D_{n+1}^\inv MB_{n}^\inv B_{n+1}
 =N.
\]
 Applying \rlem{KB-5}, then we conclude 
\begin{align}
 X^\mpi
 &=N^\mpi NRM^\mpi
 =M^\mpi MB_{n+1}^\inv B_{n}M^\mpi,&%
 X^\mpi
 &=N^\mpi LMM^\mpi
 =M^\mpi D_{n}D_{n+1}^\inv MM^\mpi.\label{L0946.2}
\end{align}
 On the other hand, from \rrem{R52SK} we recognize that \(\XFa{2n+1}{z}=MB_{n}^\inv B_{n+1}=X\) holds true.
 Comparing this to \eqref{L0946.2} completes the proof.
\eproof

\bleml{K17-7}
 Let \(\kappa\in\NOinf\) and let \(\seqska\in\Hggeqka\) with \tsHp{} \(\sHp{\kappa}\) and \tsXF{} \(\seqX{\kappa}\).
 Then:
\benui
 \il{K17-7.a} \(\XFa{2n}{z}=z\hp{2n}-\hp{2n}\ek{\XFa{2n-1}{z}}^\mpi\hp{2n}\) for all \(n\in\NO\) fulfilling \(2n\leq\kappa\) and all \(z\in\CR\).
 \il{K17-7.b} If \(\kappa\geq1\), then \(\XFa{2n+1}{z}=z\hp{2n}-\hp{2n}\ek{\XFa{2n-1}{z}}^\mpi\hp{2n}-\hp{2n+1}\) for all \(n\in\NO\) fulfilling \(2n+1\leq\kappa\) and all \(z\in\CR\).
\eenui
\elem
\bproof
 \eqref{K17-7.a} Let \(n\in\NO\) be such that \(2n\leq\kappa\) and let \(z\in\CR\).
 We first consider the case \(n=0\).
 By virtue of \rdefn{K17.1} and \eqref{XF.01}, we see then \(z\hp{0}-\hp{0}\ek{\XFa{-1}{z}}^\mpi\hp{0}=z\hp{0}-\hp{0}\cdot\Oqq^\mpi\cdot\hp{0}=z\hp{0}=\XFa{0}{z}\).
 Now suppose \(n\geq1\).
 From \rlem{L0946} we can infer \(\det\pba{n}{z}\neq0\) and
\beql{FL1}
 \ek*{\XFa{2n-1}{z}}^\mpi
 =\hp{2n-2}^\mpi\hp{2n-2}\ek{\pba{n}{z}}^\inv\pba{n-1}{z}\hp{2n-2}^\mpi.
\eeq 
 \rrem{R1019} yields \(\hp{2n}\hp{2n-2}^\mpi\hp{2n-2}=\hp{2n}\).
 Consequently, from \eqref{FL1} we perceive \(\hp{2n}\ek{\XFa{2n-1}{z}}^\mpi=\hp{2n}\ek{\pba{n}{z}}^\inv\pba{n-1}{z}\hp{2n-2}^\mpi\).
 In view of \rrem{R52SK} and \rnota{N-abcdO}, we have \(\XFa{2n}{z}=\hp{2n}\ek{\pba{n}{z}}^\inv\ek{z\pba{n}{z}-\pba{n-1}{z}\hp{2n-2}^\mpi\hp{2n}}\).
 Therefore, 
\[
 z\hp{2n}-\hp{2n}\ek*{\XFa{2n-1}{z}}^\mpi\hp{2n}
 =z\hp{2n}-\rk*{\hp{2n}\ek{\pba{n}{z}}^\inv \pba{n-1}{z}\hp{2n-2}^\mpi}\hp{2n}
 =\XFa{2n}{z}.
\]

 \eqref{K17-7.b} Suppose \(\kappa\geq1\) and that \(n\in\NO\) is such that \(2n+1\leq\kappa\).
 Then \rlem{K17-2} shows that \eqref{N47SK} holds true for all \(z\in\CR\).
 Thus, applying \rpart{K17-7.a} completes the proof of \rpart{K17-7.b}.
\eproof

 We note that \rlem{K17-7} plays an essential role for \rlem{P1056} and \rprop{P1129} below, which contain significant observations about the \tsXF{}.
 Now we are going to look at the \tsXF{} under the light of a particular Schur transform, which was introduced in \cite{MR3380267}.
 Before doing this we recall the following notion.

\bdefnnl{\zitaa{MR3380267}{\csec{8}{218}}}{K11}%
 Let \(\mG\) be a \tne{} subset of \(\C\), let \(F\colon\mG\to\Cpq\) be a matrix-valued function, and let \(A\) and \(B\) be complex \tpqa{matrices}.
 Then let \(F^\ABFt{A}{B}\colon\mG\to\Cpq\) be defined by \(F^\ABFt{A}{B}\rk{z}\defeq-A\rk{z\Iq+\ek{F\rk{z}}^\mpi A}+B\).
 The function \(F^\ABFt{A}{B}\) is called the \emph{\tABSchto{A}{B}{\(F\)}}.
\edefn

\bleml{P1056}
 Let \(\kappa\in\NOinf\), let \(\seqska\in\Hggeqka\), and let \(\sHp{\kappa}\) be the \tsHpo{\(\seqska\)}.
 Further, let \(\seqX{\kappa}\) be the \tsXFo{\(\seqska\)}.
 Then \(\XFa{2n}{z}=\XF{2n-1}^\ABFt{-\hp{2n}}{\Oqq}\rk{z}\) for all \(n\in\NO\) such that \(2n\leq\kappa\) and all \(z\in\CR\).
 Moreover \(\XFa{2n+1}{z}=\XF{2n-1}^\ABFt{-\hp{2n}}{-\hp{2n+1}}\rk{z}\) for all \(n\in\NO\) fulfilling \(2n+1\leq\kappa\) and all \(z\in\CR\).
\elem
\bproof
 This is a direct consequence of \rdefn{K11} and \rlem{K17-7}.
\eproof

 Now we state an interrelation to the class \(\RFq\).

\bpropl{P1129}
 Let \(\kappa\in\NOinf\), let \(\seqska\in\Hggeqka\), and let \(\sHp{\kappa}\) be the \tsHpo{\(\seqska\)}.
 Further, let \(\seqX{\kappa}\) be the \tsXFo{\(\seqska\)}.
 For all \(m\in\mn{0}{\kappa}\), then \(G_m\defeq\rstr_{\pip}\XF{m}\) belongs to the class \(\RFq\) and \(\beta_{m}=\hp{\eff{m}}\), where \((\alpha_{m}, \beta_{m}, \nu_{m})\) is the \tNpo{\(G_m\)} and \(\eff{m}\) is given by \eqref{SK5-11}.
\eprop
\bproof
 \rprop{DFKMT2.30N} yields \(\hp{j}^\ad=\hp{j}\) for all \(j\in\mn{0}{\kappa}\) and \(\hp{2k}\in\Cggq\) for all \(k\in\NO\) fulfilling \(2k\leq\kappa\) hold true.
 The following considerations of our proof are divided into three parts:
 
 (I) From \eqref{XF.01} we get \(G_0\rk{w}=w\hp{0}\) for all \(w\in\pip\).
 Regarding \(\hp{0}\in\Cggq\), \rthm{Theo21} yields \(G_0\in\RFOq \) and \(\beta_{0}=\hp{0}\).
 If \(\kappa=0\), then, in view of \eqref{SK5-11}, the proof is complete.
  
 (II) Suppose \(\kappa\geq1\).
 From \eqref{XF.01} we receive \(G_1\rk{w}=-\hp{1}+w\hp{0}\) for all \(w\in\pip\).
 Regarding \(\hp{0}\in\Cggq\) and \(\hp{1}^\ad=\hp{1}\), consequently \rthm{Theo21} provides \(G_1\in\RFq\) and \(\beta_1=\hp{0}\).
 If \(\kappa=1\), then, in view of \eqref{SK5-11}, the proof is finished.
 
 (III) Suppose \(\kappa\geq2\).
 According to parts~(I) and~(II) of the proof, there is an \(n\in\N\) fulfilling \(2n\leq\kappa\) such that \(G_l\in\RFq\) and \(\beta_{ l}=\hp{\eff{l}}\) hold true for each \(l\in\mn{0}{2n-1}\).
 \rlem{P1056} shows that \(\XFa{2n}{z}=\XF{2n-1}^\ABFt{-\hp{2n}}{\Oqq}\rk{z}\) is valid for all \(z\in\CR\).
 In particular, \(G_{2n}=G_{2n-1}^\ABFt{-\hp{2n}}{\Oqq}\).
 Taking into account \(G_{2n-1}\in\RFq\) as well as \(\hp{2n}\in\Cggq\) and \(\Oqq\in\CHq\), from \zitaa{MR3380267}{\cprop{8.6}{219}} we conclude \(G_{2n}\in\RFq\).
 According to \rlem{fkm12bP3.3}, we acquire
\begin{align}\label{SKM11}
 \beta_{2n}&=\lim_{y\to\infty}\rk{\iu y}^\inv G_{2n}\rk{\iu y}&
&\text{and}&
 \beta_{2n-1}&=\lim_{y\to\infty}\rk{\iu y}^\inv G_{2n-1}\rk{\iu y}.
\end{align}
 Because of \rlem{K17-6} and \(\beta_{2n-1}=\hp{2n-2}\), we obtain
\[
 \rank G_{2n-1}\rk{w}
 =\dim\ran{G_{2n-1}\rk{w}}
 =\dim\ran{\hp{2n-2}}
 =\rank\hp{2n-2}
 =\rank\beta_{2n-1} 
\]
 for all \(w\in\pip\).
 Therefore, the second equation in \eqref{SKM11} and \rprop{MPK} provide \(\lim_{y\to\infty}\rk{\iu y\ek{G_{2n-1}\rk{\iu y}}^\mpi}=\beta_{2n-1}^\mpi\) and, consequently, 
\beql{SKM13}
 \lim_{y\to\infty}\rk*{\rk{\iu y}^\inv\ek*{G_{2n-1}\rk{\iu y}}^\mpi}
 =\Oqq.
\eeq
 Using \eqref{SKM11}, \rlemp{K17-7}{K17-7.a}, and \eqref{SKM13}, we perceive
\beql{SKM}\begin{split}
 \beta_{2n}
 &=\lim_{y\to\infty}\rk{\iu y}^\inv\rk*{\iu y\hp{2n}-\hp{2n}\ek*{G_{2n-1}\rk{\iu y}}^\mpi\hp{2n}}\\
 &=\hp{2n}-\hp{2n}\ek*{\lim_{y\to\infty}\rk*{\rk{\iu y}^\inv\ek*{G_{2n-1}\rk{\iu y}}^\mpi}}\hp{2n}
 =\hp{2n}.
\end{split}\eeq
 Now assume that \(2n+1\leq\kappa\).
 From \rlem{K17-2} we recognize that \eqref{N47SK} holds true for each \(z\in\CR\).
 Thus, by virtue of \eqref{N47SK} and \(G_{2n}\in\RFq\), from \rthm{Theo21} we conclude 
\[
 G_{2n+1}\rk{w}
 =G_{2n}\rk{w}-\hp{2n+1}
 =\alpha_{2n}-\hp{2n+1}+w\beta_{2n}+\int_\R\frac{1+xw}{x-w}\nu_{2n}\rk{\dif x}
\]
 for all \(w\in\pip\).
 In view of \(\hp{2n+1}^\ad=\hp{2n+1}\), \rthm{Theo21}, and \eqref{SKM}, consequently, \(G_{2n+1}\in\RFq\) and \(\beta_{2n+1}=\beta_{2n}=\hp{2n}\) hold true.
 For each \(m\in\mn{0}{\kappa}\), in view of \eqref{SK5-11}, hence, \(G_m\in\RFq\) and \(\beta_{m}=\hp{\eff{m}}\) are proved inductively.
\eproof

\bcorl{C1213}
 Let \(\kappa\in\NOinf\), let \(\seqska\in\Hggeqka\) with \sabcd{} \(\abcd{\ev{\kappa}}\) and \tsXF{} \(\seqX{\kappa}\), and let \(m\in\mn{0}{\kappa}\).
 In view of \rprop{P1129}, denote by \(\rk{\alpha_m,\beta_m,\nu_{m}}\) the \tNpo{\(G_m\defeq\rstr_{\pip}\XF{m}\)}.
 Then \(\nu_m\rk{\R\setminus\nst{\det\pb{\ef{m}}}}=\Oqq\), where \(\ef{m}\) is given by \eqref{SK5-11}.
\ecor
\bproof
 From \rrem{SK32R} we can conclude that \(\det\pb{\ef{m}}\) is a non-zero complex polynomial.
 Consequently, \(\nst{\det\pb{\ef{m}}}\) is a finite subset of \(\C\).
 Hence, there are an integer \(\ell\in\N\) and a sequence \(\seq{\mathcal{I}_k}{k}{1}{\ell}\) of pairwise disjoint open (real, bounded or unbounded) intervals such that \(\R\setminus\nst{\det\pb{\ef{m}}}=\bigcup_{k=1}^\ell\mathcal{I}_\ell\).
 Consider an arbitrary \(k\in\mn{1}{\ell}\).
 Using \rrem{R0735}, we can infer then \(\pip\cup\mathcal{I}_k\subseteq\holpt{\XF{m}}\).
 Therefore, \(G_{m,k}\defeq\rstr_{\pip\cup\mathcal{I}_k}\XF{m}\) is a continuous continuation of \(G_m\) onto \(\pip\cup\mathcal{I}_k\).
 \rlem{L0625} furthermore shows \(G_{m,k}\rk{\mathcal{I}_k}\subseteq\CHq\).
 Since \rprop{P1129} yields \(G_m\in\RFq\), we can then apply \zitaa{MR2570113}{\crem{B.3}{823}} to obtain \(\nu_m\rk{\mathcal{I}_k}=\Oqq\).
 Hence, \(\nu_m\rk{\R\setminus\nst{\det\pb{\ef{m}}}}=\sum_{k=1}^\ell\nu_m\rk{\mathcal{I}_k}=\Oqq\) follows.
\eproof

\bcorl{C1232}
 Let \(\kappa\in\NOinf\), let \(\seqska\in\Hggeqka\) with \tsXF{} \(\seqX{\kappa}\), and let \(m\in\mn{-1}{\kappa}\).
 Then \(\frac{1}{\rim z}\rim\XFa{m}{z}\in\Cggq\) for every choice of \(z\in\CR\).
\ecor
\bproof
 For \(m=-1\) this is an immediate consequence or \rdefn{K17.1}.
 Now suppose \(m\in\mn{0}{\kappa}\).
 Because of the definition of the class \(\RFq\) and \rprop{P1129}, we realize that \(\frac{1}{\rim w}\rim\XFa{m}{w}\in\Cggq\) for all \(w\in\pip\).
 Now we consider an arbitrary \(z\in\lhp\).
 Then \(w\defeq\ko{z}\) belongs to \(\pip\) and, consequently, as already shown, \(\frac{1}{\rim w}\rim\XFa{m}{w}\in\Cggq\) holds true.
 Because of \rrem{R0735} and \rlem{L0625}, we acquire \(\ek{\XFa{m}{z}}^\ad=\XFa{m}{w}\), implying \(\rim\XFa{m}{w}=-\rim\XFa{m}{z}\).
 Thus, additionally taking into account \(\rim w=-\rim z\), we receive \(\frac{1}{\rim z}\rim\XFa{m}{z}\in\Cggq\).
 In view of \(\CR=\pip\cup\lhp\), the proof is complete.
\eproof

 Now we extend the results of \rlem{K17-6}.

\bpropl{L1347}
 Let \(\kappa\in\NOinf\), let \(\seqska\in\Hggeqka\) with \tsHp{} \(\sHp{\kappa}\) and \tsXF{} \(\seqX{\kappa}\), and let \(m\in\mn{0}{\kappa}\).
 Then, for each \(z\in\CR\), the equations
\beql{L1347.1}%
 \Ran{\ek*{\XFa{m}{z}}^\ad}
 =\Ran{\XFa{m}{z}}
 =\Ran{\rim\XFa{m}{z}}
 =\ran{\hp{\eff{m}}}
\eeq
 and
\beql{L1347.2}
 \Nul{\ek*{\XFa{m}{z}}^\ad}
 =\Nul{\XFa{m}{z}}
 =\Nul{\rim\XFa{m}{z}}
 =\nul{\hp{\eff{m}}}
\eeq
 hold true, where \(\eff{m}\) is given by \eqref{SK5-11}.
\eprop
\bproof 
 (I) First we consider an arbitrary \(w\in\pip\).
 By virtue of \rprop{P1129}, the matrix-valued function \(G_m\defeq\rstr_{\pip}\XF{m}\) belongs to \(\RFq\) with \(\beta_m=\hp{\eff{m}}\), where \((\alpha_m, \beta_m, \nu_m)\) signifies the \tNpo{\(G_m\)}.
 Thus, from \zitaa{MR2988005}{\crem{3.5}{1772}} we recognize that \(\ran{\ek{G_m\rk{w}}^\ad}=\ran{G_{m}\rk{w}}\) and \(\nul{\ek{G_{m}\rk{w}}^\ad}=\nul{G_{m}\rk{w}}\) are valid.
 Hence, \(\ran{\ek{\XFa{m}{w}}^\ad}=\ran{\XFa{m}{w}}\) and \(\nul{\ek{\XFa{m}{w}}^\ad}=\nul{\XFa{m}{w}}\).
 Setting \(\gamma_m\defeq\nu_m\rk{\R}\), furthermore \zitaa{MR2988005}{\cprop{3.7}{1772}} yields \(\ran{G_m\rk{w}}=\ran{\alpha_m}+\ran{\beta_m}+\ran{\gamma_m}\) and \(\nul{G_m\rk{w}}=\nul{\alpha_m}\cap\nul{\beta_m}\cap\nul{\gamma_m}\) as well as \(\ran{\rim G_m\rk{w}}=\ran{\beta_m}+\ran{\gamma_m}\) and \(\nul{\rim G_m\rk{w}}=\nul{\beta_m}\cap\nul{\gamma_m}\).
 Taking into account \rlem{K17-6} and \(\beta_m=\hp{\eff{m}}\), we can infer
\[%
 \Ran{\XFa{m}{w}}
 =\ran{\hp{\eff{m}}}
 =\ran{\beta_m}
 \subseteq\Ran{\rim G_m\rk{w}}
 \subseteq\Ran{G_m\rk{w}}
 =\Ran{\XFa{m}{w}}
\]
 and 
\[%
 \Nul{\XFa{m}{w}}
 =\Nul{G_{m}\rk{w}}
 \subseteq\Nul{\rim  G_m\rk{w}}
 \subseteq\nul{\beta_m}
 =\nul{\hp{\eff{m}}}
 =\Nul{\XFa{m}{w}}.
\]
 Consequently, \(\ran{\rim\XFa{m}{w}}=\ran{\rim G_{m}\rk{w}}=\ran{\hp{\eff{m}}}=\ran{\XFa{m}{w}}\) and \(\nul{\rim\XFa{m}{w}}=\nul{\rim G_{m}\rk{w}}=\nul{\hp{\eff{m}}}=\nul{\XFa{m}{w}}\).
 In particular, \eqref{L1347.1} and \eqref{L1347.2} hold true for all \(z\in\pip\).
  
 (II) Let \(z\in\lhp\).
 Then \(w\defeq\ko{z}\) belongs to \(\pip\).
 According to \rrem{R0735} and \rlem{L0625}, we get \(\XFa{m}{w}=\ek{\XFa{m}{z}}^\ad\).
 Hence, \(\XFa{m}{z}=\ek{\XFa{m}{w}}^\ad\) and \(\rim\XFa{m}{z}=-\rim\XFa{m}{w}\) follow.
 Consequently, we obtain \eqref{L1347.1} and \eqref{L1347.2} for all \(z\in\lhp\) as well by applying \cpart{(I)}.
 
 (III) Because of~(I),~(II), and \(\CR=\pip\cup\lhp\), the assertion is proved.
\eproof

\bcorl{C1236}
 Let \(\kappa\in\NOinf\), let \(\seqska\in\Hgeqka\) with \tsXF{} \(\seqX{\kappa}\), and let \(m\in\mn{0}{\kappa}\).
 Then \(\det\XFa{m}{z}\neq0\) for all \(z\in\CR\).
 Furthermore \(\rim\XFa{m}{z}\in\Cgq\) for all \(z\in\pip\).
\ecor
\bproof
 \rrem{Schr2.7} yields \(\seqska\in\Hggeqka\).
 Furthermore, we see easily that \(\hp{\eff{m}}\in\Cgq\) and in particular \(\ran{\hp{\eff{m}}}=\Cq\).
 Applying \rprop{L1347} and \rcor{C1232} completes the proof.
\eproof

\bpropl{M2140}
 Let \(\kappa\in\minf{3}\cup\set{\infi}\), let \(\seqska\in\Hggeqka\) with \tsHp{} \(\sHp{\kappa}\), \sabcd{} \(\abcd{\ev{\kappa}}\), and \tsXF{} \(\seqX{\kappa}\), and let \(n\in\N\) be such that \(2n+1\leq\kappa\).
 For every choice of \(z\) and \(w\) in \(\CR\), then \(\det\pda{n}{z}\neq0\) and \(\det\pba{n}{\ko{w}}\neq0\) as well as
\begin{multline*}
 \XFa{2n+1}{z}-\ek*{\XFa{2n+1}{w}}^\ad
 =\rk{z-\ko{w}}\hp{2n}+\hp{2n}\ek*{\pba{n}{\ko{w}}}^\inv\pb{n-1}\rk{\ko{w}}\hp{2n-2}^\mpi\\
 \lbtimes\rk*{\XFa{2n-1}{z}-\ek*{\XFa{2n-1}{w}}^\ad}\hp{2n-2}^\mpi\pda{n-1}{z}\ek*{\pda{n}{z}}^\inv\hp{2n}.
\end{multline*}
\eprop
\bproof
 According to \rrem{R1019}, we have \eqref{R1019.2}.
 Consider arbitrary \(z,w\in\CR\).
 We set \(Z\defeq\XFa{2n-1}{z}\) and \(W\defeq\XF{2n-1}\rk{\ko{w}}\).
 Regarding \eqref{SK5-11}, from \rprop{L1347} we recognize that \(\ran{Z^\ad}=\ran{Z}=\ran{\hp{2n-2}}=\ran{W}=\ran{W^\ad}\) is valid.
 Using \rrem{tsb3}, we conclude then
\begin{align}\label{M2140.1}
 Z^\mpi Z
 &=\OPu{\ran{Z^\ad}}
 =\OPu{\ran{W^\ad}}
 =W^\mpi W&%
&\text{and}&
 WW^\mpi
 &=\OPu{\ran{W}}
 =\OPu{\ran{Z}}
 =ZZ^\mpi.
\end{align}
 In view of \rrem{R0735}, \rlem{L0625} yields \(W=\ek{\XFa{2n-1}{w}}^\ad\) and
\beql{KD4}
 \XFa{2n+1}{z}-\ek*{\XFa{2n+1}{w}}^\ad
 =\XFa{2n+1}{z}-\XF{2n+1}\rk{\ko{w}}.
\eeq
 According to \rlem{L0946}, we have \(\det\pda{n}{z}\neq0\) and \(\det\pba{n}{\ko{w}}\neq0\) as well as
\begin{align}
 Z^\mpi&=\hp{2n-2}^\mpi\pda{n-1}{z}\ek*{\pda{n}{z}}^\inv\hp{2n-2}\hp{2n-2}^\mpi,&%
 W^\mpi&=\hp{2n-2}^\mpi\hp{2n-2}\ek*{\pba{n}{\ko{w}}}^\inv\pb{n-1}\rk{\ko{w}}\hp{2n-2}^\mpi.\label{KD8}
\end{align}
 From \rlemp{K17-7}{K17-7.b} we obtain \(\XFa{2n+1}{z}=z\hp{2n}-\hp{2n}Z^\mpi\hp{2n}-\hp{2n+1}\) and \(\XF{2n+1}\rk{\ko{w}}=\ko{w}\hp{2n}-\hp{2n}W^\mpi\hp{2n}-\hp{2n+1}\).
 Then, using additionally \eqref{KD4}, \eqref{M2140.1}, \eqref{KD8}, and \eqref{R1019.2}, we receive
\[\begin{split}
 &\XFa{2n+1}{z}-\ek*{\XFa{2n+1}{w}}^\ad
 =\XFa{2n+1}{z}-\XF{2n+1}\rk{\ko{w}}\\
 &=\rk{z\hp{2n}-\hp{2n}Z^\mpi\hp{2n}-\hp{2n+1}}-\rk{\ko{w}\hp{2n}-\hp{2n}W^\mpi\hp{2n}-\hp{2n+1}}\\
 &=\rk{z-\ko{w}}\hp{2n}-\hp{2n}\rk{Z^\mpi-W^\mpi}\hp{2n}
 =\rk{z-\ko{w}}\hp{2n}-\hp{2n}\rk{Z^\mpi ZZ^\mpi-W^\mpi WW^\mpi}\hp{2n}\\
 &=\rk{z-\ko{w}}\hp{2n}-\hp{2n}\rk{W^\mpi WZ^\mpi-W^\mpi ZZ^\mpi}\hp{2n}\\
 &=\rk{z-\ko{w}}\hp{2n}-\hp{2n}W^\mpi\rk{W-Z}Z^\mpi\hp{2n}\\
 &=\rk{z-\ko{w}}\hp{2n}-\hp{2n}\rk*{\hp{2n-2}^\mpi\hp{2n-2}\ek*{\pba{n}{\ko{w}}}^\inv\pb{n-1}\rk{\ko{w}}\hp{2n-2}^\mpi}\rk{W-Z}\\
 &\qquad\lbtimes\rk*{\hp{2n-2}^\mpi\pda{n-1}{z}\ek*{\pda{n}{z}}^\inv\hp{2n-2}\hp{2n-2}^\mpi}\hp{2n}\\
 &=\rk{z-\ko{w}}\hp{2n}+\hp{2n}\ek*{\pba{n}{\ko{w}}}^\inv\pb{n-1}\rk{\ko{w}}\hp{2n-2}^\mpi\rk{Z-W}\hp{2n-2}^\mpi\pda{n-1}{z}\ek*{\pda{n}{z}}^\inv\hp{2n}.
\end{split}\]
 In view of \(W=\ek{\XFa{2n-1}{w}}^\ad\), the proof is complete.
\eproof

\bcorl{M2145}
 Let \(\kappa\in\minf{3}\cup\set{\infi}\), let \(\seqska\in\Hggeqka\) with \tsHp{} \(\sHp{\kappa}\), \sabcd{} \(\abcd{\ev{\kappa}}\), and  \tsXF{} \(\seqX{\kappa}\), and let  \(n\in\N\) be such that \(2n+1\leq\kappa\).
 For all \(z\in\CR\), then
\[
 \rim\XFa{2n+1}{z}
 =\ek*{\rim\rk{z}}\hp{2n}+\hp{2n}\ek*{\pba{n}{\ko{z}}}^\inv\pba{n-1}{\ko{z}}\hp{2n-2}^\mpi\ek*{\rim\XFa{2n-1}{z}}\hp{2n-2}^\mpi\pda{n-1}{z}\ek*{\pda{n}{z}}^\inv\hp{2n}.
\]
\ecor
\bproof
 This follows by application of \rprop{M2140} with \(w=z\).
\eproof

\section{Description of the values of the solutions}\label{Cha13}
 In this section, ultimately we target an explicit description of the set
\beql{FR}
  \setaca*{F\rk{w}}{F\in\RFOqskg{2n}},
\eeq
 \tie{}, a parametrization of the set
\[%
  \setaca*{\int_\R\frac{1}{x-w}\sigma\rk{\dif x}}{\sigma\in\MggqRskg{2n}}
\]%
 with respect to an arbitrarily prescribed sequence \(\seqs{2n}\in\Hggeq{2n}\) and an arbitrarily chosen point \(w\in\pip\).
 As already mentioned, the set \eqref{FR} can be described  as matrix ball.
 In order to realize our aim, we are going to introduce a notation which will turn out to give possible parameters of the matrix ball.
 In order to justify that this notation is correct, we need a little preparation.
 Observe that the following constructions are well defined due to \rremss{SK32R}{R0735}:

\bnotal{N1326}
 Let \(\kappa\in\Ninf\) and let \(\seqska\in\Hggeqka\) with \tsHp{} \(\sHp{\kappa}\), \sabcd{} \(\abcd{\ev{\kappa}}\), and \tsXF{} \(\seqX{\kappa}\).
 For all \(n\in\NO\) such that \(2n+1\leq\kappa\) let \(\PA{n},\PB{n},\PC{n},\PD{n}\colon\CR\to\Cqq\) be defined by
\begin{align*}
 \PAa{n}{z}&\defeq\paa{n}{z}\hp{2n}^\mpi\ek*{\XFa{2n+1}{z}}^\ad-\paa{n+1}{z},&
 \PBa{n}{z}&\defeq\pba{n}{z}\hp{2n}^\mpi\ek*{\XFa{2n+1}{z}}^\ad-\pba{n+1}{z},\\
 \PCa{n}{z}&\defeq\ek*{\XFa{2n+1}{z}}^\ad\hp{2n}^\mpi\pca{n}{z}-\pca{n+1}{z},&
 \PDa{n}{z}&\defeq\ek*{\XFa{2n+1}{z}}^\ad \hp{2n}^\mpi\pda{n}{z}-\pda{n+1}{z}.
\end{align*}
\enota

 In view of \eqref{K19.0}, \eqref{K19.1}, \rlem{L0625}, and \eqref{XF.01}, for all \(z\in\CR\), we have
\begin{align}\label{ABCD0}
 \PAa{0}{z}&=-\hp{0},&
 \PBa{0}{z}&=\ko z\hp{0}^\mpi\hp{0}-z\Iq,&
 \PCa{0}{z}&=-\hp{0},&
 \PDa{0}{z}&=\ko z\hp{0}\hp{0}^\mpi-z\Iq.
\end{align}
 
\breml{R1352}
 Let \(\kappa\in\Ninf\) and let \(\seqska\in\Hggeqka\).
 Using \eqref{R1019.0}, \rremssss{A.R.A++*}{R1624}{BK8.2}{R0735}, and \rlem{L0625}, it is readily checked that \(\PCa{n}{z}=\ek{\PAa{n}{\ko{z}}}^\ad\) and \(\PDa{n}{z}=\ek{\PBa{n}{\ko{z}}}^\ad\) hold true for all \(n\in\NO\) fulfilling \(2n+1\leq\kappa\) and all \(z\in\CR\).
\erem
 
\bleml{L2412}
 Let \(\kappa\in\Ninf\), let \(\seqska\in\Hggeqka\), and let \(n\in\NO\) be such that \(2n+1\leq\kappa\).
 Then \(\det\PBa{n}{z}\neq0\) and \(\det\PDa{n}{z}\neq0\) for all \(z\in\CR\).
\elem
\bproof
 Let \(z\in\CR\).
 We consider an arbitrary \(x\in\nul{\PBa{n}{z}}\).
 Thus
\beql{L2412.1}
 \pba{n}{z}\hp{2n}^\mpi\ek{\XFa{2n+1}{z}}^\ad x
 =\pba{n+1}{z}x
\eeq
 \rlem{BK8.6} shows that
\begin{align}\label{L2412.2}
 \det\pba{n}{z}&\neq0&
&\text{and}&
 \det\pba{n+1}{z}&\neq0
\end{align}
 hold true.
 Because of \eqref{L2412.1} and \eqref{L2412.2}, we acquire
\beql{L2412.3}
 \hp{2n}^\mpi\ek{\XFa{2n+1}{z}}^\ad x
 =\ek*{\pba{n}{z}}^\inv\pba{n+1}{z}x.
\eeq
 Regarding \eqref{SK5-11}, \rprop{L1347} yields \(\ran{\ek{\XFa{2n+1}{z}}^\ad}=\ran{\hp{2n}}\).
 \rremp{tsa12}{tsa12.a} implies
\beql{L2412.4}
 \hp{2n}\hp{2n}^\mpi\ek*{\XFa{2n+1}{z}}^\ad
 =\ek*{\XFa{2n+1}{z}}^\ad.
\eeq
 By virtue of \rrem{R52SK}, we perceive \eqref{Z18}.
 Taking into account \eqref{Z18}, \eqref{L2412.3}, and \eqref{L2412.4}, we get
\[
 \XFa{2n+1}{z}x
 =\hp{2n}\ek*{\pba{n}{z}}^\inv\pba{n+1}{z}x
 =\hp{2n}\hp{2n}^\mpi\ek*{\XFa{2n+1}{z}}^\ad x
 =\ek*{\XFa{2n+1}{z}}^\ad x.
\]
 This implies \(\rk{\XFa{2n+1}{z}-\ek{\XFa{2n+1}{z}}^\ad}x=\Ouu{q}{1}\).
 Therefore, \(x\in\nul{\rim\XFa{2n+1}{z}}\).
 In view of \rprop{L1347}, we get \(x\in\nul{\ek{\XFa{2n+1}{z}}^\ad}\).
 Because of \eqref{L2412.1}, then \(\pba{n+1}{z}x=\Ouu{q}{1}\) follows.
 Thus, using \eqref{L2412.2}, we receive \(x=\Ouu{q}{1}\).
 Hence, \(\nul{\PBa{n}{z}}\subseteq\set{\Ouu{q}{1}}\).
 Consequently, \(\det\PBa{n}{z}\neq0\) for all \(z\in\CR\).
 Using \rrem{R1352}, we then conclude \(\det\PDa{n}{z}\neq0\) for all \(z\in\CR\) as well.
\eproof

\bleml{L1257}
 Let \(\kappa\in\Ninf\), let \(\seqska\in\Hggeqka\), and let \(n\in\NO\) be such that \(2n+1\leq\kappa\).
 Then \(\PDa{n}{z}\PAa{n}{z}=\PCa{n}{z}\PBa{n}{z}\) for all \(z\in\CR\).
\elem
\bproof
 Let \(z\in\CR\).
 According to \rlem{P8.14}, we acquire
\begin{align}
 \pda{n}{z}\paa{n}{z}-\pca{n}{z}\pba{n}{z}&=\Oqq,&\pda{n}{z}\paa{n+1}{z}-\pca{n}{z}\pba{n+1}{z}&=\hp{2n},\label{L1257.1}\\
 \pda{n+1}{z}\paa{n}{z}-\pca{n+1}{z}\pba{n}{z}&=-\hp{2n},&\pda{n+1}{z}\paa{n+1}{z}-\pca{n+1}{z}\pba{n+1}{z}&=\Oqq.\label{L1257.2}
\end{align}
 Regarding \eqref{SK5-11}, \rprop{L1347} yields \(\ran{\ek{\XFa{2n+1}{z}}^\ad}=\ran{\hp{2n}}\) and \(\nul{\ek{\XFa{2n+1}{z}}^\ad}=\nul{\hp{2n}}\).
 Using \rrem{tsa12}, then
\begin{align}\label{L1257.3}
 \hp{2n}\hp{2n}^\mpi\ek*{\XFa{2n+1}{z}}^\ad&=\ek*{\XFa{2n+1}{z}}^\ad&
&\text{and}&
 \ek*{\XFa{2n+1}{z}}^\ad\hp{2n}^\mpi\hp{2n}&=\ek*{\XFa{2n+1}{z}}^\ad
\end{align}
 follow.
 From \eqref{L1257.1}, \eqref{L1257.2}, and \eqref{L1257.3}, thus we get
\[\begin{split}
 &\rk*{\ek*{\XFa{2n+1}{z}}^\ad\hp{2n}^\mpi\pda{n}{z}-\pda{n+1}{z}}\rk*{\paa{n}{z}\hp{2n}^\mpi\ek*{\XFa{2n+1}{z}}^\ad-\paa{n+1}{z}}\\
 &\qquad-\rk*{\ek*{\XFa{2n+1}{z}}^\ad\hp{2n}^\mpi\pca{n}{z}-\pca{n+1}{z}}\rk*{\pba{n}{z}\hp{2n}^\mpi\ek*{\XFa{2n+1}{z}}^\ad-\pba{n+1}{z}}\\
 &=\ek*{\XFa{2n+1}{z}}^\ad\hp{2n}^\mpi\ek*{\pda{n}{z}\paa{n}{z}-\pca{n}{z}\pba{n}{z}}\hp{2n}^\mpi\ek*{\XFa{2n+1}{z}}^\ad\\
 &\qquad-\ek*{\XFa{2n+1}{z}}^\ad\hp{2n}^\mpi\ek{\pda{n}{z}\paa{n+1}{z}-\pca{n}{z}\pba{n+1}{z}}\\
 &\qquad-\ek*{\pda{n+1}{z}\paa{n}{z}-\pca{n+1}{z}\pba{n}{z}}\hp{2n}^\mpi\ek*{\XFa{2n+1}{z}}^\ad+\pda{n+1}{z}\paa{n+1}{z}-\pca{n+1}{z}\pba{n+1}{z}\\
 &=-\ek*{\XFa{2n+1}{z}}^\ad\hp{2n}^\mpi\hp{2n}-\rk{-\hp{2n}}\hp{2n}^\mpi\ek*{\XFa{2n+1}{z}}^\ad
 =\Oqq.
\end{split}\]
 In view of \rnota{N1326}, the proof is complete.
\eproof

\bleml{L1313}
 Let \(\kappa\in\NOinf\), let \(\seqska\in\Hggeqka\) with \tsHp{} \(\sHp{\kappa}\), \sabcd{} \(\abcd{\ev{\kappa}}\), and \tsXF{} \(\seqX{\kappa}\), and let \(n\in\NO\) be such that \(2n\leq\kappa\).
 For each \(z\in\CR\), then
\( 
 \det\rk{\ek{\XFa{2n}{z}}^\ad\hp{2n}^\mpi\pda{n}{z}-\pdoa{n+1}{z}}
 \neq0
\).
\elem
\bproof
 We consider an arbitrary \(z\in\CR\).
 \rprop{L1237-A} shows that the \tnatext{} \(\seqsh{2n+1}\) of \(\seqs{2n}\) belongs to \(\Hggeq{2n+1}\).
 Let \(\seqXh{2n+1}\) be the \tsXFo{\(\seqsh{2n+1}\)}.
 Then \rlem{L1520} provides \(\XFha{2n+1}{z}=\XFa{2n}{z}\).
 Let \(\shHp{2n+1}\) be the \tsHpo{\(\seqsh{2n+1}\)}.
 Then \rrem{R26-SK1} yields \(\hhp{2n}=\hp{2n}\).
 Let \(\habcd{n+1}\) be the \sabcdo{\(\seqsh{2n+1}\)}.
 Then \rlem{L1237-B} shows \(\pdh{n}=\pd{n}\) and \(\pdh{n+1}=\pdo{n+1}\).
 Consequently, we receive
\[
 \ek*{\XFa{2n}{z}}^\ad\hp{2n}^\mpi\pda{n}{z}-\pdoa{n+1}{z}
 =\ek*{\XFha{2n+1}{z}}^\ad \hhp{2n}^\mpi\pdha{n}{z}-\pdha{n+1}{z}. 
\]
 Regarding \rnota{N1326}, applying \rlem{L2412} to the sequence \(\seqsh{2n+1}\) completes the proof.
\eproof

 Now we introduce the central object of this section.
 Observe that the following constructions are well defined due to \rcor{C1232} and \rlemsss{BK8.6}{L1313}{L2412}, regarding \rnota{N1326}:

\bnotal{CB15}
 Let \(\kappa\in\NOinf\), let \(\seqska\in\Hggeqka\), let \(\sHp{\kappa}\) be the \tsHpo{\(\seqska\)}, let \(\abcd{\ev{\kappa}}\) be the \sabcdo{\(\seqska\)}, and let \(\seqX{\kappa}\) be the \tsXFo{\(\seqska\)}.
\benui
 \il{CB15.a} For each \(n\in\NO\) such that \(2n\leq \kappa\), let \(\lrk{2n},\rrk{2n},\mpk{2n}\colon\CR\to\Cqq\) be defined by
\begin{align*}
 \lrka{2n}{z}&\defeq\ek*{\pda{n}{z}}^\inv\hp{2n}\sqrt{\rk{\rim z}^\inv\rim\XFa{2n}{z}}^\mpi,\\
 \rrka{2n}{z}&\defeq\sqrt{\rk{\rim z}^\inv\rim\XFa{2n}{z}}^\mpi\hp{2n}\ek*{\pba{n}{z}}^\inv,
\intertext{and}
 \mpka{2n}{z}&\defeq-\rk*{\ek*{\XFa{2n}{z}}^\ad\hp{2n}^\mpi\pda{n}{z}-\pdoa{n+1}{z}}^\inv\rk*{\ek*{\XFa{2n}{z}}^\ad\hp{2n}^\mpi\pca{n}{z}-\pcoa{n+1}{z}}.
\end{align*}
 \il{CB15.b} For each \(n\in\NO\) such that \(2n+1\leq\kappa\), let \(\lrk{2n+1},\rrk{2n+1},\mpk{2n+1}\colon\CR\to\Cqq\) be given by
\begin{align*}
 \lrka{2n+1}{z}&\defeq\ek{\pda{n}{z}}^\inv\hp{2n}\sqrt{\rk{\rim z}^\inv\rim\XFa{2n+1}{z}}^\mpi,\\
 \rrka{2n+1}{z}&\defeq\sqrt{\rk{\rim z}^\inv\rim\XFa{2n+1}{z}}^\mpi\hp{2n}\ek{\pba{n}{z}}^\inv,
\intertext{and}
 \mpka{2n+1}{z}&\defeq-\rk*{\ek*{\XFa{2n+1}{z}}^\ad\hp{2n}^\mpi\pda{n}{z}-\pda{n+1}{z}}^\inv\rk*{\ek*{\XFa{2n+1}{z}}^\ad\hp{2n}^\mpi\pca{n}{z}-\pca{n+1}{z}}.
\end{align*}
\eenui
\enota

 Recall that \(\Kpq\) stands for the set of all contractive complex \tpqa{matrices}.
 The set
\[
 \cmb{M}{A}{B}
 \defeq\setaca{M+AKB}{K\in\Kpq}
\]
 signifies the (closed) matrix ball with \emph{center} \(M\), \emph{left semi-radius} \(A\), and \emph{right semi-radius} \(B\) with respect to given matrices \(M\in\Cpq\), \(A\in\Cpp\), and \(B\in\Cqq\).
 The ambient theory dates back to Yu.~L.~Shmul\cprime yan \cite{MR0273377}, who, moreover, examined the operator case in the context of Hilbert spaces.
 Observe that the particular case of matrices is elaborated in \zitaa{MR1152328}{\csec{1.5}{44--51}}.
 Using \rnotap{CB15}{CB15.a}, now we are able to formulate the central result of this paper, which is a generalization of Kovalishina's result \zitaa{MR703593}{\S2}, who studied the non-degenerate case that the given sequence \(\seqs{2n}\) belongs to \(\Hgq{2n}\).

\bthml{T45CD}
 Let \(n\in\NO\) and let \(\seqs{2n}\in\Hggeq{2n}\).
 For each \(w\in\pip\), then
\[
 \setaca*{F\rk{w}}{F\in\RFOqskg{2n}}
 =\Cmb{\mpka{2n}{w}}{\rk{w-\ko w}^\inv\lrka{2n}{w}}{\rrka{2n}{w}},
\]
 where \(\lrk{2n}\), \(\rrk{2n}\), and \(\mpk{2n}\) are given by \rnotap{CB15}{CB15.a}.
\ethm

 The main goal of this section is to state a proof of \rthm{T45CD}.
 Before doing this, we look at \rnota{CB15}. %
 At first view it looks a little bit surprising why the \rnota{CB15} is introduced for the two different indices \(2n\) and \(2n+1\), because \rnota{CB15} will be applied later to describe the Weyl matrix ball associated with a sequence  \(\seqs{2n}\in\Hggeq{2n}\).
 It turns out soon (see \rlemss{C220N}{L3N10} below) that the corresponding matrices introduced in \rpartss{CB15.a}{CB15.b} of \rnota{CB15} coincide.
 The final steps for the proof of \rthm{T45CD} are \rpropss{L1204}{L1301}.
 They are formulated in terms of \rnotap{CB15}{CB15.a}.
 However, in the corresponding proofs we use the matrices in the form expressed by \rnotap{CB15}{CB15.b}.
 The strategy of our proof of \rthm{T45CD} is in a wider sense inspired by the proof of \zitaa{MR2656833}{\cthm{1.1}{}}, which contains the description of Weyl matrix balls associated with a finite \tqqa{Carath\'eodory} sequence.
 Against this background the terminology for the parameters of the Weyl matrix balls was already chosen in a similar way as in \cite{MR2656833}.
 In both cases, the first obstacle was to guess the right conjecture about the parameters of the Weyl matrix ball we are striving for.
 As in the Carath\'eodory case the conjectured parameters of the Weyl matrix balls are essentially determined by a rational \tqqa{matrix}-valued function and some constant matrices.
 In our case, this corresponds to the rational \tqqa{matrix}-valued function \(\XF{2n}\) introduced in \rdefn{K17.1} and the \tHp{} \(\hp{2n}\) given by \rdefn{CM2.1.64}.
 In the following, we will use the terms listed in \rnota{CB15} without always referring to their definition.

\breml{ABM-tru}
 Let \(\kappa\in\NOinf\) and let \(\seqska\in\Hggeqka\).
 Regarding \rnotass{CB15}{N-abcdO}, then one can see from \rremsss{CM2.1.65}{BK8.1}{XF-tru} that for each \(m\in\mn{0}{\kappa}\), the functions \(\lrk{m}\), \(\rrk{m}\), and \(\mpk{m}\) are built only from the matrices \(\su{0},\su{1},\dotsc,\su{m}\) and does not depend on the matrices \(\su{j}\) with \(j\geq m+1\).
 Moreover, \rlemss{C220N}{L3N10} below even show that \(\lrk{m}\), \(\rrk{m}\), and \(\mpk{m}\) are independent of \(\su{j}\) with \(j\geq\eff{m}+1\), where \(\eff{m}\) is given by \eqref{SK5-11}.
\erem
 
\bleml{L0810}
 Let \(\kappa\in\NOinf\) and let \(\seqska\in\Hggeqka\), and let \(z\in\CR\).
 Then
\begin{align*}
 \lrka{0}{z}&=\sqrt{\su{0}},&
 \rrka{0}{z}&=\sqrt{\su{0}},&
&\text{and}&
 \mpka{0}{z}&=\rk{\ko z-z}^\inv\su{0}.
\end{align*}
\elem
\bproof
 First observe that \(\su{0}\in\Cggq\) is valid, because of \(\seqska\in\Hggeqka\) and that \(\ko z-z\neq0\) holds true, due to \(z\in\CR\).
 According to \eqref{Hp.01}, we have \(\hp{0}=\su{0}\).
 By virtue of \eqref{XF.01}, furthermore \(\XF{0}\rk{z}=z\su{0}\).
 Consequently, \(\ek{\XFa{0}{z}}^\ad=\ko z\su{0}\) and \(\rim\XF{0}\rk{z}=\rim\rk{z}\su{0}\).
 In particular, \(\rk{\rim z}^\inv\rim\XF{0}\rk{z}=\su{0}\).
 \rrem{A.R.A+sqrt} yields \(\sqrt{\su{0}}^\mpi\su{0}=\sqrt{\su{0}}\) and \(\su{0}\sqrt{\su{0}}^\mpi=\sqrt{\su{0}}\).
 Thus, taking additionally into account \rnotap{CB15}{CB15.a} and \eqref{K19.0}, we perceive \(\lrka{0}{z}=\ek{\pda{0}{z}}^\inv\hp{0}\sqrt{\rk{\rim z}^\inv\rim\XF{0}\rk{z}}^\mpi=\su{0}\sqrt{\su{0}}^\mpi=\sqrt{\su{0}}\) and \(\rrka{0}{z}=\sqrt{\rk{\rim z}^\inv\rim\XFa{0}{z}}^\mpi\hp{0}\ek{\pba{0}{z}}^\inv=\sqrt{\su{0}}^\mpi\su{0}=\sqrt{\su{0}}\).
 Regarding \eqref{K19.0} and \eqref{abcdO-1}, we get furthermore
\beql{L0810.1}
 \ek*{\XFa{0}{z}}^\ad\hp{0}^\mpi\pda{0}{z}-\pdoa{1}{z}
 =\ko z\su{0}\hp{0}^\mpi\cdot\Iq-z\Iq
 =\ko z\su{0}\su{0}^\mpi-z\Iq
\eeq
 and
\beql{L0810.2}
 \ek*{\XFa{0}{z}}^\ad\hp{0}^\mpi\pca{0}{z}-\pcoa{1}{z}
 =\ko z\su{0}\hp{0}^\mpi\cdot\Oqq-\hp{0}
 =-\su{0}.
\eeq
 \rlem{L1313} yields that the matrix on the left-hand side of \eqref{L0810.1} is non-singular.
 Consequently, the matrix on the right-hand side of \eqref{L0810.1} is non-singular as well.
 Since \(\rk{\ko z\su{0}\su{0}^\mpi-z\Iq}\su{0}=\rk{\ko z-z}\su{0}\) holds true, we thus conclude \(\rk{\ko z-z}^\inv\su{0}=\rk{\ko z\su{0}\su{0}^\mpi-z\Iq}^\inv\su{0}\).
 Using additionally \rnotap{CB15}{CB15.a}, \eqref{L0810.1}, and \eqref{L0810.2}, we get finally \(\mpka{0}{z}=-\rk{\ko z\su{0}\su{0}^\mpi-z\Iq}^\inv\rk{-\su{0}}={\rk{\ko z-z}}^\inv\su{0}\).
\eproof

 Now we are going to verify the announced statements containing the coincidence of the corresponding matrices introduced in \rpartss{CB15.a}{CB15.b} of \rnota{CB15}, respectively.
 
\bleml{C220N}
 Let \(\kappa\in\Ninf \), let \(\seqska\in\Hggeqka\), let \(n\in\NO\) be such that \(2n+1\leq\kappa\), and let \(z\in\CR\).
 Then \(\lrka{2n}{z}=\lrka{2n+1}{z}\) and \(\rrka{2n}{z}=\rrka{2n+1}{z}\).
\elem
\bproof
 Let \(\sHp{\kappa}\) be the \tsHpo{\(\seqska\)}.
 From \rrem{R1019} we can infer \(\rim\hp{2n+1}=\Oqq\).
 Using \rlem{K17-2}, then we get \(\rim\XFa{2n}{z}=\rim\XFa{2n+1}{z}\).
 In view of \rnota{CB15}, the proof is complete.
\eproof

\bleml{L0819}
 Let \(\kappa\in\Ninf\), let \(\seqska\in\Hggeqka\), and let \(n\in\NO\) be such that \(2n+1\leq\kappa\).
 Let \(\PA{n},\PB{n},\PC{n},\PD{n}\) be given by \rnota{N1326}.
 For all \(z\in\CR\), then the matrices \(\PBa{n}{z}\) and \(\PDa{n}{z}\) are invertible and, furthermore, \(\mpka{2n+1}{z}=-\PAa{n}{z}\ek{\PBa{n}{z}}^\inv\) as well as \(\mpka{2n+1}{z}=-\ek{\PDa{n}{z}}^\inv\PCa{n}{z}\) hold true.
\elem
\bproof
 Let \(z\in\CR\).
 \rlem{L2412} shows that the matrices \(\PBa{n}{z}\) and \(\PDa{n}{z}\) are invertible.
 Because of \rlem{L1257}, we have \(\PDa{n}{z}\PAa{n}{z}=\PCa{n}{z}\PBa{n}{z}\).
 Consequently, we conclude \(-\PAa{n}{z}\ek{\PBa{n}{z}}^\inv=-\ek{\PDa{n}{z}}^\inv\PCa{n}{z}\).
 By virtue of \rnotass{CB15}{N1326}, we furthermore see \(\mpka{2n+1}{z}=-\ek{\PDa{n}{z}}^\inv\PCa{n}{z}\), which completes the proof.
\eproof

\bleml{L3N10}
 Let \(\kappa\in\Ninf\), let \(\seqska\in\Hggeqka\), let \(n\in\NO\) be such that \(2n+1\leq\kappa\), and let \(z\in\CR\).
 Then \(\mpka{2n}{z}=\mpka{2n+1}{z}\).
\elem
\bproof
 Let \(\sHp{\kappa}\) be the \tsHpo{\(\seqska\)}, let \(\abcd{\ev{\kappa}}\) be the \sabcdo{\(\seqska\)}, and let \(\seqX{\kappa}\) be the \tsXFo{\(\seqska\)}.
 Let \(\PCoa{n}{z}\defeq\ek{\XFa{2n}{z}}^\ad\hp{2n}^\mpi\pca{n}{z}-\pcoa{n+1}{z}\) and \(\PDoa{n}{z}\defeq\ek{\XFa{2n}{z}}^\ad\hp{2n}^\mpi\pda{n}{z}-\pdoa{n+1}{z}\).
 According to \rlem{L1313}, we infer \(\det\PDoa{n}{z}\neq0\), whereas \rlem{L2412} yields \(\det\PDa{n}{z}\neq0\).
 We see from \rnotass{CB15}{N1326} that
\begin{align}\label{IWN}
 \mpka{2n}{z}&=-\ek*{\PDoa{n}{z}}^\inv\PCoa{n}{z}&
 &\text{and}&
 \mpka{2n+1}{z}&=-\ek*{\PDa{n}{z}}^\inv\PCa{n}{z} 
\end{align}
 hold true.
 \rlem{K17-2} provides \(\XFa{2n}{z}-\XFa{2n+1}{z}=\hp{2n+1}\).
 Since \rrem{R1019} yields \(\hp{2n+1}^\ad =\hp{2n+1}\), consequently, \(\ek{\XFa{2n}{z}}^\ad-\ek{\XFa{2n+1}{z}}^\ad=\hp{2n+1}\) follows.
 Using \rrem{R1358}, we can furthermore infer that
\begin{align*}
 \pcoa{n+1}{z}-\pca{n+1}{z}&=\hp{2n+1}\hp{2n}^\mpi\pca{n}{z}&
&\text{and}&
 \pdoa{n+1}{z}-\pda{n+1}{z}&=\hp{2n+1}\hp{2n}^\mpi\pda{n}{z}
\end{align*}
 are valid. 
 Regarding \rnota{N1326}, hence we obtain
\[%
 \PCoa{n}{z}-\PCa{n}{z} 
 =\rk*{\ek*{\XFa{2n}{z}}^\ad-\ek*{\XFa{2n+1}{z}}^\ad}\hp{2n}^\mpi\pca{n}{z}-\ek*{\pcoa{n+1}{z}-\pca{n+1}{z}}
 =\Oqq
\]%
 and, analogously, \(\PDoa{n}{z}-\PDa{n}{z}=\Oqq.\)
 In view of \eqref{IWN}, the proof is then complete.
\eproof

 We continue providing further statements that address essential properties of the objects defined in \rnota{CB15}.
 Moreover, these results prepare to embracing the capability of establishing the desired parametrization.

\bleml{CL18}
 Let \(\kappa\in\NOinf\) and let \(\seqska\in\Hggeqka\).
 For each \(m\in\mn{0}{\kappa}\), then the matrix-valued functions \(\lrk{m},\rrk{m},\mpk{m}\colon\CR\to\Cqq\) given by \rnota{CB15} are continuous.
\elem
\bproof
 Let \(m\in\mn{0}{\kappa}\) and let \(\seqX{\kappa}\) be the \tsXFo{\(\seqska\)}.
 Using \rrem{R0735}, we can infer that \(f,g\colon\CR\to\Cqq\) defined by \(f\rk{z}\defeq\ek{\XFa{m}{z}}^\ad\) and \(g\rk{z}\defeq\rk{\rim z}^\inv\rim\XFa{m}{z}\), respectively, both are continuous.
 According to \rcor{C1232}, we have furthermore \(g\rk{z}\in\Cggq\) for all \(z\in\CR\).
 Because of \rrem{ZR2}, then \(r\colon\CR\to\Cggq\) defined by \(r\rk{z}\defeq\sqrt{g\rk{z}}\) is continuous.
 Let \(\sHp{\kappa}\) be the \tsHpo{\(\seqska\)}.
 From \rprop{L1347} we can conclude \(\ran{g\rk{z}}=\ran{\hp{\eff{m}}}\) for all \(z\in\CR\).
 In view of \rrem{A.R.r-sqrt}, then \(\ran{r\rk{z}}=\ran{\hp{\eff{m}}}\) for all \(z\in\CR\) follows.
 In particular, \(\rank r\rk{z}=\rank\hp{\eff{m}}=\rank r\rk{w}\) for every choice of \(z\) and \(w\) in \(\CR\).
 Consequently, using \rlem{LemC1}, we receive that \(\rho\colon\CR\to\Cqq\) given by \(\rho\rk{z}\defeq\ek{r\rk{z}}^\mpi\) is continuous.
 For our following considerations, we note that \(m=2n\) or \(m=2n+1\), where \(n\defeq\ef{m}\) is given by \eqref{SK5-11}.
 We first show that \(\lrk{m}\) and \(\rrk{m}\) are continuous.
 \rrem{SK32R} and \rlem{BK8.6} provide that \(\pb{n}\) and \(\pd{n}\) are matrix polynomials which fulfill \(\det\pba{n}{z}\neq0\) and \(\det\pda{n}{z}\neq0\) for all \(z\in\CR\).
 Therefore, \rlem{LemC1} yields that \(\eta,\theta\colon\CR\to\Cqq\) defined by \(\eta\rk{z}\defeq\ek{\pba{n}{z}}^\inv\) and \(\theta\rk{z}\defeq\ek{\pda{n}{z}}^\inv\), respectively, both are continuous.
 Taking into account that \rnota{CB15} implies \(\lrka{m}{z}=\theta\rk{z}\hp{2n}\rho\rk{z}\) and \(\rrka{m}{z}=\rho\rk{z}\hp{2n}\eta\rk{z}\)
 for all \(z\in\CR\), we obtain that both \(\lrk{m}\) and \(\rrk{m}\) are continuous.
 We now show that \(\mpk{m}\) is continuous.
 Consider the case \(m=2n\).
 Since \(f\) is continuous, we conclude from \rrem{SK32R} and \rnota{N-abcdO} that \(c,d\colon\CR\to\Cqq\) defined by \(c\rk{z}\defeq f\rk{z}\hp{2n}^\mpi\pca{n}{z}-\pcoa{n+1}{z}\) and \(d\rk{z}\defeq f\rk{z}\hp{2n}^\mpi\pda{n}{z}-\pdoa{n+1}{z}\), respectively, both are continuous.
 \rlem{L1313} shows that \(\det d\rk{z}\neq0\) holds true for all \(z\in\CR\).
 By virtue of \rlem{LemC1}, then \(\delta\colon\CR\to\Cqq\) given by \(\delta\rk{z}\defeq\ek{d\rk{z}}^\inv\) is continuous as well.
 Since we realize from \rnotap{CB15}{CB15.a} that \(\mpk{2n}=-\delta c\) is valid, thus we obtain that \(\mpk{2n}\) is continuous.
 Using \rlem{L2412} instead of \rlem{L1313}, the case \(m=2n+1\) can be treated analogously.
\eproof

 Now we verify that, for each \(z\in\CR\), the spaces \(\nul{\lrka{m}{z}}\) and \(\ran{\rrka{m}{z}}\) are completely determined by \(\hp{\eff{m}}\) and, consequently, independent of \(z\).

\bpropl{CL19}
 Let \(\kappa\in\NOinf\), let \(\seqska\in\Hggeqka\) with \tsHp{} \(\sHp{\kappa}\), and let \(z\in\CR\).
 For all \(m\in\mn{0}{\kappa}\), then
\begin{align}\label{BG}
 \nul{\lrka{m}{z}}&=\nul{\hp{\eff{m}}}&
&\text{as well as}&
 \ran{\rrka{m}{z}}&=\ran{\hp{\eff{m}}}
\end{align}
 and, in particular, \(\rank\lrka{m}{z}=\rank\hp{\eff{m}}\) and \(\rank\rrka{m}{z}=\rank\hp{\eff{m}}\), where \(\lrk{m}\) and \(\rrk{m}\) are given by \rnota{CB15} and where \(\eff{m}\) is given by \eqref{SK5-11}.
\eprop
\bproof
 Our proof is divided into four parts:
 
 (I) First we consider the case that \(m=2n+1\) with some \(n\in\NO\).
 According to \eqref{SK5-11}, then \(\eff{m}=2n\).
 We use the notation introduced in the proof of \rlem{CL18}.
 Regarding \rnotap{CB15}{CB15.b} and \rlem{BK8.6}, then \rrem{R1241} shows that
\begin{align}
 \Nul{\lrka{2n+1}{z}}&=\Nul{\ek*{\pda{n}{z}}^\inv\hp{2n}\rho\rk{z}}=\Nul{\hp{2n}\rho\rk{z}}\label{N9V}\\
\intertext{and}
 \Ran{\rrka{2n+1}{z}}&=\Ran{\rho\rk{z}\hp{2n}\ek*{\pba{n}{z}}^\inv}=\Ran{\rho\rk{z}\hp{2n}}\label{M2V}
\end{align}
 are valid.
 According to the proof of \rlem{CL18} and the notation therein, we also get \(\ek{r\rk{z}}^\ad=r\rk{z}\) and \(\ran{r\rk{z}}=\ran{\hp{\eff{m}}}\).
 Hence, using \(\eff{m}=2n\) and \rrem{A.R.r-mpi}, we receive
\begin{align}
 \Ran{\rho\rk{z}}
 =\Ran{\ek*{r\rk{z}}^\mpi}
 =\Ran{\ek*{r\rk{z}}^\ad}
 =\Ran{r\rk{z}}
 =\ran{\hp{2n}}\label{LBM}
\intertext{and}
 \Nul{\rho\rk{z}}
 =\Nul{\ek*{r\rk{z}}^\mpi}
 =\Nul{\ek*{r\rk{z}}^\ad}
 =\Nul{r\rk{z}}
 =\nul{\hp{2n}}.\label{NVB}
\end{align}
 Since \rrem{R1019} yields \(\hp{2n}^\ad=\hp{2n}\), in view of \eqref{LBM} and \eqref{NVB}, then \rlem{C2} provides \(\nul{\hp{2n}\rho\rk{z}}=\nul{\rho\rk{z}}\) and \(\ran{\rho\rk{z}\hp{2n}}=\ran{\rho\rk{z}}\).
 Consequently, in view of \eqref{N9V} and \eqref{NVB}, we obtain \(\nul{\lrka{2n+1}{z}}=\nul{\hp{2n}\rho\rk{z}}=\nul{\rho\rk{z}}=\nul{\hp{2n}}\) and, because of \eqref{M2V} and \eqref{LBM}, analogously, \(\ran{\rrka{2n+1}{z}}=\ran{\hp{2n}}\).
 Regarding \eqref{SK5-11}, thus \eqref{BG} is proved in the case that \(m\) is a positive odd integer.
 
 (II) Consider now the case that \(m=2n\) with some \(n\in\NO\).
 \rprop{L1237-A} shows that the \tnatext{} \(\seqsh{2n+1}\) of \(\seqs{2n}\) belongs to \(\Hggeq{2n+1}\).
 By virtue of \rdefn{D.nat-ext} and \rrem{ABM-tru}, we have \(\lrhka{2n}{z}=\lrka{2n}{z}\) and \(\rrhka{2n}{z}=\rrka{2n}{z}\), where \(\lrhk{2n}\) and \(\rrhk{2n}\) are built according to \rnotap{CB15}{CB15.a} from the sequence \(\seqsh{2n+1}\).
 \rlem{C220N} provides \(\lrhka{2n}{z}=\lrhka{2n+1}{z}\) and \(\rrhka{2n}{z}=\rrhka{2n+1}{z}\), where \(\lrhk{2n+1}\) and \(\rrhk{2n+1}\) are built according to \rnotap{CB15}{CB15.b} from the sequence \(\seqsh{2n+1}\).
 Let \(\shHp{2n+1}\) be the \tsHpo{\(\seqsh{2n+1}\)}.
 According to \rrem{R26-SK1}, then \(\hhp{2n}=\hp{2n}\).
 Using \cpart{(I)} of the proof, consequently
\[ 
 \Nul{\lrka{2n}{z}}
 =\Nul{\lrhka{2n}{z}}
 =\Nul{\lrhka{2n+1}{z}}
 =\nul{\hhp{2n}}
 =\nul{\hp{2n}}
\]
 and, analogously, \(\ran{\rrka{2n}{z}}=\ran{\hp{2n}}\).
 Regarding \eqref{SK5-11}, thus we realize that \eqref{BG} is proved in the case that \(m\) is a non-negative even integer as well.
 
 (III) Summarizing parts~(I) and~(II), we see that \eqref{BG} is proved for all \(m\in\mn{0}{\kappa}\).
 
 (IV) For all \(m\in\mn{0}{\kappa}\), from \eqref{BG} we conclude finally \(\rank\lrka{m}{z}=q-\dim\nul{\lrka{m}{z}}=q-\dim\nul{\hp{\eff{m}}}=\rank\hp{\eff{m}}\) and, obviously, \(\rank\rrka{m}{z}=\rank\hp{\eff{m}}\).
\eproof

 Now we state some technical auxiliary results.

\bleml{L1433}
 Let \(n\in\NO\), let \(\seqs{2n+1}\in\Hggeq{2n+1}\), and let \(z\in\CR\).
 Then the matrix \(\PDa{n}{z}\) is invertible, the matrix \(\rk{\rim z}^\inv\rim\XFa{2n+1}{z}\) is \tnnH{}, and the equation
\beql{L1433.A}
 \rk{\ko z-z}\ek*{\PDa{n}{z}}^\inv\hp{2n}\hp{2n}^\mpi
 =\lrka{2n+1}{z}\sqrt{\rk{\rim z}^\inv\rim\XFa{2n+1}{z}}^\mpi
\eeq
 holds true, where \(\PD{n}\) is given by \rnota{N1326} and \(\lrk{2n+1}\) is given by \rnotap{CB15}{CB15.b}.
\elem
\bproof
 From \rlem{L2412} we see that \(\det\PDa{n}{z}\neq0\).
 By virtue of \rcor{C1232}, we acquire that \(B\defeq\rk{\rim z}^\inv\rim\XFa{2n+1}{z}\) satisfies \(B\in\Cggq\).
 According to \rpropp{K17-5}{K17-5.a}, we have \(\det\pda{n}{z}\neq0\) and \(\XFa{2n+1}{z}=\pda{n+1}{z}\ek{\pda{n}{z}}^\inv\hp{2n}\).
 \rprop{L1347} yields \(\nul{\ek{\XFa{2n+1}{z}}^\ad}=\nul{\hp{2n}}\), which, in view of \rremp{tsa12}{tsa12.b}, implies \(\ek{\XFa{2n+1}{z}}^\ad\hp{2n}^\mpi\hp{2n}=\ek{\XFa{2n+1}{z}}^\ad\).
 Taking additionally into account \rnota{N1326}, we hence get
\[\begin{split}
 \PDa{n}{z}\ek*{\pda{n}{z}}^\inv\hp{2n}
 &=\ek*{\XFa{2n+1}{z}}^\ad\hp{2n}^\mpi\hp{2n}-\pda{n+1}{z}\ek{\pda{n}{z}}^\inv\hp{2n}\\
 &=\ek*{\XFa{2n+1}{z}}^\ad-\XFa{2n+1}{z}
 =-2\iu\rim\rk*{\XFa{2n+1}{z}}
 =-2\iu\rim\rk{z}B
 =\rk{\ko z-z}B
\end{split}\]
 and, consequently, 
\beql{FR3}
 \ek*{\pda{n}{z}}^\inv\hp{2n}
 =\rk{\ko z-z}\ek*{\PDa{n}{z}}^\inv B.
\eeq
 By virtue of \rprop{L1347}, we acquire \(\ran{B}=\ran{\hp{2n}}\).
 Thus, from \rrem{tsb3} we can conclude \(BB^\mpi=\hp{2n}\hp{2n}^\mpi\).
 Additionally, using \eqref{FR3}, \rrem{A.R.A+>}, and \rnotap{CB15}{CB15.b}, we obtain finally
\[\begin{split}
 \rk{\ko z-z}\ek*{\PDa{n}{z}}^\inv\hp{2n}\hp{2n}^\mpi
 &=\rk{\ko z-z}\ek*{\PDa{n}{z}}^\inv BB^\mpi
 =\ek*{\pda{n}{z}}^\inv\hp{2n}B^\mpi\\
 &=\ek*{\pda{n}{z}}^\inv\hp{2n}\sqrt{B}^\mpi\sqrt{B}^\mpi
 =\lrka{2n+1}{z}\sqrt{B}^\mpi.\qedhere
\end{split}\]
\eproof

\bleml{L0852}
 Let \(n\in\NO\), let \(\seqs{2n+1}\in\Hggeq{2n+1}\), and let \(z\in\CR\).
 Let \(\cpP\in\Cqq\) be such that \(\ran{\cpP}\subseteq\ran{\hp{2n}}\) and let \(\cpQ\in\Cqq\).
 Furthermore, let
\begin{align}
 \STA{2n+1}&\defeq\paa{n}{z}\hp{2n}^\mpi \cpP+\paa{n+1}{z}\cpQ&
&\text{and}&
 \STB{2n+1}&\defeq\pba{n}{z}\hp{2n}^\mpi \cpP+\pba{n+1}{z}\cpQ.\label{L0852.AB}
\end{align}
 Then the matrices \(\PDa{n}{z}\) and \(\pba{n}{z}\) are invertible, the matrix \(\rk{\rim z}^\inv\rim\XFa{2n+1}{z}\) is \tnnH{}, and the identities
\beql{L0852.C}
 \rk{z-\ko z}\rk*{\ek*{\PDa{n}{z}}^\inv\PCa{n}{z}\STB{2n+1}-\STA{2n+1}}
 =\lrka{2n+1}{z}\sqrt{\rk{\rim z}^\inv\rim\XFa{2n+1}{z}}^\mpi\rk*{\cpP+\ek*{\XFa{2n+1}{z}}^\ad \cpQ}
\eeq
 and 
\[
 \hp{2n}\ek*{\pba{n}{z}}^\inv\STB{2n+1}
 =\cpP+\XFa{2n+1}{z}\cpQ
\]
 hold true, where \(\PC{n}\) and \(\PD{n}\) are given by \rnota{N1326} and \(\lrk{2n+1}\) is given by \rnotap{CB15}{CB15.b}.
\elem
\bproof
 From \rlem{L1433} we realize that \(\PDa{n}{z}\) is invertible, that \(\rk{\rim z}^\inv\rim\XFa{2n+1}{z}\) is \tnnH{}, and that \eqref{L1433.A} holds true.
 By virtue \rlem{P8.14}, we receive that
\begin{align}\label{L0852.1}
 \pda{n}{z}\paa{n}{z}&=\pca{n}{z}\pba{n}{z}&
&\text{and}&
 \pda{n+1}{z}\paa{n+1}{z}&=\pca{n+1}{z}\pba{n+1}{z},%
\end{align}
 as well as
\begin{align}\label{L0852.2}
 \pda{n}{z}\paa{n+1}{z}-\pca{n}{z}\pba{n+1}{z}&=\hp{2n}&%
&\text{and}&
 \pda{n+1}{z}\paa{n}{z}-\pca{n+1}{z}\pba{n}{z}&=-\hp{2n}%
\end{align}
 hold true.
 \rprop{L1347} yields \(\ran{\ek{\XFa{2n+1}{z}}^\ad}=\ran{\hp{2n}}\) and \(\nul{\ek{\XFa{2n+1}{z}}^\ad}=\nul{\hp{2n}}\).
 Thus, in view of \rrem{tsa12}, we obtain
\begin{align}\label{L0852.3}%
 \hp{2n}\hp{2n}^\mpi\ek*{\XFa{2n+1}{z}}^\ad&=\ek*{\XFa{2n+1}{z}}^\ad&
&\text{and}&
 \ek*{\XFa{2n+1}{z}}^\ad\hp{2n}^\mpi\hp{2n}&=\ek*{\XFa{2n+1}{z}}^\ad.
\end{align}
 Using \rnota{N1326}, \eqref{L0852.AB}, \eqref{L0852.1}, \eqref{L0852.2}, and \eqref{L0852.3}, we perceive
\beql{L0852.4}\begin{split}%
 \PCa{n}{z}\STB{2n+1}-\PDa{n}{z}\STA{2n+1}
 &=\rk*{\ek*{\XFa{2n+1}{z}}^\ad\hp{2n}^\mpi\pca{n}{z}-\pca{n+1}{z}}\rk*{\pba{n}{z}\hp{2n}^\mpi \cpP+\pba{n+1}{z}\cpQ}\\
 &\qquad-\rk*{\ek*{\XFa{2n+1}{z}}^\ad\hp{2n}^\mpi\pda{n}{z}-\pda{n+1}{z}}\rk*{\paa{n}{z}\hp{2n}^\mpi \cpP+\paa{n+1}{z}\cpQ}\\
 &=\ek*{\XFa{2n+1}{z}}^\ad\hp{2n}^\mpi\ek*{\pca{n}{z}\pba{n}{z}-\pda{n}{z}\paa{n}{z}}\hp{2n}^\mpi \cpP\\
 &\qquad+\ek*{\XFa{2n+1}{z}}^\ad\hp{2n}^\mpi\ek*{\pca{n}{z}\pba{n+1}{z}-\pda{n}{z}\paa{n+1}{z}}\cpQ\\
 &\qquad-\ek*{\pca{n+1}{z}\pba{n}{z}-\pda{n+1}{z}\paa{n}{z}}\hp{2n}^\mpi \cpP\\
 &\qquad-\ek*{\pca{n+1}{z}\pba{n+1}{z}-\pda{n+1}{z}\paa{n+1}{z}}\cpQ\\
 &=-\ek*{\XFa{2n+1}{z}}^\ad\hp{2n}^\mpi\hp{2n}\cpQ-\hp{2n}\hp{2n}^\mpi \cpP
 =-\ek*{\XFa{2n+1}{z}}^\ad \cpQ-\hp{2n}\hp{2n}^\mpi \cpP\\
 &=-\hp{2n}\hp{2n}^\mpi\ek*{\XFa{2n+1}{z}}^\ad \cpQ-\hp{2n}\hp{2n}^\mpi \cpP
 =-\hp{2n}\hp{2n}^\mpi\rk*{\cpP+\ek*{\XFa{2n+1}{z}}^\ad \cpQ}.
\end{split}\eeq
 Because of \(\ek{\PDa{n}{z}}^\inv\PCa{n}{z}\STB{2n+1}-\STA{2n+1}=\ek{\PDa{n}{z}}^\inv\ek{\PCa{n}{z}\STB{2n+1}-\PDa{n}{z}\STA{2n+1}}\)   and \eqref{L0852.4} as well as \eqref{L1433.A}, we conclude that \eqref{L0852.C} holds true.
 By virtue of the assumption \(\ran{\cpP}\subseteq\ran{\hp{2n}}\) and \rremp{tsa12}{tsa12.a}, we acquire \(\hp{2n}\hp{2n}^\mpi \cpP=\cpP\).
 According to \rremp{R52SK}{R52SK.a}, we have \(\det\pba{n}{z}\neq0\) and \(\XFa{2n+1}{z}=\hp{2n}\ek{\pba{n}{z}}^\inv\pba{n+1}{z}\).
 Regarding \eqref{L0852.AB}, then it follows
\[
 \hp{2n}\ek*{\pba{n}{z}}^\inv\STB{2n+1} 
 =\hp{2n}\hp{2n}^\mpi \cpP+\hp{2n}\ek*{\pba{n}{z}}^\inv\pba{n+1}{z}\cpQ
 =\cpP+\XFa{2n+1}{z}\cpQ.\qedhere
\]
\eproof

\bleml{L0827}
 Let \(n\in\NO\) and let \(\seqs{2n+1}\in\Hggeq{2n+1}\).
 Let \(P\defeq\OPu{\ran{\hp{2n}}}\) and \(Q\defeq\OPu{\nul{\hp{2n}}}\).
 Let \(C\in\Cqq\), let \(z\in\CR\), and let \(E\defeq-\ek{\XFa{2n+1}{z}}^\ad\) and \(B\defeq\rk{\rim z}^\inv\rim E\).
 Then \(B\in\Cggq\) and the matrices \(\cpP\defeq E\sqrt{B}^\mpi-E^\ad\sqrt{B}^\mpi CP\) and \(\cpQ\defeq\sqrt{B}^\mpi-\sqrt{B}^\mpi CP+Q\) fulfill \(\det\rk{\pba{n}{z}\hp{2n}^\mpi \cpP +\pba{n+1}{z}\cpQ}\neq0\).
\elem
\bproof
 Clearly, \(B=\rk{\rim z}^\inv\rim\XFa{2n+1}{z}\).
 \rcor{C1232} then shows that \(B\in\Cggq\).
 Obviously, we have \(-E^\ad=\XFa{2n+1}{z}\).
 Therefore, \(\cpP=\XFa{2n+1}{z}\sqrt{B}^\mpi CP-\ek{\XFa{2n+1}{z}}^\ad\sqrt{B}^\mpi\).
 Hence,
\begin{multline}\label{N368H}
 \pba{n}{z}\hp{2n}^\mpi\cpP +\pba{n+1}{z}\cpQ
 =\ek*{\pba{n}{z}\hp{2n}^\mpi\XFa{2n+1}{z}-\pba{n+1}{z}}\sqrt{B}^\mpi CP\\
-\rk*{\pba{n}{z}\hp{2n}^\mpi\ek{\XFa{2n+1}{z}}^\ad-\pba{n+1}{z}}\sqrt{B}^\mpi+\pba{n+1}{z}Q.
\end{multline}
 We consider an arbitrary \(v\in\nul{\pba{n}{z}\hp{2n}^\mpi \cpP+\pba{n+1}{z}\cpQ}\).
 From \eqref{N368H} then
\begin{multline}\label{L0827.1}
 \ek*{\pba{n}{z}\hp{2n}^\mpi\XFa{2n+1}{z}-\pba{n+1}{z}}\sqrt{B}^\mpi CPv+\pba{n+1}{z}Q v\\
 =\rk*{\pba{n}{z}\hp{2n}^\mpi\ek{\XFa{2n+1}{z}}^\ad-\pba{n+1}{z}}\sqrt{B}^\mpi v
\end{multline}
 follows.
 Since \rrem{R52SK} yields \(\det\pba{n}{z}\neq0\) and \(\XFa{2n+1}{z}=\hp{2n}\ek{\pba{n}{z}}^\inv\pba{n+1}{z}\), multiplication of equation \eqref{L0827.1} from the left by \(\hp{2n}\ek{\pba{n}{z}}^\inv\) provides
\begin{multline}\label{N38}
 \ek*{\hp{2n}\hp{2n}^\mpi\XFa{2n+1}{z}-\XFa{2n+1}{z}}\sqrt{B}^\mpi CPv+\XFa{2n+1}{z}Q v\\
 =\rk*{\hp{2n}\hp{2n}^\mpi\ek*{\XFa{2n+1}{z}}^\ad-\XFa{2n+1}{z}}\sqrt{B}^\mpi v.
\end{multline}
 \rprop{L1347} provides \(\ran{\ek{\XFa{2n+1}{z}}^\ad}=\ran{\XFa{2n+1}{z}}=\ran{\hp{2n}}\) and \(\nul{\rim\XFa{2n+1}{z}}=\nul{\XFa{2n+1}{z}}=\nul{\hp{2n}}\).
 In view of \rremp{tsa12}{tsa12.a}, we obtain then \(\hp{2n}\hp{2n}^\mpi\ek{\XFa{2n+1}{z}}^\ad=\ek{\XFa{2n+1}{z}}^\ad\) and \(\hp{2n}\hp{2n}^\mpi\XFa{2n+1}{z}=\XFa{2n+1}{z}\).
 Consequently, from \eqref{N38} we get
\beql{M38}
 \XFa{2n+1}{z}Q v
 =\rk*{\ek*{\XFa{2n+1}{z}}^\ad-\XFa{2n+1}{z}}\sqrt{B}^\mpi v.
\eeq
 \rrem{R.P} yields \(\ran{Q}=\nul{\hp{2n}}\).
 In view of \(\nul{\XFa{2n+1}{z}}=\nul{\hp{2n}}\), hence \(\ran{Q}=\nul{\XFa{2n+1}{z}}\), implying \(\XFa{2n+1}{z}Qv=\Ouu{q}{1}\).
 Taking additionally into account \(\rim\rk{z}B=\rim\XFa{2n+1}{z}\), from \eqref{M38} we can conclude \(\Ouu{q}{1}=-2\iu\rim\rk{z}B\sqrt{B}^\mpi v\).
 Since \rrem{A.R.A+sqrt} yields \(B\sqrt{B}^\mpi=\sqrt{B}\), hence \(\Ouu{q}{1}=-2\iu\rim\rk{z}\sqrt{B}v\), which, in view of \(z\in\CR\), implies \(v\in\nul{\sqrt{B}}\).
 According to \rrem{A.R.r-sqrt}, we have \(\nul{\sqrt{B}}=\nul{B}\).
 In view of \(B=\rk{\rim z}^\inv\rim\XFa{2n+1}{z}\) and \(\nul{\rim\XFa{2n+1}{z}}=\nul{\hp{2n}}\), we thus obtain \(v\in\nul{\hp{2n}}\), implying \(Qv=v\).
 From \rrem{R.P} we see \(\nul{P}=\ran{\hp{2n}}^\oc\).
 According to \rrem{tsa2}, we have \(\ran{\hp{2n}}^\oc=\nul{\hp{2n}^\ad}\).
 Since \rrem{R1019} yields \(\hp{2n}^\ad=\hp{2n}\), then \(\nul{P}=\nul{\hp{2n}}\) follows, implying \(Pv=\Ouu{q}{1}\).
 By virtue of \rrem{A.R.r-mpi}, we see \(\nul{\sqrt{B}^\mpi}=\nul{\sqrt{B}^\ad}=\nul{\sqrt{B}}\).
 In view of \(v\in\nul{\sqrt{B}}\), thus \(\sqrt{B}^\mpi v=\Ouu{q}{1}\).
 Using additionally \(Pv=\Ouu{q}{1}\) and \(Qv=v\), from \eqref{L0827.1} we obtain then \(\pba{n+1}{z}v=\Ouu{q}{1}\).
 Since \rlem{BK8.6} yields \(\det\pba{n+1}{z}\neq0\), we get \(v=\Ouu{q}{1}\).
 Hence, \(\nul{\pba{n}{z}\hp{2n}^\mpi \cpP+\pba{n+1}{z}\cpQ}\subseteq\set{\Ouu{q}{1}}\), which completes the proof.
\eproof

 The following result plays an important role in the proofs of \rpropss{L1204}{L1301}, which provide the two sides of \rthm{T45CD}.

\bpropl{L1745}
 Let \(n\in\NO\), let \(\seqs{2n+1}\in\Hggeq{2n+1}\), and let \(z\in\CR\).
 Let \(P\defeq\OPu{\ran{\hp{2n}}}\) and \(Q\defeq\OPu{\nul{\hp{2n}}}\).
 Then there exist matrices \(\cpP,\cpQ\in\Cqq\) such that the following three conditions are fulfilled:
\baeqi{0}
 \il{L1745.i} \(\rim\rk{z}\rim\rk{\cpQ ^\ad \cpP }\in\Cggq\).
 \il{L1745.ii} \(P\cpP =\cpP \), \(\cpP P=\cpP \), and \(\cpQ P=\cpQ -Q\).
 \il{L1745.iii} \(\det\rk{\pba{n}{z}\hp{2n}^\mpi \cpP +\pba{n+1}{z}\cpQ }\neq0\).
\eaeqi
 If \(\cpP,\cpQ\in\Cqq\) are arbitrary matrices such that~\ref{L1745.i}--\ref{L1745.iii} are fulfilled, then the matrix \(\rk{\rim z}^\inv\rim\XFa{2n+1}{z}\) is \tnnH{}, the matrix
\[
 \cpC 
 \defeq\sqrt{\rk{\rim z}^\inv\rim\XFa{2n+1}{z}}^\mpi\rk*{\cpP +\ek*{\XFa{2n+1}{z}}^\ad \cpQ }\rk*{\cpP +\XFa{2n+1}{z}\cpQ }^\mpi\sqrt{\rk{\rim z}^\inv\rim\XFa{2n+1}{z}}
\]
 is contractive, and the identity
\begin{multline}\label{L1745.A}
 -\ek*{\paa{n}{z}\hp{2n}^\mpi \cpP +\paa{n+1}{z}\cpQ }\ek*{\pba{n}{z}\hp{2n}^\mpi \cpP +\pba{n+1}{z}\cpQ }^\inv\\
 =\mpka{2n+1}{z}+\rk{z-\ko z}^\inv \lrka{2n+1}{z}\cpC \rrka{2n+1}{z}
\end{multline}
 holds true, where \(\lrk{2n+1}\), \(\rrk{2n+1}\), and \(\mpk{2n+1}\) are given by \rnotap{CB15}{CB15.b}.
\eprop
\bproof
 \rrem{R1019} yields \(\hp{2n}^\ad=\hp{2n}\).
 According to \rrem{tsa2}, we have \(\nul{\hp{2n}^\ad}=\ran{\hp{2n}}^\oc\).
 Consequently, \(\nul{\hp{2n}}=\ran{\hp{2n}}^\oc\) follows.
 Using \rrem{A.R.0<P<1}, we can thus conclude \(P+Q=\Iq\).
 From \rlem{BK8.6} we easily see then that, \teg{}, \(\cpP=\Oqq\) and \(\cpQ=\Iq\) fulfill~\ref{L1745.i}--\ref{L1745.iii}.
 
 Now let \(\cpP,\cpQ\in\Cqq\) be arbitrarily chosen such that~\ref{L1745.i}--\ref{L1745.iii} are fulfilled.
 In order to apply \rlem{L1104}, we set \(E\defeq-\ek{\XFa{2n+1}{z}}^\ad\), \(b\defeq\rim\rk{z}\), and \(B\defeq b^\inv\rim E\).
 Clearly, then
\begin{align}\label{L1745.2}
 \XFa{2n+1}{z}&=-E^\ad,&
 \ek*{\XFa{2n+1}{z}}^\ad&=-E,&
&\text{and}&
 \rk{\rim z}^\inv\rim\XFa{2n+1}{z}&=B.
\end{align}
 Because of~\ref{L1745.ii} and \rrem{R.P}, we discern
\beql{L1745.1}
 \ran{\cpP}
 =\ran{P\cpP}
 \subseteq\ran{P}
 =\ran{\hp{2n}}.
\eeq
 Let \(\STA{2n+1}\) and \(\STB{2n+1}\) be given by \eqref{L0852.AB}.
 Regarding \eqref{L1745.1} and \eqref{L1745.2}, from \rlem{L0852} we can conclude then that the matrices \(\PDa{n}{z}\) and \(\pba{n}{z}\) are both invertible, that \(B\) is \tnnH{}, and that
\beql{L1745.11A}
 \rk{z-\ko z}\rk*{\ek*{\PDa{n}{z}}^\inv\PCa{n}{z}\STB{2n+1}-\STA{2n+1}}
 =\lrka{2n+1}{z}\sqrt{B}^\mpi\rk{\cpP-E\cpQ}
\eeq
 and
\beql{L1745.11B}
 \hp{2n}\ek*{\pba{n}{z}}^\inv\STB{2n+1}
 =\cpP-E^\ad\cpQ
\eeq
 hold true.
 In view of \eqref{L1745.2} and \rnotap{CB15}{CB15.b}, we have
\begin{align}\label{L1745.12}
 \cpC &=\sqrt{B}^\mpi\rk{\cpP-E\cpQ}\rk{\cpP-E^\ad\cpQ}^\mpi\sqrt{B}& 
&\text{and}&
 \rrka{2n+1}{z}&=\sqrt{B}^\mpi\hp{2n}\ek*{\pba{n}{z}}^\inv.
\end{align}
 \rprop{L1347} yields \(\ran{\rim\XFa{2n+1}{z}}=\ran{\ek{\XFa{2n+1}{z}}^\ad}=\ran{\XFa{2n+1}{z}}=\ran{\hp{2n}}\) and \(\nul{\rim\XFa{2n+1}{z}}=\nul{\ek{\XFa{2n+1}{z}}^\ad}=\nul{\XFa{2n+1}{z}}=\nul{\hp{2n}}\).
 Regarding \eqref{L1745.2}, consequently,
\begin{align}\label{L1745.4}
 \ran{B}&=\ran{E}=\ran{E^\ad}=\ran{\hp{2n}}&
&\text{and}&
 \nul{B}&=\nul{E}=\nul{E^\ad}=\nul{\hp{2n}}
\end{align}
 follow.
 Taking into account \eqref{L0852.AB} and~\ref{L1745.iii}, we perceive \(\det\STB{2n+1}\neq0\).
 Therefore, using \eqref{L1745.11B} and \rrem{R1241}, we obtain
\beql{L1745.5}
 \ran{\cpP-E^\ad\cpQ }
 =\Ran{\hp{2n}\ek*{\pba{n}{z}}^\inv\STB{2n+1}}
 =\ran{\hp{2n}}.
\eeq
 The combination of \eqref{L1745.1}, \eqref{L1745.4}, and \eqref{L1745.5} yields \(\ran{\cpP}\subseteq\ran{E}=\ran{B}=\ran{\cpP-E^\ad\cpQ}\).
 Regarding that \(z\in\CR\) implies \(b\in\R\setminus\set{0}\) and taking additionally into account the first equation in \eqref{L1745.12} and \(B\in\Cggq\), we hence can apply \rlem{L1104} and~\ref{L1745.i} to obtain that \(\cpC \) is contractive.
 From \rnotass{CB15}{N1326} we conclude
\beql{L1745.6}
 \mpka{2n+1}{z}
 =-\ek*{\PDa{n}{z}}^\inv\PCa{n}{z}.
\eeq
 From \rrem{R.P} and \eqref{L1745.4} we see \(\ran{Q}=\nul{\hp{2n}}=\nul{E}=\nul{E^\ad}\), implying \(EQ=\Oqq\) and \(E^\ad Q=\Oqq\).
 Hence, using additionally~\ref{L1745.ii}, we infer
\beql{L1745.16}
 \rk{\cpP-E\cpQ}P
 =\cpP P-E\cpQ P
 =\cpP-E\rk{\cpQ -Q}
 =\cpP-E\cpQ+EQ
 =\cpP-E\cpQ
\eeq
 and
\[
 \rk{\cpP-E^\ad\cpQ}P
 =\cpP P-E^\ad\cpQ P
 =\cpP-E^\ad\rk{\cpQ -Q}
 =\cpP-E^\ad\cpQ+E^\ad Q
 =\cpP-E^\ad\cpQ.
\]
 By virtue of \(P+Q=\Iq\), we can thus conclude \(\rk{\cpP-E^\ad\cpQ}Q=\Oqq\).
 In view of \(\ran{Q}=\nul{\hp{2n}}\), we can infer then \(\nul{\hp{2n}}\subseteq\nul{\cpP-E^\ad\cpQ}\).
 Furthermore, because of \eqref{L1745.5} and a familiar result of linear algebra, we acquire \(\dim\nul{\hp{2n}}=q-\dim\ran{\hp{2n}}=q-\dim\ran{\cpP-E^\ad\cpQ}=\dim\nul{\cpP-E^\ad\cpQ}<\infty\).
 Therefore, \(\nul{\hp{2n}}=\nul{\cpP-E^\ad\cpQ}\) follows.
 Using \rrem{tsa2}, we can infer then \(\ran{\hp{2n}^\ad}=\ran{\rk{\cpP-E^\ad\cpQ}^\ad}\).
 Taking additionally into account \(\hp{2n}^\ad=\hp{2n}\) and \rrem{tsb3}, we obtain \(P=\OPu{\ran{\hp{2n}}}=\OPu{\ran{\hp{2n}^\ad}}=\OPu{\ran{\rk{\cpP-E^\ad\cpQ}^\ad}}=\rk{\cpP-E^\ad\cpQ}^\mpi\rk{\cpP-E^\ad\cpQ}\).
 Consequently, using \eqref{L1745.16}, we recognize that
\beql{L1745.14}
  \rk{\cpP-E\cpQ}\rk{\cpP-E^\ad\cpQ}^\mpi\rk{\cpP-E^\ad\cpQ}
  =\rk{\cpP-E\cpQ}P
  =\cpP-E\cpQ
\eeq
 is fulfilled.
 Regarding \(\det\STB{2n+1}\neq0\), from \eqref{L1745.11B} we receive
\beql{L1745.13}
 \rk{\cpP-E^\ad\cpQ}\STB{2n+1}^\inv
 =\hp{2n}\ek*{\pba{n}{z}}^\inv.
\eeq
 \rrem{tsb3} yields \(\sqrt{B}\sqrt{B}^\mpi=\OPu{\ran{\sqrt{B}}}\).
 From \rrem{A.R.r-sqrt} and \eqref{L1745.4}, we discern \(\ran{\sqrt{B}}=\ran{B}=\ran{\hp{2n}}\).
 Consequently,
\beql{L1745.21}
 \sqrt{B}\sqrt{B}^\mpi\hp{2n}
 =\OPu{\ran{\sqrt{B}}}\hp{2n}
 =\OPu{\ran{\hp{2n}}}\hp{2n}
 =\hp{2n}.
\eeq
 Regarding that \(z\in\CR\) implies \(z-\ko z\neq0\) and using \(\det\STB{2n+1}\neq0\), \eqref{L1745.6}, \eqref{L1745.11A}, \eqref{L1745.14}, \eqref{L1745.13}, \eqref{L1745.21}, and \eqref{L1745.12}, finally we obtain
\[\begin{split}
  -\STA{2n+1}\STB{2n+1}^\inv-\mpka{2n+1}{z}
  &=-\STA{2n+1}\STB{2n+1}^\inv-\rk*{-\ek*{\PDa{n}{z}}^\inv\PCa{n}{z}}
  =\rk*{\ek*{\PDa{n}{z}}^\inv\PCa{n}{z}\STB{2n+1}-\STA{2n+1}}\STB{2n+1}^\inv\\
  &=\rk{z-\ko z}^\inv \lrka{2n+1}{z}\sqrt{B}^\mpi\rk{\cpP-E\cpQ}\STB{2n+1}^\inv\\
  &=\rk{z-\ko z}^\inv \lrka{2n+1}{z}\sqrt{B}^\mpi\rk{\cpP-E\cpQ}\rk{\cpP-E^\ad\cpQ}^\mpi\rk{\cpP-E^\ad\cpQ}\STB{2n+1}^\inv\\
  &=\rk{z-\ko z}^\inv \lrka{2n+1}{z}\sqrt{B}^\mpi\rk{\cpP-E\cpQ}\rk{\cpP-E^\ad\cpQ}^\mpi\hp{2n}\ek*{\pba{n}{z}}^\inv\\
  &=\rk{z-\ko z}^\inv \lrka{2n+1}{z}\sqrt{B}^\mpi\rk{\cpP-E\cpQ}\rk{\cpP-E^\ad\cpQ}^\mpi\sqrt{B}\sqrt{B}^\mpi\hp{2n}\ek*{\pba{n}{z}}^\inv\\
  &=\rk{z-\ko z}^\inv \lrka{2n+1}{z}\cpC \rrka{2n+1}{z}.
\end{split}\]
 Regarding~\ref{L1745.iii} and \eqref{L0852.AB}, thus \eqref{L1745.A} is proved as well.
\eproof

 We keep in mind that we target to achieve a parametrization of the set \(\setaca{F\rk{w}}{F\in\RFOqskg{2n}}\) with respect to an arbitrarily prescribed sequence \(\seqs{2n}\in\Hggeq{2n}\) and an arbitrarily chosen \(w\in\pip\).
 The next proposition provides the proof of one half of \rthm{T45CD}.

\bpropl{L1204}
 Let \(n\in\NO\), let \(\seqs{2n}\in\Hggeq{2n}\), and let \(F\in\RFOqskg{2n}\).
 For each \(w\in\pip\), then 
\beql{L1204.A}
 F\rk{w}
 \in\Cmb{\mpka{2n}{w}}{\rk{w-\ko w}^\inv\lrka{2n}{w}}{\rrka{2n}{w}},
\eeq
 where \(\lrk{2n}\), \(\rrk{2n}\), and \(\mpk{2n}\) are given by \rnotap{CB15}{CB15.a}.
\eprop
\bproof
 Let \(\sHp{2n}\) be the \tsHpo{\(\seqs{2n}\)} and let \(\abcd{n}\) be the \sabcdo{\(\seqs{2n}\)}.
 Using \rthmp{CM123}{CM123.b}, we perceive that there exists a pair \(\copa{\phi}{\psi}\in\PRFqa{\hp{2n}}\) such that both \(\phi\) and \(\psi\) are holomorphic in \(\pip\) and that \eqref{L1204.1} and \eqref{L1204.2} 
 hold true for all \(z\in\pip\).
 \rrem{R1019} yields \(\hp{2n}^\ad=\hp{2n}\).
 Hence, setting \(P\defeq\OPu{\ran{\hp{2n}}}\) and \(Q\defeq\OPu{\nul{\hp{2n}}}\), from \rlem{CfolA15} we obtain that there exists a pair \(\copa{\cpfP}{\cpfQ}\in\cpcl{\phi}{\psi}\) such that all four equations in \eqref{ST37N} are valid.
 In particular, in view of \rrem{Ma09_1.11}, we recognize that \(\copa{\cpfP}{\cpfQ}\in\PRFq\) and that there exist a \tqqa{matrix}-valued function \(g\) meromorphic in \(\pip\) and a discrete subset \(\mD_1\) of \(\pip\) such that \(\phi\), \(\psi\), \(\cpfP\), \(\cpfQ\), and \(g\) are holomorphic in \(\pip\setminus\mD_1\) and that
\begin{align}\label{L1204.3}%
 \cpfP\rk{z}&=\phi\rk{z}g\rk{z},&
 \cpfQ\rk{z}&=\psi\rk{z}g\rk{z},&
&\text{and}&
 \det g\rk{z}&\neq0
\end{align}
 as well as 
\begin{align}\label{L1204.4}%
 P\cpfP\rk{z}&=\cpfP\rk{z},&
 \cpfP\rk{z}P&=\cpfP\rk{z},&
&\text{and}&
 \cpfQ\rk{z}P
 =\cpfQ\rk{z}-Q
\end{align}
 hold true for all \(z\in\pip\setminus\mD_1\).
 Using the first two equations in \eqref{L1204.3}, we obtain \(\paa{n}{z}\hp{2n}^\mpi\cpfP\rk{z}+\paoa{n+1}{z}\cpfQ\rk{z}=\ek{\paa{n}{z}\hp{2n}^\mpi\phi\rk{z}+\paoa{n+1}{z}\psi\rk{z}}g\rk{z}\) and \(\pba{n}{z}\hp{2n}^\mpi\cpfP\rk{z}+\pboa{n+1}{z}\cpfQ\rk{z}=\ek{\pba{n}{z}\hp{2n}^\mpi\phi\rk{z}+\pboa{n+1}{z}\psi\rk{z}}g\rk{z}\) for all \(z\in\pip\setminus\mD_1\).
 Taking additionally into account the last equation in \eqref{L1204.3} as well as \eqref{L1204.1} and \eqref{L1204.2}, we can conclude
\beql{L1204.5}%
 \det\rk*{\pba{n}{z}\hp{2n}^\mpi \cpfP\rk{z}+\pboa{n+1}{z}\cpfQ\rk{z}}
 =\det\rk*{\pba{n}{z}\hp{2n}^\mpi\phi\rk{z}+\pboa{n+1}{z}\psi\rk{z}}\cdot\det g\rk{z}
 \neq0
\eeq
 and  
\beql{L1204.6}\begin{split}%
 &-\ek*{\paa{n}{z}\hp{2n}^\mpi \cpfP\rk{z}+\paoa{n+1}{z}\cpfQ\rk{z}}\ek*{\pba{n}{z}\hp{2n}^\mpi \cpfP\rk{z}+\pboa{n+1}{z}\cpfQ\rk{z}}^\inv\\ 
 &=-\rk*{\ek*{\paa{n}{z}\hp{2n}^\mpi\phi\rk{z}+\paoa{n+1}{z}\psi\rk{z}}g\rk{z}}\rk*{\ek*{\pba{n}{z}\hp{2n}^\mpi\phi\rk{z}+\pboa{n+1}{z}\psi\rk{z}}g\rk{z}}^\inv\\
 &=-\ek*{\paa{n}{z}\hp{2n}^\mpi\phi\rk{z}+\paoa{n+1}{z}\psi\rk{z}}\ek*{\pba{n}{z}\hp{2n}^\mpi\phi\rk{z}+\pboa{n+1}{z}\psi\rk{z}}^\inv
 =F\rk{z}
\end{split}\eeq
 for all \(z\in\pip\setminus\mD_1\).
 Regarding \(\copa{\cpfP}{\cpfQ}\in\PRFq\) and \rdefn{def-nev-paar}, we realize that \(\cpfP\) and \(\cpfQ\) are \tqqa{matrix}-valued functions meromorphic in \(\pip\), that there is a discrete subset \(\mD_2\) of \(\pip\) such that \(\cpfP\) and \(\cpfQ\) are holomorphic in \(\pip\setminus\mD_2\), and that, in view of \rrem{IG.b}, moreover, 
\beql{L1204.7}%
 \rim\rk*{\ek*{\cpfQ\rk{z}}^\ad \cpfP\rk{z}}
 \in\Cggq 
\eeq
 is valid for each \(z\in\pip\setminus\mD_2\).
 
 Now we consider an arbitrary \(w\in\pip\).
 Obviously, \(\mD\defeq\mD_1\cup\mD_2\) is a discrete subset of \(\pip\).
 Consequently, it is feasible to choose a sequence \(\rk{z_\ell}_{\ell=1}^\infi\) of points belonging to \(\pip\setminus\mD\) converging to \(w\).
 Taking into account that \(\cpfP\) and \(\cpfQ\) are holomorphic in \(\pip\setminus\mD\) and that \eqref{L1204.4} and \eqref{L1204.7} hold true for all \(z\in\pip\setminus\mD\), for each \(\ell\in\N\), the matrices \(\cpP_\ell\defeq \cpfP\rk{z_\ell}\) and \(\cpQ_\ell\defeq \cpfQ\rk{z_\ell}\) fulfill
\begin{align}\label{L1204.8}%
 P\cpP_\ell&=\cpP_\ell,&
 \cpP_\ell P&=\cpP_\ell,&
 \cpQ_\ell P&=\cpQ_\ell-Q,&
&\text{and}&
 \rim\rk{z_\ell}\rim\rk{\cpQ_\ell^\ad \cpP_\ell}&\in\Cggq.
\end{align}
 Since \eqref{L1204.5} and \eqref{L1204.6} hold true for all \(z\in\pip\setminus\mD\), we receive 
\beql{L1204.9}%
 \det\rk*{\pba{n}{z_\ell}\hp{2n}^\mpi \cpP_\ell+\pboa{n+1}{z_\ell}\cpQ_\ell}
 \neq0
\eeq
and 
\beql{L1204.10}%
 F\rk{z_\ell}
 =-\ek*{\paa{n}{z_\ell}\hp{2n}^\mpi \cpP_\ell+\paoa{n+1}{z_\ell}\cpQ_\ell}\ek*{\pba{n}{z_\ell}\hp{2n}^\mpi \cpP_\ell+\pboa{n+1}{z_\ell}\cpQ_\ell}^\inv 
\eeq
 for all \(\ell\in\N\).
 Let \(\seqsh{2n+1}\) be the \tnatexto{\(\seqs{2n}\)} and let \(\habcd{n+1}\) be the \sabcdo{\(\seqsh{2n+1}\)}.
 \rlem{L1237-B} yields then \(\pah{n}=\pa{n}\) and \(\pbh{n}=\pb{n}\) as well as \(\pah{n+1}=\pao{n+1}\) and \(\pbh{n+1}=\pbo{n+1}\).
 Let \(\shHp{2n+1}\) be the \tsHpo{\(\seqsh{2n+1}\)}.
 According to \rrem{R26-SK1}, we have \(\hhp{2n}=\hp{2n}\).
 Consequently, regarding \eqref{L1204.9} and \eqref{L1204.10} we can conclude
\beql{L1204.11}%
 \det\rk*{\pbh{n}\rk{z_\ell}\hhp{2n}^\mpi\cpP_\ell+\pbh{n+1}\rk{z_\ell}\cpQ_\ell}
 \neq0
\eeq
 and 
\beql{L1204.12}%
 F\rk{z_\ell}
 =-\ek*{\pah{n}\rk{z_\ell}\hhp{2n}^\mpi\cpP_\ell+\pah{n+1}\rk{z_\ell}\cpQ_\ell}\ek*{\pbh{n}\rk{z_\ell}\hhp{2n}^\mpi\cpP_\ell+\pbh{n+1}\rk{z_\ell}\cpQ_\ell}^\inv 
\eeq
 for all \(\ell\in\N\).
 From \rprop{L1237-A} we infer that \(\seqsh{2n+1}\) belongs to \(\Hggeq{2n+1}\).
 Regarding \(\hhp{2n}=\hp{2n}\), we see that \(P=\OPu{\ran{\hhp{2n}}}\) and \(Q=\OPu{\nul{\hhp{2n}}}\).
 Thus, in view of \eqref{L1204.8}, \eqref{L1204.11}, and \eqref{L1204.12}, the application of \rprop{L1745} provides that, for all \(\ell\in\N\), there exists a contractive complex \tqqa{matrix} \(\cpC _\ell\) such that
\beql{L1204.13}%
 F\rk{z_\ell}
 =\mphk{2n+1}\rk{z_\ell}+\rk{z_\ell-\ko{z_\ell}}^\inv\lrhk{2n+1}\rk{z_\ell}\cpC _\ell\rrhk{2n+1}\rk{z_\ell}
\eeq
 holds true, where \(\lrhk{2n+1}\), \(\rrhk{2n+1}\), and \(\mphk{2n+1}\) are built according to \rnotap{CB15}{CB15.b} from the sequence \(\seqsh{2n+1}\).
 Now we continue by a natural limit procedure.
 Since \((\cpC _\ell)_{\ell=1}^\infi\) is a sequence of contractive complex \tqqa{matrices}, \rrem{KTK} yields the existence of a subsequence \((\cpC _{\ell_m})_{m=1}^\infi\) and a contractive complex \tqqa{matrix} \(\cpC\) such that \(\lim_{m\to\infty}\cpC_{\ell_m}=\cpC\).
 Since \(F\) belongs to \(\RFOqskg{2n}\), we recognize that \(F\) is holomorphic in \(\pip\) and, in particular, continuous.
 Thus, \(\lim_{m\to\infty} z_{\ell_m}=w\) implies \(\lim_{m\to\infty} F\rk{z_{\ell_m}}=F\rk{w}\).
 Regarding \(\seqsh{2n+1}\in\Hggeq{2n+1}\), from \rlem{CL18} we discern that all the matrix-valued functions \(\lrhk{2n+1}\), \(\rrhk{2n+1}\), and \(\mphk{2n+1}\) are continuous in \(\CR\).
 Thus, \(\lim_{m\to\infty}\lrhka{2n+1}{z_{\ell_m}}=\lrhka{2n+1}{w}\) and \(\lim_{m\to\infty}\rrhka{2n+1}{z_{\ell_m}}=\rrhka{2n+1}{w}\) as well as \(\lim_{m\to\infty}\mphka{2n+1}{z_{\ell_m}}=\mphka{2n+1}{w}\).
 Consequently, letting \(m\to\infty\), from \eqref{L1204.13} we perceive \(F\rk{w}=\mphka{2n+1}{w}+\rk{w-\ko w}^\inv\lrhka{2n+1}{w}\cpC \rrhka{2n+1}{w}\).
 In view of \(\seqsh{2n+1}\in\Hggeq{2n+1}\), from \rlem{C220N} we obtain \(\lrhk{2n}=\lrhk{2n+1}\) and \(\rrhk{2n}=\rrhk{2n+1}\), whereas \rlem{L3N10} yields \(\mphk{2n}=\mphk{2n+1}\), where \(\lrhk{2n}\), \(\rrhk{2n}\), and \(\mphk{2n}\) are built according to \rnotap{CB15}{CB15.a} from the sequence \(\seqsh{2n+1}\).
 By virtue of \rdefn{D.nat-ext} and \rrem{ABM-tru}, we have \(\lrhk{2n}=\lrk{2n}\) and \(\rrhk{2n}=\rrk{2n}\) as well as \(\mphk{2n}=\mpk{2n}\).
 Therefore, we can finally conclude \(F\rk{w}=\mpka{2n}{w}+\rk{w-\ko w}^\inv\lrka{2n}{w}\cpC \rrka{2n}{w}\).
 Since \(\cpC \) is contractive, then \eqref{L1204.A} is proved. 
\eproof

 Now we are able to state and prove the missing other side of the proof of \rthm{T45CD}.
 As in the proof of \rprop{L1204} we will see that the \tnatext{} \(\seqsh{2n+1}\) of \(\seqs{2n}\) plays an essential role.

\bpropl{L1301}
 Let \(n\in\NO\) and let \(\seqs{2n}\in\Hggeq{2n}\).
 Furthermore, let \(w\in\pip\) and let \(C\) be a contractive complex \tqqa{matrix}.
 Then there exists a rational matrix-valued function \(F\) which belongs to \(\RFOqskg{2n}\) and which fulfills \(F\rk{w}=\mpka{2n}{w}+{\rk{w-\ko w}}^\inv\lrka{2n}{w}C\rrka{2n}{w}\), where \(\lrk{2n}\), \(\rrk{2n}\), and \(\mpk{2n}\) are given by \rnotap{CB15}{CB15.a}.
\eprop
\bproof
 \rprop{L1237-A} shows that the \tnatext{} \(\seqsh{2n+1}\) of \(\seqs{2n}\) belongs to \(\Hggeq{2n+1}\).
 Let \(\seqXh{2n+1}\) be the \tsXFo{\(\seqsh{2n+1}\)}.
 Let \(\Ec\defeq-\ek{\XFha{2n+1}{w}}^\ad\), let \(b\defeq\rim\rk{w}\), let \(\Bc\defeq b^\inv\rim\Ec\), and let \(\Pc\defeq\OPu{\ran{\Ec}}\) and \(\Qc\defeq\OPu{\nul{\Ec}}\).
 Clearly, then
\begin{align}\label{L1301.1}
 \XFha{2n+1}{w}&=-\Ec^\ad,&
 \ek*{\XFha{2n+1}{w}}^\ad&=-\Ec,&
&\text{and}&
 \rk{\rim w}^\inv\rim\XFha{2n+1}{w}&=\Bc.
\end{align}
 From \rcor{C1232} we can infer then \(\Bc\in\Cggq\).
 Obviously, the (constant) matrix-valued functions \(\phi,\psi\colon\pip\to\Cqq\) defined by
\begin{align}\label{L1301.2}%
 \phi\rk{z}&\defeq\Ec\sqrt{\Bc}^\mpi-\Ec^\ad\sqrt{\Bc}^\mpi C\Pc&
&\text{and}&
 \psi\rk{z}&\defeq\sqrt{\Bc}^\mpi-\sqrt{\Bc}^\mpi C\Pc+\Qc,
\end{align}
 respectively, fulfill
\beql{L1301.11} 
 \holpt{\phi}
 =\holpt{\psi}
 =\pip.
\eeq
 Let \(\shHp{2n+1}\) be the \tsHpo{\(\seqsh{2n+1}\)}.
 \rprop{L1347} yields \(\ran{\rim\XFha{2n+1}{z}}=\ran{\ek{\XFha{2n+1}{z}}^\ad}=\ran{\hhp{2n}}\) and \(\nul{\rim\XFha{2n+1}{z}}=\nul{\ek{\XFha{2n+1}{z}}^\ad}=\nul{\hhp{2n}}\) for all \(z\in\pip\).
 By virtue of \eqref{L1301.1}, consequently,
\begin{align}\label{L1301.3}
 \ran{\Bc}&=\ran{\Ec}=\ran{\hhp{2n}}&
&\text{and}&
 \nul{\Bc}&=\nul{\Ec}=\nul{\hhp{2n}}
\end{align}
 follow.
 Regarding that \(w\in\pip\) implies \(b\in(0,\infp)\) and taking into account \(B\in\Cggq\) and \(\ran{\Bc}=\ran{\Ec}\) as well as \eqref{L1301.2} and \(C\in\Kqq\), we can apply \rlem{L2043} to discern that
\begin{align}\label{L1301.4}
 \rank\Mat{\phi\rk{z}\\\psi\rk{z}}&=q&
&\text{and}&
 \rim\rk*{\ek*{\psi\rk{z}}^\ad \phi\rk{z}}&\in\Cggq,
\end{align}
 that
\begin{align}\label{L1301.5}
 P\phi\rk{z}&=\phi\rk{z},&
 \phi\rk{z}P&=\phi\rk{z},&
&\text{and}&
 \psi\rk{z}P&=\psi\rk{z}-Q,
\end{align}
 and that
\beql{L1301.6}
 \sqrt{B}^\mpi\ek*{\phi\rk{z}-E\psi\rk{z}}\ek*{\phi\rk{z}-E^\ad \psi\rk{z}}^\mpi\sqrt{B}
 =PCP
\eeq
 hold true for all \(z\in\pip\).
 In view of \rdefn{def-nev-paar} and \rrem{IG.b}, from \eqref{L1301.11} and \eqref{L1301.4} we recognize that the pair \(\copa{\phi}{\psi}\) belongs to \(\PRFq\).
 Let \(\sHp{2n}\) be the \tsHpo{\(\seqs{2n}\)}.
 \rrem{R26-SK1} yields \(\hhp{2n}=\hp{2n}\).
 Regarding \eqref{L1301.3}, consequently
\begin{align}\label{L1301.0}
 P&=\OPu{\ran{E}}=\OPu{\ran{\hhp{2n}}}=\OPu{\ran{\hp{2n}}}&
&\text{and}&
 Q&=\OPu{\nul{E}}=\OPu{\nul{\hhp{2n}}}=\OPu{\nul{\hp{2n}}}.
\end{align}
 From \eqref{L1301.0} and \eqref{L1301.5} we receive \(\OPu{\ran{\hp{2n}}}\phi\rk{z}=\phi\rk{z}\) for all \(z\in\pip\).
 According to \rnota{CbezA5}, hence \(\copa{\phi}{\psi}\in\PRFqa{\hp{2n}}\).
 Let \(\abcd{n}\) be the \sabcdo{\(\seqs{2n}\)}.
 In view \rnota{N-abcdO}, let \(\pars{n}\), \(\pbrs{n}\), \(\paors{n+1}\), and \(\pbors{n+1}\) be the restrictions of \(\pa{n}\), \(\pb{n}\), \(\pao{n+1}\), and \(\pbo{n+1}\) onto \(\pip\), respectively.
 From \rthmp{CM123}{CM123.a} we can then infer that the function \(\det\rk{\pbrs{n}\hp{2n}^\mpi\phi+\pbors{n+1}\psi}\) does not vanish identically in \(\pip\) and that the matrix-valued function
\beql{L1301.7}%
 F
 \defeq-\rk{\pars{n}\hp{2n}^\mpi\phi+\paors{n+1}\psi}\rk{\pbrs{n}\hp{2n}^\mpi\phi+\pbors{n+1}\psi}^\inv 
\eeq
 belongs to \(\RFOqskg{2n}\).
 Since from \eqref{L1301.2} we recognize that \(\phi\) and \(\psi\) are constant matrix-valued functions and since \(\pars{n}\), \(\pbrs{n}\), \(\paors{n+1}\), and \(\pbors{n+1}\) are restrictions of matrix polynomials, we see that \(F\) is a rational matrix-valued function.
 Let
\begin{align}\label{L1301.10}
 \cpP&\defeq\phi\rk{w}&
&\text{and}&
 \cpQ&\defeq\psi\rk{w}.
\end{align}
 By virtue of \eqref{L1301.10} and \eqref{L1301.2}, we see then that
\begin{align}\label{L1301.9}
 \cpP&=\Ec\sqrt{\Bc}^\mpi-\Ec^\ad\sqrt{\Bc}^\mpi C\Pc&
&\text{and}&
 \cpQ&=\sqrt{\Bc}^\mpi-\sqrt{\Bc}^\mpi C\Pc+\Qc
\end{align}
 hold true.
 Now we are going to justify that all assumptions for the application of \rprop{L1745} to the sequence \(\seqsh{2n+1}\) are satisfied.
 Let \(\habcd{n+1}\) be the \sabcdo{\(\seqsh{2n+1}\)}.
 According to \eqref{L1301.0}, we have \(\Pc=\OPu{\ran{\hhp{2n}}}\) and \(\Qc=\OPu{\nul{\hhp{2n}}}\).
 Using additionally \(\seqsh{2n+1}\in\Hggeq{2n+1}\) and \(\Ec=-\ek{\XFha{2n+1}{w}}^\ad\) as well as \(\Bc=\rk{\rim w}^\inv\rim\Ec\) and \eqref{L1301.9}, we can apply \rlem{L0827} to obtain
\beql{L1301.8}
 \det\rk*{\pbha{n}{w}\hhp{2n}^\mpi\cpP+\pbha{n+1}{w}\cpQ}
 \neq0.
\eeq
 \rlem{L1237-B} yields \(\pah{n}=\pa{n}\) and \(\pbh{n}=\pb{n}\) as well as \(\pah{n+1}=\pao{n+1}\) and \(\pbh{n+1}=\pbo{n+1}\).
 Taking additionally into account the definition of \(\pars{n}\), \(\pbrs{n}\), \(\paors{n+1}\), and \(\pbors{n+1}\) as well as \(\hhp{2n}=\hp{2n}\), from \eqref{L1301.11}, \eqref{L1301.10}, and \eqref{L1301.8} we can conclude \(\det\rk{\pbrsa{n}{w}\hp{2n}^\mpi\phi\rk{w}+\pborsa{n+1}{w}\psi\rk{w}}\neq0\) and, in view of \eqref{L1301.7}, consequently
\begin{multline}\label{L1301.15}%
 -\ek*{\paha{n}{w}\hhp{2n}^\mpi\cpP+\paha{n+1}{w}\cpQ}\ek*{\pbha{n}{w}\hhp{2n}^\mpi\cpP+\pbha{n+1}{w}\cpQ}^\inv\\
 =-\ek*{\parsa{n}{w}\hp{2n}^\mpi\phi\rk{w}+\paorsa{n+1}{w}\psi\rk{w}}\ek*{\pbrsa{n}{w}\hp{2n}^\mpi\phi\rk{w}+\pborsa{n+1}{w}\psi\rk{w}}^\inv
 =F\rk{w}.
\end{multline}
 In view of \(\rim\rk{w}\in(0,\infp)\) and \eqref{L1301.10}, from \eqref{L1301.4} and \eqref{L1301.5}, we can conclude \(\rim\rk{w}\rim\rk{\cpQ^\ad\cpP}\in\Cggq\) as well as \(P\cpP=\cpP\), \(\cpP P=\cpP\), and \(\cpQ P=\cpQ-Q\).
 Regarding \(\seqsh{2n+1}\in\Hggeq{2n+1}\) as well as \(\Pc=\OPu{\ran{\hhp{2n}}}\) and \(\Qc=\OPu{\nul{\hhp{2n}}}\) and using additionally \eqref{L1301.8}, we can thus apply \rprop{L1745} to infer that the matrix
\beql{L1301.12}
 \mathbf{K} 
 \defeq\sqrt{\rk{\rim w}^\inv\rim\XFha{2n+1}{w}}^\mpi\rk*{\cpP+\ek*{\XFha{2n+1}{w}}^\ad\cpQ}\rk*{\cpP+\XFha{2n+1}{w}\cpQ}^\mpi\sqrt{\rk{\rim w}^\inv\rim\XFha{2n+1}{w}}
\eeq
 fulfills
\begin{multline}\label{L1301.13}%
 -\ek*{\paha{n}{w}\hhp{2n}^\mpi\cpP +\paha{n+1}{w}\cpQ }\ek*{\pbha{n}{w}\hhp{2n}^\mpi \cpP +\pbha{n+1}{w}\cpQ }^\inv\\
 =\mphka{2n+1}{w}+\rk{w-\ko w}^\inv \lrhka{2n+1}{w}\mathbf{K}  \rrhka{2n+1}{w},
 \end{multline}
 where \(\lrhk{2n+1}\), \(\rrhk{2n+1}\), and \(\mphk{2n+1}\) are built according to \rnotap{CB15}{CB15.b} from the sequence \(\seqsh{2n+1}\).
 By virtue of \eqref{L1301.12}, \eqref{L1301.1}, \eqref{L1301.10}, and \eqref{L1301.6}, we discern
\beql{L1301.17}
 \mathbf{K} 
 =\sqrt{B}^\mpi\ek*{\phi\rk{w}-E\psi\rk{w}}\ek*{\phi\rk{w}-E^\ad \psi\rk{w}}^\mpi\sqrt{B}
 =\Pc C\Pc.
\eeq
 In view of \(\seqsh{2n+1}\in\Hggeq{2n+1}\), from \rlem{C220N} we obtain \(\lrhka{2n}{w}=\lrhka{2n+1}{w}\) and \(\rrhka{2n}{w}=\rrhka{2n+1}{w}\), whereas \rlem{L3N10} yields \(\mphka{2n}{w}=\mphka{2n+1}{w}\), where \(\lrhk{2n}\), \(\rrhk{2n}\), and \(\mphk{2n}\) are built according to \rnotap{CB15}{CB15.a} from the sequence \(\seqsh{2n+1}\).
 Because of \rdefn{D.nat-ext} and \rrem{ABM-tru}, we can discern \(\lrhka{2n}{w}=\lrka{2n}{w}\) and \(\rrhka{2n}{w}=\rrka{2n}{w}\) as well as \(\mphka{2n}{w}=\mpka{2n}{w}\).
 Consequently,
\begin{align}\label{L1301.14}%
 \lrhka{2n+1}{w}&=\lrka{2n}{w},&
 \rrhka{2n+1}{w}&=\rrka{2n}{w},&
&\text{and}&
 \mphka{2n+1}{w}&=\mpka{2n}{w}.
\end{align}
 From \rprop{CL19} we get \(\nul{\lrka{2n}{w}}=\nul{\hp{2n}}\) and \(\ran{\rrka{2n}{w}}=\ran{\hp{2n}}\).
 According to \eqref{L1301.0}, we have \(\Pc=\OPu{\ran{\hp{2n}}}\) and \(\Qc=\OPu{\nul{\hp{2n}}}\).
 Clearly, then \(P\rrka{2n}{w}=\rrka{2n}{w}\) follows.
 By virtue of \rrem{R.P}, we see furthermore \(\ran{P}=\ran{\hp{2n}}\) and \(\ran{Q}=\nul{\hp{2n}}\).
 Hence, we get \(\lrka{2n}{w}Q=\Oqq\).
 \rrem{tsa2} provides \(\nul{\hp{2n}^\ad}=\ran{\hp{2n}}^\oc\).
 Since \rrem{R1019} yields \(\hp{2n}^\ad=\hp{2n}\), then \(\nul{\hp{2n}}=\ran{\hp{2n}}^\oc\) follows.
 Thus, \rrem{A.R.0<P<1} yields \(P+Q=\Iq\), implying \(\lrka{2n}{w}P=\lrka{2n}{w}\).
 Consequently, we infer \(\lrka{2n}{w}\Pc C\Pc\rrka{2n}{w}=\lrka{2n}{w}C\rrka{2n}{w}\).
 Using additionally \eqref{L1301.15}, \eqref{L1301.13}, \eqref{L1301.17}, and \eqref{L1301.14}, we finally conclude
\[\begin{split}
 F\rk{w}
 &=-\ek*{\paha{n}{w}\hhp{2n}^\mpi\cpP+\paha{n+1}{w}\cpQ}\ek*{\pbha{n}{w}\hhp{2n}^\mpi\cpP+\pbha{n+1}{w}\cpQ}^\inv\\
 &=\mphka{2n+1}{w}+\rk{w-\ko w}^\inv \lrhka{2n+1}{w}\mathbf{K}  \rrhka{2n+1}{w}\\
 &=\mpka{2n}{w}+\rk{w-\ko w}^\inv \lrka{2n}{w}\Pc C\Pc\rrka{2n}{w}
 =\mpka{2n}{w}+\rk{w-\ko w}^\inv \lrka{2n}{w}C\rrka{2n}{w}.\qedhere
\end{split}\]
\eproof

 Note that \rprop{L1301} provides the remarkable additional information that each matrix from the matrix ball \(\cmb{\mpka{2n}{w}}{\rk{w-\ko w}^\inv\lrka{2n}{w}}{\rrka{2n}{w}}\) can be attained as value of a rational \tqqa{matrix}-valued function belonging to \(\RFOqskg{2n}\).
 In conclusion, we embrace the opportunity to prove the principal \rthm{T45CD}: 

\bproof[Proof of \rthm{T45CD}]
 Combine \rpropss{L1204}{L1301}.
\eproof

 We add the following result, which describes, for arbitrarily given \(w\in\pip\), explicitly a rational matrix-valued function \(F\) belonging to \(\RFOqsg{2n}\), the value \(F(w)\) of which coincides with the center of the matrix ball considered in \rthm{T45CD}.

\bpropl{P1004}
 Let \(n\in\NO\) and let \(\seqs{2n}\in\Hggeq{2n}\) with \tsHp{} \(\sHp{2n}\), \sabcd{} \(\abcd{n}\), and \tsXF{} \(\seqX{2n}\).
 Denote by \(\pars{n}\), \(\pbrs{n}\), \(\paors{n+1}\), and \(\pbors{n+1}\) the restrictions of \(\pa{n}\), \(\pb{n}\), \(\pao{n+1}\), and \(\pbo{n+1}\) onto \(\pip\), respectively, where \(\pao{n+1}\) and \(\pbo{n+1}\) are defined in \rnota{N-abcdO}.
 Let \(w\in\pip\) and let \(G\colon\pip\to\Cqq\) be defined by \(G\rk{z}\defeq-\ek{\XFa{2n}{w}}^\ad\).
 Then \(G\in\Pqevena{\hp{2n}}\), the function \(\det\rk{\pbrs{n}\hp{2n}^\mpi G+\pbors{n+1}}\) does not vanish identically, and \(F\defeq-\rk{\pars{n}\hp{2n}^\mpi G+\paors{n+1}}\rk{\pbrs{n}\hp{2n}^\mpi G+\pbors{n+1}}^\inv\) is a rational matrix-valued function belonging to \(\RFOqsg{2n}\) fulfilling \(F\rk{w}=\mpka{2n}{w}\), where \(\mpk{2n}\) is defined in \rnotap{CB15}{CB15.a}.
\eprop
\bproof
 From \rcor{C1232} we can infer \(\rim\XFa{2n}{w}\in\Cggq\).
 Setting \(E\defeq-\ek{\XFa{2n}{w}}^\ad\), hence \(\rim E=-\rim\rk{\ek{\XFa{2n}{w}}^\ad}=\rim\XFa{2n}{w}\in\Cggq\) follows.
 \rprop{L1347} yields \(\nul{\ek{\XFa{2n}{w}}^\ad}=\nul{\hp{2n}}\).
 Consequently, \(\nul{E}=\nul{\hp{2n}}\).
 Using \rlem{L1654}, we can then infer \(G\in\Pqevena{\hp{2n}}\).
 Thus, \rthmp{T1153}{T1153.a} shows that \(\det\rk{\pbrs{n}\hp{2n}^\mpi G+\pbors{n+1}}\) does not vanish identically, and that \(F\) belongs to \(\RFOqsg{2n}\).
 From \rprop{L1237-A} we see that the \tnatext{} \(\seqsh{2n+1}\) of \(\seqs{2n}\) belongs to \(\Hggeq{2n+1}\).
 Hence, we can apply \rlem{L0819} to the sequence \(\seqsh{2n+1}\) to obtain \(\det\PBha{n}{w}\neq0\) and \(\mphka{2n+1}{w}=-\PAha{n}{w}\ek{\PBha{n}{w}}^\inv\), where \(\mphk{2n+1}\) and \(\PAh{n},\PBh{n}\) are built according to \rnotass{CB15}{N1326}, \tresp{}, from the sequence \(\seqsh{2n+1}\).
 Let \(\habcd{n+1}\) be the \sabcdo{\(\seqsh{2n+1}\)}.
 \rlem{L1237-B} yields then \(\pah{n}=\pa{n}\) and \(\pbh{n}=\pb{n}\) as well as \(\pah{n+1}=\pao{n+1}\) and \(\pbh{n+1}=\pbo{n+1}\).
 Let \(\shHp{2n+1}\) be the \tsHpo{\(\seqsh{2n+1}\)}.
 \rrem{R26-SK1} yields then \(\hhp{2n}=\hp{2n}\).
 Let \(\seqXh{2n+1}\) be the \tsXFo{\(\seqsh{2n+1}\)}.
 \rlem{L1520} shows that \(\XFha{2n+1}{w}=\XFa{2n}{w}\) holds true.
 Regarding additionally \rnota{N1326}, we thus get
\begin{multline*}
 \parsa{n}{w}\hp{2n}^\mpi G\rk{w}+\paorsa{n+1}{w}
 =-\paa{n}{w}\hp{2n}^\mpi\ek*{\XFa{2n}{w}}^\ad+\paoa{n+1}{w}\\
 =-\paha{n}{w}\hhp{2n}^\mpi\ek*{\XFha{2n+1}{w}}^\ad+\paha{n+1}{w}
 =-\rk*{\paha{n}{w}\hhp{2n}^\mpi\ek*{\XFha{2n+1}{w}}^\ad-\paha{n+1}{w}}
 =-\PAha{n}{w}
\end{multline*}
 and, similarly, \(\pbrsa{n}{w}\hp{2n}^\mpi G\rk{w}+\pborsa{n+1}{w}=-\PBha{n}{w}\).
 Consequently, the inequality \(\det\rk{\pbrsa{n}{w}\hp{2n}^\mpi G\rk{w}+\pborsa{n+1}{w}}\neq0\) and the equation
\[\begin{split}
 F\rk{w}
 &=-\ek*{\parsa{n}{w}\hp{2n}^\mpi G\rk{w}+\paorsa{n+1}{w}}\ek*{\pbrsa{n}{w}\hp{2n}^\mpi G\rk{w}+\pborsa{n+1}{w}}^\inv\\
 &=-\PAha{n}{w}\ek*{\PBha{n}{w}}^\inv
 =\mphka{2n+1}{w}
\end{split}\]
 follow.
 In view of \(\seqsh{2n+1}\in\Hggeq{2n+1}\), we can apply \rlem{L3N10} to the sequence \(\seqsh{2n+1}\) to obtain \(\mphk{2n}=\mphk{2n+1}\), where \(\mphk{2n}\) is built according to \rnotap{CB15}{CB15.a} from the sequence \(\seqsh{2n+1}\).
 By virtue of \rdefn{D.nat-ext} and \rrem{ABM-tru}, we have \(\mphk{2n}=\mpk{2n}\).
 Therefore, we can finally conclude \(F\rk{w}=\mpka{2n}{w}\).
\eproof

 For the following designation, we refer to \rnota{CB15}.
 There is a connection with the quantities introduced there.
 Observe that the following constructions are well defined due to \rlemss{BK8.6}{L2412}, regarding \rnota{N1326}:

\bnotal{N2150M}
 Let \(\kappa\in\NOinf\) and let \(\seqska\in\Hggeqka\) with \tsHp{} \(\sHp{\kappa}\), \sabcd{} \(\abcd{\ev{\kappa}}\), and \tsXF{} \(\seqX{\kappa}\).
\benui
 \il{N2150M.a} For all \(n\in\NO\) such that \(2n\leq\kappa\), let \(\lpsdrk{2n},\rpsdrk{2n}\colon\CR\to\Cqq\) be defined by
\begin{align*}
 \lpsdrka{2n}{z}&\defeq\abs{z-\ko z}^\inv\ek*{\pda{n}{z}}^\inv\hp{2n}\ek*{\rk{\rim z}^\inv\rim\XFa{2n}{z}}^\mpi\hp{2n}\ek*{\pda{n}{z}}^\invad
\intertext{and}
 \rpsdrka{2n}{z}&\defeq\abs{z-\ko z}^\inv\ek*{\pba{n}{z}}^\invad\hp{2n}\ek*{\rk{\rim z}^\inv\rim\XFa{2n}{z}}^\mpi\hp{2n}\ek*{\pba{n}{z}}^\inv.
\end{align*}
 \il{N2150M.b} If \(\kappa\geq1\), for all \(n\in\NO\) such that \(2n+1\leq\kappa\), let \(\lpsdrk{2n+1},\rpsdrk{2n+1}\colon\CR\to\Cqq\) be defined by
\begin{align*}
 \lpsdrka{2n+1}{z}&\defeq\abs{z-\ko z}^\inv\ek*{\pda{n}{z}}^\inv\hp{2n}\ek*{\rk{\rim z}^\inv\rim\XFa{2n+1}{z}}^\mpi\hp{2n}\ek*{\pda{n}{z}}^\invad
\intertext{and}
 \rpsdrka{2n+1}{z}&\defeq\abs{z-\ko z}^\inv\ek*{\pba{n}{z}}^\invad\hp{2n}\ek*{\rk{\rim z}^\inv\rim\XFa{2n+1}{z}}^\mpi\hp{2n}\ek*{\pba{n}{z}}^\inv.
\end{align*}
\eenui
\enota

\breml{R228S}
 Let \(\kappa\in\NOinf\), let \(\seqska\in\Hggeqka\), and let \(m\in\mn{0}{\kappa}\).
 Using \eqref{R1019.0} in \rrem{R1019} as well as \rrem{A.R.A+>}, it is readily checked that the matrix-valued functions \(\lrk{m},\rrk{m}\colon\CR\to\Cqq\) given by \rnota{CB15} and \(\lpsdrk{m},\rpsdrk{m}\colon\CR\to\Cqq\) given by \rnota{N2150M} are connected via \(\lpsdrka{m}{z}=\abs{z-\ko z}^\inv\ek{\lrka{m}{z}}\ek{\lrka{m}{z}}^\ad\) and \(\rpsdrka{m}{z}={\abs{z-\ko z}}^\inv\ek{\rrka{m}{z}}^\ad\ek{\rrka{m}{z}}\) for all \(z\in\CR\).
 In particular, \(\lpsdrka{m}{z}\in\Cggq\) and \(\rpsdrka{m}{z}\in\Cggq\) as well as \(\ek{\rk{z-\ko z}^\inv\lrka{m}{z}}\ek{\rk{z-\ko z}^\inv\lrka{m}{z}}^\ad=\abs{z-\ko z}^\inv\sqrt{\lpsdrka{m}{z}}\sqrt{\lpsdrka{m}{z}}^\ad\) and \(\ek{\rrka{m}{z}}^\ad\rrka{m}{z}=\abs{z-\ko z}\sqrt{\rpsdrka{m}{z}}^\ad\sqrt{\rpsdrka{m}{z}}\) hold true for all \(z\in\CR\).
\erem

\breml{LRC-tru}
 Let \(\kappa\in\NOinf\) and let \(\seqska\in\Hggeqka\).
 Regarding \rnota{N2150M}, then one can see from \rremsss{CM2.1.65}{BK8.1}{XF-tru} that for each \(m\in\mn{0}{\kappa}\), the functions \(\lpsdrk{m}\) and \(\rpsdrk{m}\) are built only from the matrices \(\su{0},\su{1},\dotsc,\su{m}\) and does not depend on the matrices \(\su{j}\) with \(j\geq m+1\).
 Moreover, the following \rlem{L0844} even shows that \(\lpsdrk{m}\) and \(\rpsdrk{m}\) are independent of \(\su{j}\) with \(j\geq\eff{m}+1\), where \(\eff{m}\) is given by \eqref{SK5-11}.
\erem
 
\bleml{L0844}
 Let \(\kappa\in\Ninf\), let \(\seqska\in\Hggeqka\), and let \(n\in\NO\) be such that \(2n+1\leq\kappa\).
 Then \(\lpsdrka{2n}{z}=\lpsdrka{2n+1}{z}\) and \(\rpsdrka{2n}{z}=\rpsdrka{2n+1}{z}\) for all \(z\in\CR\). 
\elem
\bproof
 Use \rrem{R228S} and \rlem{C220N}.
\eproof

 Now we cite a result for operator balls due to Yu.~L.~Shmul\cprime{}yan.

\bpropnl{\zitaa{MR0273377}{\cthm{1.3}{868}}}{d92t152}
 Let \(M\in\Cpq\), let \(A_1, A_2\in\Cpp\), and let \(B_1, B_2\in\Cqq\).
 Then the following statements are equivalent:
\baeqi{0}
 \il{d92t152.i} \(\cmb{M}{A_1}{B_1}=\cmb{M}{A_2}{B_2}\).
 \il{d92t152.ii} One of the following two conditions is satisfied:
\begin{enumerate}
 \item Each of the pairs \(\rk{A_1,B_1}\) and \(\rk{A_2,B_2}\) contains a zero matrix.
 \item There exists a positive real number \(\rho\) such that the equations \(A_1A_1^\ad=\rho A_2A_2^\ad\) and \(B_1^\ad B_1=\rho^\inv B_2^\ad B_2\) are satisfied.
\end{enumerate}
\eaeqi
\eprop

 Note that there is a detailed proof of \rprop{d92t152} also in \zitaa{MR1152328}{\cthm{1.5.2}{49}}.
 
 Now we are able to reformulate \rthm{T45CD} with a new representation of the matrix ball in question.

\bthml{F49}
 Let \(n\in\NO\) and let \(\seqs{2n}\in\Hggeq{2n}\).
 For all \(w\in\pip\), then
\[
 \setaca*{F\rk{w}}{F\in\RFOqskg{2n}}
 =\Cmb{\mpka{2n}{w}}{\sqrt{\lpsdrk{2n}\rk{w}}}{\sqrt{\rpsdrk{2n}\rk{w}}}.
\]
\ethm
\bproof
 Regarding \rthm{T45CD} and \rrem{R228S}, the assertion follows by applying \rprop{d92t152}.
\eproof

 Now we derive a connection between the right and left semi-radii of the matrix ball described in \rthm{F49}.
 
\bpropl{L2024}
 Let \(\kappa\in\NOinf\), let \(\seqska\in\Hggeqka\), and let \(m\in\mn{0}{\kappa}\).
 For all \(z\in\CR\), then \(\rpsdrka{m}{z}=\lpsdrka{m}{\ko{z}}\), where \(\lpsdrk{m}\) and \(\rpsdrk{m}\) are given by \rnota{N2150M}.
\eprop
\bproof
 Let \(\sHp{\kappa}\) be the \tsHpo{\(\seqska\)}, let \(\abcd{\ev{\kappa}}\) be the \sabcdo{\(\seqska\)}, and let \(\seqX{\kappa}\) be the \tsXFo{\(\seqska\)}.
 We consider an arbitrary \(z\in\CR\).
 According to \rnota{N2150M}, we have \(\lpsdrka{m}{\ko z}=\abs{\ko z-z}^\inv\ek{\pda{n}{\ko z}}^\inv\hp{2n}\ek{\rk{\rim\ko z}^\inv\rim\XFa{m}{\ko z}}^\mpi\hp{2n}\ek{\pda{n}{\ko z}}^\invad\) and \(\rpsdrka{m}{z}=\abs{z-\ko z}^\inv\ek{\pba{n}{z}}^\invad\hp{2n}\ek{\rk{\rim z}^\inv\rim\XFa{m}{z}}^\mpi\hp{2n}\ek{\pba{n}{z}}^\inv\), where \(n\defeq\ef{m}\) is given by \eqref{SK5-11}.
 Because of \rrem{R1624}, we have \(\su{j}^\ad=\su{j}\) for all \(j\in\mn{0}{\kappa}\).
 Consequently, we can apply \rrem{BK8.2} to get \(\pda{n}{\ko{z}}=\ek{\pba{n}{z}}^\ad\).
 Regarding \rrem{R0735}, \rlem{L0625} provides \(\XFa{m}{\ko{z}}=\ek{\XFa{m}{z}}^\ad\), implying \(\rim\XFa{m}{\ko{z}}=-\rim\XFa{m}{z}\).
 Taking additionally into account \(\rim\ko z=-\rim z\), we consequently perceive
\[\begin{split}
 \abs{\ko z-z}\lpsdrka{m}{\ko{z}}
 &=\ek*{\pda{n}{\ko z}}^\inv\hp{2n}\ek*{\rk{\rim\ko z}^\inv\rim\XFa{m}{\ko z}}^\mpi\hp{2n}\ek*{\pda{n}{\ko z}}^\invad\\
 &=\ek*{\pba{n}{z}}^\invad\hp{2n}\ek*{\rk{\rim z}^\inv\rim\XFa{m}{z}}^\mpi\hp{2n}\ek*{\pba{n}{z}}^\inv
 =\abs{z-\ko z}\rpsdrka{m}{z},
\end{split}\]
 which, in view of \(\abs{\ko z-z}=\abs{z-\ko z}\neq0\) yields the assertion.
\eproof

\section{Some connections to I.~V.~Kovalishina's investigations in the non-degenerate case}\label{S1031}
 In this section, we first sketch the approach developed by I.~V.~Kovalishina \cite{MR703593} who considered the non-degenerate  case of a sequence \(\seqs{2n}\in\Hgq{2n}\).
 Actually, she studied simultaneously the Hamburger moment problem, the matricial Carath\'eodory problem, and the matricial Nevanlinna--Pick problem in the framework of V.~P.~Potapov's method of Fundamental Matrix Inequalities (FMI~method).
 What concerns the computation of the corresponding Weyl matrix balls, she concentrated on the Nevanlinna--Pick case.
 The computation of the Weyl matrix balls for the Hamburger case along the line of I.~V.~Kovalishina was realized by Chen/Hu~\zitaa{MR1740433}{\cpage{77}} and also in an unpublished handwritten manuscript by Yu.~M.~Dyukarev \cite{Dyu04}.
 The bridge to our approach is built by representing the canonical \tqqa{blocks} of the resolvent matrix of Kovalishina in terms of orthogonal matrix polynomials of the first and the second type.
 These \tqqa{matrix} polynomials turn out to be strongly interrelated to our \tabcd{} introduced in \rdefn{K19} in terms of which the parameters of the corresponding Weyl matrix balls are written (see \rnota{CB15}, \rthm{T45CD}).
 In order to recall I.~V.~Kovalishina's description of the set \(\RFOqskg{2n}\), we need some notations.

 Let \(\rv{0}\defeq\Iq\) and let \(\T{0}\defeq\Oqq\).
 Furthermore, for all \(k\in\N\), let \(\rv{k}\defeq\smat{\Iq\\\Ouu{kq}{q}}\) and let \(\T{k}\defeq\smat{\Ouu{q}{kq}&\Oqq\\\Iu{kq}&\Ouu{kq}{q}}\).
 Obviously, for all \(k\in\NO\) and all \(z\in\C\), the matrix \(\Iu{\rk{k+1}q}-z\T{k}\) is invertible.
 In particular, for all \(k\in\NO\), the function \(\RT{k}\colon\C\to\Coo{\rk{k+1}q}{\rk{k+1}q}\) given by \(\RTa{k}{z}\defeq{\rk{\Iu{\rk{k+1}q}-z\T{k}}}^\inv\) is well defined.
 One can easily see, that for all \(z\in\C\) and all \(n\in\NO\), we have
\beql{Rblock}
 \RTa{n}{z}
 =
 \Mat{
  \Iq&\Oqq&\Oqq&\hdots&\Oqq&\Oqq\\
  z\Iq&\Iq&\Oqq&\hdots&\Oqq&\Oqq\\
  z^2\Iq&z\Iq&\Iq&&\Oqq&\Oqq\\
  \vdots&\vdots&\vdots&\ddots&&\vdots\\
  z^{n-1}\Iq&z^{n-2}\Iq&z^{n-3}\Iq&\hdots&\Iq&\Oqq\\
  z^n\Iq&z^{n-1}\Iq&z^{n-2}\Iq&\hdots&z\Iq&\Iq
 }
\eeq
 and, in view of \eqref{NMB}, hence
\begin{align}\label{Rv=E}
 \RTa{n}{z}\rv{n}&=\ek*{\eua{n}{\ko z}}^\ad&
 &\text{and}&
 \rv{n}^\ad\ek*{\RTa{n}{\ko z}}^\ad&=\eua{n}{z}.
\end{align}
 Let \(\kappa\in\NOinf\) and let \(\seqska\) be a sequence of complex \tpqa{matrices}.
 Then let
\begin{align}\label{u}
 \ru{0}&\defeq\Opq&
 &\text{and}&
 \ru{n}
 &\defeq
 \Mat{
  \Opq\\
  -\yuu{0}{n-1}
 }
\end{align}
 for all \(n\in\N\) with \(n-1\leq\kappa\). 
 
\bnotal{N.rU}
 Let \(n\in\NO\) and let \(\seqs{2n}\in\Hgq{2n}\).
 Then let \(\rU{n}\colon\C\to\Coo{2q}{2q}\) be defined by \(\rUa{n}{z}\defeq\Iu{2q}+\iu z\mat{\ru{n},\rv{n}}^\ad\ek*{\RTa{n}{\ko z}}^\ad\Hu{n}^\inv\mat{\ru{n},\rv{n}}\Jimq\).
\enota

 Now we are able to present Kovalishina's description of the set \(\RFOqskg{2n}\).

\bthmnl{\cite{MR703593}}{T1111}
 Let \(n\in\NO\) and let \(\seqs{2n}\in\Hgq{2n}\).
 Further, let \(\smat{U_{n;11}&U_{n;12}\\U_{n;21}&U_{n;22}}\) be the \tqqa{\tbr{}} of the restriction of \(\rU{n}\) onto \(\pip\).
 Then:
\benui
 \il{T1111.a} For each \(\copa{\phi}{\psi}\in\PRFq\), the function \(\det\rk{U_{n;21}\phi+U_{n;22}\psi}\) does not vanish identically and the matrix-valued function \(F\defeq\rk{U_{n;11}\phi+U_{n;12}\psi}\rk{U_{n;21}\phi+U_{n;22}\psi}^\inv\) belongs to \(\RFOqskg{2n}\).
 \il{T1111.b} For each \(F\in\RFOqskg{2n}\), there exists a pair \(\copa{\phi}{\psi}\in\PRFq\) such that \(F=\rk{U_{n;11}\phi+U_{n;12}\psi}\rk{U_{n;21}\phi+U_{n;22}\psi}^\inv\).
 \il{T1111.c} For each \(\ell\in\set{1,2}\) let \(\copa{\phi_\ell}{\psi_\ell}\in\PRFq\) and let \(F_\ell\defeq\rk{U_{n;11}\phi_\ell+U_{n;12}\psi_\ell}\rk{U_{n;21}\phi_\ell+U_{n;22}\psi_\ell}^\inv\).
 Then \(F_1=F_2\) if and only if \(\cpcl{\phi_1}{\psi_1}=\cpcl{\phi_2}{\psi_2}\).
\eenui
\ethm

 The function \(\rU{n}\) introduced in \rnota{N.rU} satisfies a useful functional equation.

\blemnl{\cite{MR703593}}{L1720} 
 Let \(n\in\NO\) and let \(\seqs{2n}\in\Hgq{2n}\).
 Then the matrix-valued function \(\rU{n}\) given by \rnota{N.rU} is a matrix polynomial of degree at most \(n+1\) and the following statements hold true:
\benui
 \il{L1720.a} For every choice of \(w,z\in\C\),
\[
 \Jimq-\rUa{n}{z}\Jimq\ek*{\rUa{n}{w}}^\ad
 =\iu\rk{\ko w-z}\mat{\ru{n},\rv{n}}^\ad\ek*{\RTa{n}{\ko z}}^\ad\Hu{n}^\inv\RTa{n}{\ko w}\mat{\ru{n},\rv{n}}.
\]
 \il{L1720.b} For \(w\in\pip\), the matrix \(\Jimq-\rUa{n}{w}\Jimq\ek*{\rUa{n}{w}}^\ad\) is \tnnH{}.
 In particular, \(\Jimq-\rUa{n}{x}\Jimq\ek*{\rUa{n}{x}}^\ad=\Ouu{2q}{2q}\) is fulfilled for all \(x\in\R\).
 \il{L1720.c} For all \(z\in\C\), the matrix \(\rUa{n}{z}\) is non-singular and \(\ek{\rUa{n}{z}}^\inv=\Jimq\ek{\rUa{n}{\ko z}}^\ad\Jimq\) holds true.
 In particular, \(\rU{n}^\inv\) is a matrix polynomial of degree at most \(n+1\), which satisfies the representation \(\ek{\rUa{n}{z}}^\inv=\Iu{2q}-\iu z\mat{\ru{n},\rv{n}}^\ad\Hu{n}^\inv\RTa{n}{z}\mat{\ru{n},\rv{n}}\Jimq\) for all \(z\in\C\).
\eenui
\elem

 Note that a detailed proof of \rlem{L1720} is given, \teg, in \zitaa{Roe03}{\clem{8.7}{89}}.
 
\blemnl{\cite{MR703593,MR1740433,Dyu04}}{L1737}
 Under the assumptions of \rlem{L1720}, one can see from \rlem{L1720} that \(\rW{n}\colon\C\to\Coo{2q}{2q}\) defined by \(\rWa{n}{z}\defeq\ek{\rUa{n}{z}}^\invad\rk{-\Jimq}\ek{\rUa{n}{z}}^\inv\) fulfills \(\det\rWa{n}{z}\neq0\) and \(\ek{\rWa{n}{z}}^\inv=\Jimq\rWa{n}{\ko z}\Jimq\) for all \(z\in\C\).
 Moreover, \rlem{L1720} shows that, for all \(z\in\C\), the \tqqa{\tbr}
\[
 \rWa{n}{z}
 =
 \Mat{
  -\rRa{n}{z}&\rSa{n}{z}\\
  \ek{\rSa{n}{z}}^\ad&-\rTa{n}{z}
 }
\]
 holds true, where
\begin{align*}
 \rRa{n}{z}&\defeq\iu\rk{\ko z-z}\rv{n}^\ad\ek*{\RTa{n}{z}}^\ad\Hu{n}^\inv\ek*{\RTa{n}{z}}\rv{n},\\
 \rSa{n}{z}&\defeq\iu\Iq+\iu\rk{\ko z-z}\rv{n}^\ad\ek*{\RTa{n}{z}}^\ad\Hu{n}^\inv\ek*{\RTa{n}{z}}\ru{n},
\intertext{and}
 \rTa{n}{z}&\defeq\iu\rk{\ko z-z}\ru{n}^\ad\ek*{\RTa{n}{z}}^\ad\Hu{n}^\inv\ek*{\RTa{n}{z}}\ru{n}. 
\end{align*}
 In particular, for all \(w\in\pip\), the matrices \(\rRa{n}{w}\) and \(-\rRa{n}{\ko w}\) are \tpH{}, where \(-\ek{\rRa{n}{\ko w}}^\inv=\ek{\rSa{n}{w}}^\ad\ek{\rRa{n}{w}}^\inv\rSa{n}{w}-\rTa{n}{w}\).
\elem
 
 Note that a detailed proof of \rlem{L1737} is stated in \zitaa{Tau14}{\clemss{10.6}{265}{10.8}{269}}.

\bnotal{N.rCdg}
 Let \(n\in\NO\) and let \(\seqs{2n}\in\Hgq{2n}\).
 Then let \(\rC{n},\rg{n},\rd{n}\colon\pip\to\Coo{q}{q}\) be defined by
\[
 \rCa{n}{w}
 \defeq\ek*{\rRa{n}{w}}^\inv\rSa{n}{w}
\]
 and by
\begin{align*}
 \rga{n}{w}&\defeq\ek*{\rRa{n}{w}}^\inv&
&\text{and}&
 \rda{n}{w}&\defeq\ek*{\rSa{n}{w}}^\ad\ek*{\rRa{n}{w}}^\inv\rSa{n}{w}-\rTa{n}{w}.
\end{align*}
\enota

\breml{R1542}
 Let \(n\in\NO\) and let \(\seqs{2n}\in\Hgq{2n}\).
 According to \rlem{L1737}, for all \(w\in\pip\), then \(\rga{n}{w},\rda{n}{w}\in\Cgq\).
\erem

\breml{R9.10}
 Under the assumptions of \rlem{L1720}, from \rlem{L1737}, for all \(z\in\C\), we conclude
\[
 \rWa{n}{z}\Jimq\rWa{n}{\ko z}
 =\rWa{n}{z}\Jimq\rWa{n}{\ko z}\Jimq^2
 =\rWa{n}{z}\ek*{\rWa{n}{z}}^\inv\Jimq
 =\Jimq
\]
 and, in view of \eqref{JQ}, in particular
\[\begin{split}
 \Oqq
 &=\mat{\Iq,\Oqq}\Jimq\Mat{\Iq\\\Oqq}
 =\mat{\Iq,\Oqq}\rWa{n}{z}\Jimq\rWa{n}{\ko z}\Mat{\Iq\\\Oqq}\\
 &=\iu\rk*{\rRa{n}{z}\ek*{\rSa{n}{\ko z}}^\ad-\rSa{n}{z}\rRa{n}{\ko z}},
\end{split}\]
 which implies \(\rCa{n}{w}=\ek{\rRa{n}{w}}^\inv\rSa{n}{w}=\ek{\rSa{n}{\ko w}}^\ad\ek{\rRa{n}{\ko w}}^\inv\) for all \(w\in\pip\).
\erem

 Using \rnota{N.rCdg} and \rrem{R1542}, now we are able for \(w\in\pip\) to describe the set \(\setaca{F\rk{w}}{\RFOqskg{2n}}\) as a matrix ball.
 
\bthmnl{\cite{MR703593,MR1740433,Dyu04}}{T1546}
 Let \(n\in\NO\) and let \(\seqs{2n}\in\Hgq{2n}\).
 For all \(w\in\pip\), then
\[
 \setaca*{F\rk{w}}{F\in\RFOqskg{2n}}
 =\Cmb{\rCa{n}{w}}{\sqrt{\rga{n}{w}}}{\sqrt{\rda{n}{w}}}.
\]
\ethm

 From \rthm{T1546}, which is dedicated to the non-degenerate case, \rthm{F49}, where the general case is considered, and \rprop{d92t152} it is clear that there is an explicit connection between the parameters of the matrix balls in question, where the centers coincide necessarily by a result due to Shmul\cprime{}yan \zitaa{MR0273377}{\csec{1}{}} (see also \zitaa{MR1152328}{\ccor{1.5.1}{47}}).
 The main aim of this section is now to verify that the parameters of the matrix ball occurring in \rthm{T1546} coincide even with the matrices occurring in \rthm{F49}.
 This requires some considerations.
 
 Let \(\kappa\in\NOinf\) and let \(\seqska\in\Hggeqka\).
 Then, in view of \rnota{B.N.deg}, for all \(n\in\NO\) with \(2n-1\leq\kappa\), let
\beql{yn}
 \yu{n}
 \defeq\cvuo{n}{\pb{n}}.
\eeq
 Regarding \rrem{SK32R}, in particular \(\yu{0}=\Iq\).
 Because of \rremss{SK32R}{B.R.P=euY}, \eqref{Rv=E}, and \eqref{yn}, for all \(n\in\NO\) such that \(2n-1\leq\kappa\) and all \(z\in\C\), we can infer
\beql{PbR}
 \pba{n}{z}
 =\rv{n}^\ad\ek*{\RTa{n}{\ko z}}^\ad\yu{n}
\eeq
 and, using additionally \rremss{R1624}{BK8.2}, hence
\beql{PdR}
 \pda{n}{z}
 =\yu{n}^\ad\RTa{n}{z}\rv{n}.
\eeq

\bleml{L0952}
 Let \(\kappa\in\NOinf\), let \(\seqska\in\Hggeqka\), and let \(\abcd{\ev{\kappa}}\) be the \sabcdo{\(\seqska\)}.
 For all \(k\in\NO\) with \(2k-1\leq\kappa\), then \(\paa{k}{z}=-\ru{k}^\ad\ek{\RTa{k}{\ko z}}^\ad\yu{k}\) and \(\pca{k}{z}=-\yu{k}^\ad\RTa{k}{z}\ru{k}\) for all \(z\in\C\).
\elem
\bproof
 Let \(k\in\NO\) with \(2k-1\leq\kappa\) and let \(z\in\C\).
 In view of \rprop{P1203}, then \(\pa{k}=\pb{k}^\secra{s}\). 
 Using \eqref{u}, \rrem{SK32R}, \rnota{B.N.sec}, and \eqref{K19.0}, in the case \(k=0\) we easily obtain then
 \[
  -\ru{0}^\ad\ek*{\RTa{0}{\ko z}}^\ad\yu{0}
  =\Oqq
  =\pb{0}^\secra{s}(z)
  =\paa{0}{z}.
 \]
 Now suppose \(k\geq1\).
 In view of \eqref{yz} and \rrem{R1624}, we have \(\yuu{0}{k-1}^\ad=\zuu{0}{k-1}\).
 Using \eqref{yz}, \eqref{Rblock}, \eqref{NMB}, and \eqref{M.N.S}, we conclude
 \[
  \zuu{0}{k-1}\ek*{\RTa{k-1}{\ko z}}^\ad
  =\row\Seq{\sum_{\ell=0}^jz^\ell\su{j-\ell}}{j}{0}{k-1}
  =\eua{k-1}{z}\SUu{k-1}.
 \]
 From \rrem{SK32R}, \rnota{B.N.sec}, and \eqref{yn} we know
 \(
  \pb{k}^\secra{s}(z)
  =\eua{k-1}{z}\mat{\Ouu{kq}{q},\SUu{k-1}}\yu{k}
 \).
 Taking additionally into account \eqref{u} and \eqref{Rblock}, then
 \[
  \begin{split}
   -\ru{k}^\ad\ek*{\RTa{k}{\ko z}}^\ad\yu{k}
   &=-\mat{\Oqq,-\zuu{0}{k-1}}
   \Mat{
    \uk&\uk\\
    \Ouu{kq}{q}&\ek{\RTa{k-1}{\ko z}}^\ad
   }\yu{k}\\
   &=\mat*{\Oqq,\zuu{0}{k-1}\ek{\RTa{k-1}{\ko z}}^\ad}\yu{k}   
   =\mat*{\Oqq,\eua{k-1}{z}\SUu{k-1}}\yu{k}\\
   &=\eua{k-1}{z}\mat{\Ouu{kq}{q},\SUu{k-1}}\yu{k}
   =\pb{k}^\secra{s}(z)
   =\paa{k}{z}
  \end{split}
 \]
 follows. 
 Using \rremss{R1624}{BK8.2}, we can easily get the second identity from the first one.
\eproof

\breml{N9.3}
 Let \(\kappa\in\Ninf\) and let \(\seqs{2\kappa}\in\Hggq{2\kappa}\).
 In view of \rpropss{143.T1336}{B.P.oMP-lGs} and the notation given in \eqref{yn}, for all \(n\in\mn{1}{\kappa}\), we have \(\yu{n}=\tmat{-\Xu{n}\\\Iq}\) with some matrix \(\Xu{n}\) belonging to \(\tu{n}\defeq\setaca{r\in\Coo{nq}{q}}{\Hu{n-1}r=\yuu{n}{2n-1}}\).
 If \(\seqs{2\kappa}\in\Hgq{2\kappa}\), then, for all \(n\in\mn{1}{\kappa}\), the matrix \(\Hu{n-1}\) is non-singular and \(\yu{n}=\rY{n}\) where
\beql{YY}
 \rY{n}
 \defeq
 \Mat{
  -\Hu{n-1}^\inv\yuu{n}{2n-1}\\
  \Iq
 }.
\eeq
\erem

 Setting \(\rY{0}\defeq\Iq\) and using \eqref{YY}, we see from \rprop{P1644} that the following construction is well defined:

\bnotal{N.pn}
 Let \(\kappa\in\NOinf\) and let \(\seqs{2\kappa}\in\Hgq{2\kappa}\) with \tsHp{} \(\sHp{2\kappa}\).
 Then, for all \(k\in\mn{0}{\kappa}\), let \(\pna{k}, \pnb{k},\pnc{k}, \pnd{k}\colon\C\to\Cqq\) be defined by
\begin{align*}
 \pnaa{k}{z}&\defeq-\ru{k}^\ad\ek*{\RTa{k}{\ko z}}^\ad\rY{k}\sqrt{\hp{2k}}^\inv,&
 \pnba{k}{z}&\defeq\rv{k}^\ad\ek*{\RTa{k}{\ko z}}^\ad\rY{k}\sqrt{\hp{2k}}^\inv,\\
 \pnca{k}{z}&\defeq-\sqrt{\hp{2k}}^\inv\rY{k}^\ad\RTa{k}{z}\ru{k},&
 \pnda{k}{z}&\defeq\sqrt{\hp{2k}}^\inv\rY{k}^\ad\RTa{k}{z}\rv{k}.
\end{align*}
\enota

\breml{R48lr}
 Under the assumption of \rnota{N.pn}, for every choice of \(k\in\mn{0}{\kappa}\) and \(z\in\C\), we have obviously \(\pnca{k}{z}=\ek{\pnaa{k}{\ko z}}^\ad\) and \(\pnda{k}{z}=\ek{\pnba{k}{\ko z}}^\ad\).
\erem

\bleml{R1704}
 Let \(\kappa\in\NOinf\) and let \(\seqs{2\kappa}\in\Hgq{2\kappa}\) with \tsHp{} \(\sHp{2\kappa}\) and \sabcd{} \(\abcd{\kappa}\).
 For all \(k\in\mn{0}{\kappa}\), then
\begin{align*}
 \pna{k}\sqrt{\hp{2k}}&=\pa{k},&
 \pnb{k}\sqrt{\hp{2k}}&=\pb{k},&
 \sqrt{\hp{2k}}\pnc{k}&=\pc{k},&
 \sqrt{\hp{2k}}\pnd{k}&=\pd{k}.
\end{align*}
\elem
\bproof
 In view of \rnota{N.pn} and \rrem{N9.3}, the assertion follows from \eqref{PbR}, \eqref{PdR}, and \rlem{L0952}.
\eproof

 An important step in the approach due to Kovalishina is now to write the matrices introduced in \rlem{L1737} in terms of the matrix polynomials from \rnota{N.pn}

\bpropnl{\cite{MR703593}}{P1550}
 Let \(n\in\NO\) and let \(\seqs{2n}\in\Hgq{2n}\).
 For all \(z\in\C\), then
\begin{align*}
 \rRa{n}{z}&=\iu\rk{\ko z-z}\sum_{k=0}^n\ek*{\pnda{k}{z}}^\ad\pnda{k}{z},&
 \rSa{n}{z}&=\iu\Iq-\iu\rk{\ko z-z}\sum_{k=0}^n\ek*{\pnda{k}{z}}^\ad\pnca{k}{z},
\end{align*}
 and
\[
 \rTa{n}{z}
 =\iu\rk{\ko z-z}\sum_{k=0}^n\ek*{\pnca{k}{z}}^\ad\pnca{k}{z}.
\]
\eprop

 Note that detailed proofs of \rprop{P1550} are given in \cite{Dyu04} and \zitaa{Tau14}{\cbem{10.16}{280}}.
 As a consequence of \rprop{P1550} and \rrem{R48lr} we will express the parameters of the matrix ball occurring in \rthm{T1546} in terms of the matrix polynomials from \rnota{N.pn}.

\bcorl{C1645}
 Let \(n\in\NO\) and let \(\seqs{2n}\in\Hgq{2n}\).
 For all \(w\in\pip\), then
\begin{align*}
 \rCa{n}{w}&=\rk*{\iu\rk{\ko w-w}\sum_{k=0}^n\ek*{\pnda{k}{w}}^\ad\pnda{k}{w}}^\inv\rk*{\iu\Iq-\iu\rk{\ko w-w}\sum_{k=0}^n\ek*{\pnda{k}{w}}^\ad\pnca{k}{w}},\\
 \rCa{n}{w}&=\rk*{\iu\Iq-\iu\rk{\ko w-w}\sum_{k=0}^n\pnaa{k}{w}\ek*{\pnba{k}{w}}^\ad}\rk*{\iu\rk{\ko w-w}\sum_{k=0}^n\pnba{k}{w}\ek*{\pnba{k}{w}}^\ad}^\inv
\end{align*}
 and
\begin{align*}
 \rga{n}{w}&=\rk*{\iu\rk{\ko w-w}\sum_{k=0}^n\ek*{\pnda{k}{w}}^\ad\pnda{k}{w}}^\inv,&
 \rda{n}{w}&=\rk*{\iu\rk{\ko w-w}\sum_{k=0}^n\pnba{k}{w}\ek*{\pnba{k}{w}}^\ad}^\inv.
\end{align*}
\ecor

 A proof of \rcor{C1645} can be easily obtained using \rlem{L1737} and \rremss{R9.10}{R48lr}.
 We omit the details.
 The final considerations are aimed to look at above investigations under the view of \rsec{Cha13}.

\bleml{L0824}
 Let \(\kappa\in\Ninf\), let \(\seqska\in\Hgeqka\) with \tsHp{} \(\sHp{\kappa}\) and \sabcd{} \(\abcd{\ev{\kappa}}\), and let \(n\in\NO\) be such that \(2n+1\leq\kappa\).
 Then \(\det\hp{2n}\neq0\) and, for each \(z\in\CR\), moreover
\begin{align}
 \PAa{n}{z}\hp{2n}^\inv\ek*{\pba{n}{z}}^\ad&=\paa{n}{z}\hp{2n}^\inv\ek*{\pba{n+1}{z}}^\ad-\paa{n+1}{z}\hp{2n}^\inv\ek*{\pba{n}{z}}^\ad,\label{K1}\\
 \PBa{n}{z}\hp{2n}^\inv\ek*{\pba{n}{z}}^\ad&=\pba{n}{z}\hp{2n}^\inv\ek*{\pba{n+1}{z}}^\ad-\pba{n+1}{z}\hp{2n}^\inv\ek*{\pba{n}{z}}^\ad,\notag\\
 \ek*{\pda{n}{z}}^\ad\hp{2n}^\inv\PCa{n}{z}&=\ek*{\pda{n+1}{z}}^\ad\hp{2n}^\inv\pca{n}{z}-\ek*{\pda{n}{z}}^\ad\hp{2n}^\inv\pca{n+1}{z},\notag\\
 \ek*{\pda{n}{z}}^\ad\hp{2n}^\inv\PDa{n}{z}&=\ek*{\pda{n+1}{z}}^\ad\hp{2n}^\inv\pda{n}{z}-\ek*{\pda{n}{z}}^\ad\hp{2n}^\inv\pda{n+1}{z}.\notag
\end{align}
\elem
\bproof
 First observe that \rrem{Schr2.7} yields \(\seqska\in\Hggeqka\).
 Clearly, we have \(\seqs{2n}\in\Hgq{2n}\).
 From \rprop{P1644} we see then \(\hp{2n}\in\Cgq\).
 In particular, \(\hp{2n}^\ad=\hp{2n}\) and \(\det\hp{2n}\neq0\).
 Taking into account \rremp{R52SK}{R52SK.a} and \(\hp{2n}^\ad=\hp{2n}\), we can infer \(\det\pba{n}{z}\neq0\) and \(\ek{\XFa{2n+1}{z}}^\ad=\ek{\pba{n+1}{z}}^\ad\ek{\pba{n}{z}}^\invad\hp{2n}\).
 Consequently, \(\ek{\XFa{2n+1}{z}}^\ad\hp{2n}^\inv\ek*{\pba{n}{z}}^\ad=\ek{\pba{n+1}{z}}^\ad\).
 Regarding additionally \rnota{N1326} and \(\det\hp{2n}\neq0\), we can conclude then
\[\begin{split}
 \PAa{n}{z}\hp{2n}^\inv\ek*{\pba{n}{z}}^\ad
 &=\rk*{\paa{n}{z}\hp{2n}^\inv\ek*{\XFa{2n+1}{z}}^\ad-\paa{n+1}{z}}\hp{2n}^\inv\ek*{\pba{n}{z}}^\ad\\
 &=\paa{n}{z}\hp{2n}^\inv\ek*{\pba{n+1}{z}}^\ad-\paa{n+1}{z}\hp{2n}^\inv\ek*{\pba{n}{z}}^\ad
\end{split}\]
 and, analogously, \(\PBa{n}{z}\hp{2n}^\inv\ek{\pba{n}{z}}^\ad=\pba{n}{z}\hp{2n}^\inv\ek{\pba{n+1}{z}}^\ad-\pba{n+1}{z}\hp{2n}^\inv\ek{\pba{n}{z}}^\ad\).
 Using \rpropp{K17-5}{K17-5.a} instead of \rremp{R52SK}{R52SK.a}, the remaining equations can be obtained similarly.
\eproof

\bleml{L2038}
 Let \(\kappa\in\minf{3}\cup\set{\infi}\), let \(\seqska\in\Hgeqka\) with \tsHp{} \(\sHp{\kappa}\) and \sabcd{} \(\abcd{\ev{\kappa}}\), and let \(n\in\N\) be such that \(2n+1\leq\kappa\).
 For all \(z\in\CR\), then
\begin{align*}
 \PAa{n}{z}\hp{2n}^\inv\ek*{\pba{n}{z}}^\ad&=\rk{\ko z-z}\paa{n}{z}\hp{2n}^\inv\ek*{\pba{n}{z}}^\ad+\PAa{n-1}{z}\hp{2n-2}^\inv\ek*{\pba{n-1}{z}}^\ad,\\
 \PBa{n}{z}\hp{2n}^\inv\ek*{\pba{n}{z}}^\ad&=\rk{\ko z-z}\pba{n}{z}\hp{2n}^\inv\ek*{\pba{n}{z}}^\ad+\PBa{n-1}{z}\hp{2n-2}^\inv\ek*{\pba{n-1}{z}}^\ad,\\
 \ek*{\pda{n}{z}}^\ad\hp{2n}^\inv\PCa{n}{z}&=\rk{\ko z-z}\ek*{\pda{n}{z}}^\ad\hp{2n}^\inv\pca{n}{z}+\ek*{\pda{n-1}{z}}^\ad\hp{2n-2}^\inv\PCa{n-1}{z},\\
 \ek*{\pda{n}{z}}^\ad\hp{2n}^\inv\PDa{n}{z}&=\rk{\ko z-z}\ek*{\pda{n}{z}}^\ad\hp{2n}^\inv\pda{n}{z}+\ek*{\pda{n-1}{z}}^\ad\hp{2n-2}^\inv\PDa{n-1}{z}.
\end{align*}
\elem
\bproof
 First observe that \rrem{Schr2.7} yields \(\seqska\in\Hggeqka\).
 \rrem{R1019} then shows \(\hp{j}\in\CHq\) for all \(j\in\mn{0}{\kappa}\).
 Clearly, we have \(\seqs{2n}\in\Hgq{2n}\).
 From \rprop{P1644} we see then \(\hp{2n},\hp{2n-2}\in\Cgq\).
 In particular, \(\det\hp{2n}\neq0\) and \(\det\hp{2n-2}\neq0\).
 Let \(z\in\CR\).
 Taking into account \rdefn{K19} and \(\hp{2n+1},\hp{2n},\hp{2n-2}\in\CHq\) as well as \(\det\hp{2n}\neq0\) and \(\det\hp{2n-2}\neq0\), we can conclude \(\paa{n+1}{z}=\paa{n}{z}\rk{z\Iq-\hp{2n}^\inv\hp{2n+1}}-\paa{n-1}{z}\hp{2n-2}^\inv\hp{2n}\) and \(\ek{\pba{n+1}{z}}^\ad=\rk{\ko z\Iq-\hp{2n+1}\hp{2n}^\inv}\ek{\pba{n}{z}}^\ad-\hp{2n}\hp{2n-2}^\inv\ek{\pba{n-1}{z}}^\ad\).
 Consequently,
\[
 \paa{n+1}{z}\hp{2n}^\inv\ek*{\pba{n}{z}}^\ad
 =z\paa{n}{z}\hp{2n}^\inv\ek*{\pba{n}{z}}^\ad-\paa{n}{z}\hp{2n}^\inv\hp{2n+1}\hp{2n}^\inv\ek*{\pba{n}{z}}^\ad-\paa{n-1}{z}\hp{2n-2}^\inv\ek*{\pba{n}{z}}^\ad
\]
 and
\[
 \paa{n}{z}\hp{2n}^\inv\ek*{\pba{n+1}{z}}^\ad
 =\ko z\paa{n}{z}\hp{2n}^\inv\ek*{\pba{n}{z}}^\ad-\paa{n}{z}\hp{2n}^\inv\hp{2n+1}\hp{2n}^\inv\ek*{\pba{n}{z}}^\ad-\paa{n}{z}\hp{2n-2}^\inv\ek*{\pba{n-1}{z}}^\ad
\]
 follow.
 Using \rlem{L0824}, we can infer that \eqref{K1} and
\[
 \PAa{n-1}{z}\hp{2n-2}^\inv\ek*{\pba{n-1}{z}}^\ad
 =\paa{n-1}{z}\hp{2n-2}^\inv\ek*{\pba{n}{z}}^\ad-\paa{n}{z}\hp{2n-2}^\inv\ek*{\pba{n-1}{z}}^\ad
\]
 hold true.
 Thus, we obtain
\[\begin{split}
 \PAa{n}{z}\hp{2n}^\inv\ek*{\pba{n}{z}}^\ad
 &=\rk*{\ko z\paa{n}{z}\hp{2n}^\inv\ek*{\pba{n}{z}}^\ad-\paa{n}{z}\hp{2n}^\inv\hp{2n+1}\hp{2n}^\inv\ek*{\pba{n}{z}}^\ad-\paa{n}{z}\hp{2n-2}^\inv\ek*{\pba{n-1}{z}}^\ad}\\
 &\;-\rk*{z\paa{n}{z}\hp{2n}^\inv\ek*{\pba{n}{z}}^\ad-\paa{n}{z}\hp{2n}^\inv\hp{2n+1}\hp{2n}^\inv\ek*{\pba{n}{z}}^\ad-\paa{n-1}{z}\hp{2n-2}^\inv\ek*{\pba{n}{z}}^\ad}\\
 &=\rk{\ko z-z}\paa{n}{z}\hp{2n}^\inv\ek*{\pba{n}{z}}^\ad-\paa{n}{z}\hp{2n-2}^\inv\ek*{\pba{n-1}{z}}^\ad+\paa{n-1}{z}\hp{2n-2}^\inv\ek*{\pba{n}{z}}^\ad\\
 &=\rk{\ko z-z}\paa{n}{z}\hp{2n}^\inv\ek*{\pba{n}{z}}^\ad+\PAa{n-1}{z}\hp{2n-2}^\inv\ek*{\pba{n-1}{z}}^\ad.
\end{split}\]
 In view of \rdefn{K19} and \rlem{L0824} the remaining assertions can be obtained similarly.
\eproof
 
\bleml{L0609}
 Let \(\kappa\in\Ninf\) and let \(\seqska\in\Hgeqka\) with \tsHp{} \(\sHp{\kappa}\) and \sabcd{} \(\abcd{\ev{\kappa}}\).
 For all \(n\in\NO\) such that \(2n+1\leq\kappa\) and all \(z\in\CR\), then
\begin{align*}
 \PAa{n}{z}\hp{2n}^\inv\ek*{\pba{n}{z}}^\ad&=-\Iq+\rk{\ko z-z}\sum_{k=0}^n\paa{k}{z}\hp{2k}^\inv\ek*{\pba{k}{z}}^\ad,\\
 \PBa{n}{z}\hp{2n}^\inv\ek*{\pba{n}{z}}^\ad&=\rk{\ko z-z}\sum_{k=0}^n\pba{k}{z}\hp{2k}^\inv\ek*{\pba{k}{z}}^\ad,\\
 \ek*{\pda{n}{z}}^\ad\hp{2n}^\inv\PCa{n}{z}&=-\Iq+\rk{\ko z-z}\sum_{k=0}^n\ek*{\pda{k}{z}}^\ad\hp{2k}^\inv\pca{k}{z},
\intertext{and}
 \ek*{\pda{n}{z}}^\ad\hp{2n}^\inv\PDa{n}{z}&=\rk{\ko z-z}\sum_{k=0}^n\ek*{\pda{k}{z}}^\ad\hp{2k}^\inv\pda{k}{z}.
\end{align*}
\elem
\bproof
 First observe that \rrem{Schr2.7} yields \(\seqska\in\Hggeqka\).
 Clearly, we have \(\seqs{0}\in\Hgq{0}\).
 From \rprop{P1644} we get then \(\hp{0}\in\Cgq\).
 In particular, \(\det\hp{0}\neq0\).
 Let \(z\in\CR\).
 By virtue of \eqref{ABCD0}, we can conclude \(\PAa{0}{z}=-\hp{0}\), \(\PBa{0}{z}=\rk{\ko z-z}\Iq\), \(\PCa{0}{z}=-\hp{0}\), and \(\PDa{0}{z}=\rk{\ko z-z}\Iq\).
 Taking additionally into account \eqref{K19.0}, then \(\PAa{0}{z}\hp{0}^\inv\ek{\pba{0}{z}}^\ad=-\Iq+\rk{\ko z-z}\paa{0}{z}\hp{0}^\inv\ek{\pba{0}{z}}^\ad\) and \(\PBa{0}{z}\hp{0}^\inv\ek{\pba{0}{z}}^\ad=\rk{\ko z-z}\pba{0}{z}\hp{0}^\inv\ek{\pba{0}{z}}^\ad\) as well as \(\ek{\pda{0}{z}}^\ad\hp{0}^\inv\PCa{0}{z}=-\Iq+\rk{\ko z-z}\ek{\pda{0}{z}}^\ad\hp{0}^\inv\pca{0}{z}\) and \(\ek{\pda{0}{z}}^\ad\hp{0}^\inv\PDa{0}{z}=\rk{\ko z-z}\ek{\pda{0}{z}}^\ad\hp{0}^\inv\pda{0}{z}\) 
 follow.
 Combining these equations with \rlem{L2038} completes the proof.
\eproof

\bpropl{L0926}
 Let \(\kappa\in\Ninf\), let \(\seqska\in\Hgeqka\) with \tsHp{} \(\sHp{\kappa}\) and \sabcd{} \(\abcd{\ev{\kappa}}\), let \(n\in\NO\) be such that \(2n+1\leq\kappa\), and let \(z\in\CR\).
 Then the matrices \(\lpsdrka{2n+1}{z}\) and \(\rpsdrka{2n+1}{z}\) are invertible and fulfill
\begin{align}
 \ek*{\lpsdrka{2n+1}{z}}^\inv&=\frac{\abs{z-\ko z}}{z-\ko z}\rk*{\ek*{\pda{n}{z}}^\ad\hp{2n}^\inv\pda{n+1}{z}-\ek*{\pda{n+1}{z}}^\ad\hp{2n}^\inv\pda{n}{z}}\label{L0926.1}
\intertext{and}
 \ek*{\rpsdrka{2n+1}{z}}^\inv&=\frac{\abs{z-\ko z}}{z-\ko z}\rk*{\pba{n+1}{z}\hp{2n}^\inv\ek*{\pba{n}{z}}^\ad-\pba{n}{z}\hp{2n}^\inv\ek*{\pba{n+1}{z}}^\ad}\label{L0926.2}.
\end{align}
\eprop
\bproof
 First observe that \rrem{Schr2.7} yields \(\seqska\in\Hggeqka\).
 Clearly, we have \(\seqs{2n}\in\Hgq{2n}\).
 From \rprop{P1644} we see then \(\hp{2n}\in\Cgq\).
 In particular, \(\hp{2n}^\ad=\hp{2n}\) and \(\det\hp{2n}\neq0\).
 In view of \(\det\hp{2n}\neq0\) and \eqref{SK5-11}, we can infer from \rprop{L1347} that \(\det\rim\XFa{2n+1}{z}\neq0\).
 Regarding additionally \rnotap{N2150M}{N2150M.b} and \(\det\hp{2n}\neq0\), we can then conclude \(\det\lpsdrka{2n+1}{z}\neq0\) and
\[
 \ek*{\lpsdrka{2n+1}{z}}^\inv
 =\abs{z-\ko z}\ek*{\pda{n}{z}}^\ad\hp{2n}^\inv\ek*{\rk{\rim z}^\inv\rim\XFa{2n+1}{z}}\hp{2n}^\inv\pda{n}{z}
\]
 Taking into account \rpropp{K17-5}{K17-5.a} and \(\hp{2n}^\ad=\hp{2n}\), we can infer \(\det\pda{n}{z}\neq0\) as well as \eqref{GR2} and \(\ek{\XFa{2n+1}{z}}^\ad=\hp{2n}\ek{\pda{n}{z}}^\invad\ek{\pda{n+1}{z}}^\ad\).
 Hence, using additionally \(\rim\XFa{2n+1}{z}=\frac{1}{2\iu}\rk{\XFa{2n+1}{z}-\ek{\XFa{2n+1}{z}}^\ad}\), we obtain \eqref{L0926.1}.
 Applying \rremp{R52SK}{R52SK.a} instead of \rpropp{K17-5}{K17-5.a}, inequality \(\det\rpsdrka{2n+1}{z}\neq0\) and equation \eqref{L0926.2} can be obtained similarly.
\eproof

\bpropl{L0627}
 Let \(\kappa\in\Ninf\), let \(\seqska\in\Hgeqka\) with \tsHp{} \(\sHp{\kappa}\) and \sabcd{} \(\abcd{\ev{\kappa}}\), let \(n\in\NO\) be such that \(2n+1\leq\kappa\), and let \(z\in\CR\).
 Then the matrices \(\lpsdrka{2n+1}{z}\) and \(\rpsdrka{2n+1}{z}\) are invertible and fulfill
\begin{align*}
 \ek*{\lpsdrka{2n+1}{z}}^\inv&=\abs{z-\ko z}\sum_{k=0}^n\ek*{\pda{k}{z}}^\ad\hp{2k}^\inv\pda{k}{z}&
&\text{and}&
 \ek*{\rpsdrka{2n+1}{z}}^\inv&=\abs{z-\ko z}\sum_{k=0}^n\pba{k}{z}\hp{2k}^\inv\ek*{\pba{k}{z}}^\ad.
\end{align*}
\eprop
\bproof
 Using \rprop{L0926} and \rlemss{L0824}{L0609}, we can infer \(\det\lpsdrka{2n+1}{z}\neq0\) and
\[\begin{split}
 \ek*{\lpsdrka{2n+1}{z}}^\inv
 &=\frac{\abs{z-\ko z}}{z-\ko z}\rk*{\ek*{\pda{n}{z}}^\ad\hp{2n}^\inv\pda{n+1}{z}-\ek*{\pda{n+1}{z}}^\ad\hp{2n}^\inv\pda{n}{z}}\\
 &=-\frac{\abs{z-\ko z}}{z-\ko z}\ek*{\pda{n}{z}}^\ad\hp{2n}^\inv\PDa{n}{z}
 =\abs{z-\ko z}\sum_{k=0}^n\ek*{\pda{k}{z}}^\ad\hp{2k}^\inv\pda{k}{z}.
\end{split}\]
 The remaining assertions can be obtained similarly.
\eproof

\bpropl{P1229}
 Let \(n\in\NO\), let \(\seqs{2n}\in\Hgq{2n}\) with \tsHp{} \(\sHp{2n}\) and \sabcd{} \(\abcd{n}\), and let \(z\in\CR\).
 Then the matrices \(\lpsdrka{2n}{z}\) and \(\rpsdrka{2n}{z}\) are invertible and fulfill
\begin{align*}
 \ek*{\lpsdrka{2n}{z}}^\inv&=\abs{z-\ko z}\sum_{k=0}^n\ek*{\pda{k}{z}}^\ad\hp{2k}^\inv\pda{k}{z}&
&\text{and}&
 \ek*{\rpsdrka{2n}{z}}^\inv&=\abs{z-\ko z}\sum_{k=0}^n\pba{k}{z}\hp{2k}^\inv\ek*{\pba{k}{z}}^\ad.
\end{align*}
\eprop
\bproof
 First observe that \rrem{Schr2.7} yields \(\seqs{2n}\in\Hggeq{2n}\).
 \rprop{L1237-A} then shows that the \tnatext{} \(\seqsh{2n+1}\) of \(\seqs{2n}\) belongs to \(\Hggeq{2n+1}\).
 Let \(\shHp{2n+1}\) be the \tsHpo{\(\seqsh{2n+1}\)}.
 Then \rrem{R26-SK1} yields \(\hhp{j}=\hp{j}\) for all \(j\in\mn{0}{2n}\).
 From \rprop{P1644} we see \(\hp{2k}\in\Cgq\) for all \(k\in\mn{0}{n}\).
 Hence, we can infer \(\hhp{2k}\in\Cgq\) for all \(k\in\mn{0}{n}\).
 Furthermore, \rrem{R1019} provides \(\hhp{j}^\ad=\hhp{j}\) for all \(j\in\mn{0}{2n+1}\).
 Using \rprop{P1644}, then one can check that \(\seqsh{2n+1}\) belongs to \(\Hgeq{2n+1}\).
 Thus, we can apply \rprop{L0627} to the sequence \(\seqsh{2n+1}\) to obtain \(\det\lpsdrkha{2n+1}{z}\neq0\) and \(\ek{\lpsdrkha{2n+1}{z}}^\inv=\abs{z-\ko z}\sum_{k=0}^{n}\ek{\pdha{k}{z}}^\ad\hhp{2k}^\inv\pdha{k}{z}\), where the matrix-valued function \(\lpsdrkh{2n+1}\) is built according to \rnotap{N2150M}{N2150M.b} from the sequence \(\seqsh{2n+1}\) and \(\habcd{n+1}\) is the \sabcdo{\(\seqsh{2n+1}\)}.
 According to \rdefn{D.nat-ext}, we have \(\shu{j}=\su{j}\) for all \(j\in\mn{0}{2n}\).
 Hence, \rrem{BK8.1} yields \(\pdh{k}=\pd{k}\) for all \(k\in\mn{0}{n}\) and \rrem{LRC-tru} shows \(\lpsdrkha{2n}{z}=\lpsdrka{2n}{z}\).
 Furthermore, \rlem{L0844} provides \(\lpsdrkha{2n}{z}=\lpsdrkha{2n+1}{z}\).
 Taking additionally into account that \(\hhp{j}=\hp{j}\) holds true for all \(j\in\mn{0}{2n}\), we can conclude then that \(\lpsdrka{2n}{z}\) is invertible and fulfills \(\ek{\lpsdrka{2n}{z}}^\inv=\abs{z-\ko z}\sum_{k=0}^{n}\ek{\pda{k}{z}}^\ad\hp{2k}^\inv\pda{k}{z}\).
 The remaining assertions can be obtained similarly.
\eproof

\bleml{L0723}
 Let \(\kappa\in\Ninf\), let \(\seqska\in\Hgeqka\) with \tsHp{} \(\sHp{\kappa}\) and \sabcd{} \(\abcd{\ev{\kappa}}\), let \(n\in\NO\) be such that \(2n+1\leq\kappa\), and let \(z\in\CR\).
 Then
\begin{align*}
 \mpka{2n+1}{z}&=-\rk*{\rk{z-\ko z}\sum_{k=0}^n\ek*{\pda{k}{z}}^\ad\hp{2k}^\inv\pda{k}{z}}^\inv\rk*{\Iq+\rk{z-\ko z}\sum_{k=0}^n\ek*{\pda{k}{z}}^\ad\hp{2k}^\inv\pca{k}{z}}
\intertext{and}
 \mpka{2n+1}{z}&=-\rk*{\Iq+\rk{z-\ko z}\sum_{k=0}^n\paa{k}{z}\hp{2k}^\inv\ek*{\pba{k}{z}}^\ad}\rk*{\rk{z-\ko z}\sum_{k=0}^n\pba{k}{z}\hp{2k}^\inv\ek*{\pba{k}{z}}^\ad}^\inv.
\end{align*}
\elem
\bproof
 First observe that \rrem{Schr2.7} yields \(\seqska\in\Hggeqka\).
 Therefore, we can apply \rlem{L0819} to infer \(\det\PDa{n}{z}\neq0\) and \(\mpka{2n+1}{z}=-\ek{\PDa{n}{z}}^\inv\PCa{n}{z}\).
 Clearly, we have \(\seqs{2n}\in\Hgq{2n}\).
 From \rprop{P1644} we see then \(\hp{2k}\in\Cgq\), and, in particular, \(\det\hp{2k}\neq0\) for all \(k\in\mn{0}{n}\).
 Furthermore, \rlem{BK8.6} yields \(\det\pda{n}{z}\neq0\).
 Taking additionally into account \rlem{L0609}, we can thus infer
\[\begin{split}
 \mpka{2n+1}{z}
 &=-\ek*{\PDa{n}{z}}^\inv\PCa{n}{z}
 =-\rk*{-\ek*{\pda{n}{z}}^\ad\hp{2n}^\inv\PDa{n}{z}}^\inv\rk*{-\ek*{\pda{n}{z}}^\ad\hp{2n}^\inv\PCa{n}{z}}\\
 &=-\rk*{\rk{z-\ko z}\sum_{k=0}^n\ek*{\pda{k}{z}}^\ad\hp{2k}^\inv\pda{k}{z}}^\inv\rk*{\Iq+\rk{z-\ko z}\sum_{k=0}^n\ek*{\pda{k}{z}}^\ad\hp{2k}^\inv\pca{k}{z}}.
\end{split}\]
 The remaining identity can be obtained similarly.
\eproof

\bpropl{P1157}
 Let \(n\in\NO\), let \(\seqs{2n}\in\Hgq{2n}\) with \tsHp{} \(\sHp{2n}\) and \sabcd{} \(\abcd{n}\), and let \(z\in\CR\).
 Then
\begin{align*}
 \mpka{2n}{z}&=-\rk*{\rk{z-\ko z}\sum_{k=0}^n\ek*{\pda{k}{z}}^\ad\hp{2k}^\inv\pda{k}{z}}^\inv\rk*{\Iq+\rk{z-\ko z}\sum_{k=0}^n\ek*{\pda{k}{z}}^\ad\hp{2k}^\inv\pca{k}{z}}
\intertext{and}
 \mpka{2n}{z}&=-\rk*{\Iq+\rk{z-\ko z}\sum_{k=0}^n\paa{k}{z}\hp{2k}^\inv\ek*{\pba{k}{z}}^\ad}\rk*{\rk{z-\ko z}\sum_{k=0}^n\pba{k}{z}\hp{2k}^\inv\ek*{\pba{k}{z}}^\ad}^\inv.
\end{align*}
\eprop
\bproof
 First observe that \rrem{Schr2.7} yields \(\seqs{2n}\in\Hggeq{2n}\).
 \rprop{L1237-A} then shows that the \tnatext{} \(\seqsh{2n+1}\) of \(\seqs{2n}\) belongs to \(\Hggeq{2n+1}\).
 Let \(\shHp{2n+1}\) be the \tsHpo{\(\seqsh{2n+1}\)}.
 Then \rrem{R26-SK1} yields \(\hhp{j}=\hp{j}\) for all \(j\in\mn{0}{2n}\).
 As in the proof of \rprop{P1229} we can infer \(\seqsh{2n+1}\in\Hgeq{2n+1}\).
 Thus, we can apply \rlem{L0723} to the sequence \(\seqsh{2n+1}\) to obtain \(\mphka{2n+1}{z}=-\rk{\rk{z-\ko z}\sum_{k=0}^n\ek{\pdha{k}{z}}^\ad\hhp{2k}^\inv\pdha{k}{z}}^\inv\rk{\Iq+{\rk{z-\ko z}}\sum_{k=0}^n\ek{\pdha{k}{z}}^\ad\hhp{2k}^\inv\pcha{k}{z}}\), where the matrix-valued function \(\mphk{2n+1}\) is built according to \rnotap{CB15}{CB15.b} from the sequence \(\seqsh{2n+1}\) and \(\habcd{n+1}\) is the \sabcdo{\(\seqsh{2n+1}\)}.
 According to \rdefn{D.nat-ext}, we have \(\shu{j}=\su{j}\) for all \(j\in\mn{0}{2n}\).
 Hence, \rrem{BK8.1} yields \(\pch{k}=\pc{k}\) and \(\pdh{k}=\pd{k}\) for all \(k\in\mn{0}{n}\) and \rrem{ABM-tru} shows \(\mphka{2n}{z}=\mpka{2n}{z}\).
 Furthermore, \rlem{L3N10} provides \(\mphka{2n}{z}=\mphka{2n+1}{z}\).
 Taking additionally into account that \(\hhp{j}=\hp{j}\) holds true for all \(j\in\mn{0}{2n}\), we can conclude then \(\mpka{2n}{z}=-\rk{\rk{z-\ko z}\sum_{k=0}^n\ek{\pda{k}{z}}^\ad\hp{2k}^\inv\pda{k}{z}}^\inv\rk{\Iq+\rk{z-\ko z}\sum_{k=0}^n\ek{\pda{k}{z}}^\ad\hp{2k}^\inv\pca{k}{z}}\).
 The remaining identity can be obtained similarly.
\eproof
 
 Now we are able to see that the parameters of the matrix balls stated in \rthmss{T1546}{F49} indeed coincide:

\bpropl{P1712}
 Let \(n\in\NO\) and let \(\seqs{2n}\in\Hgq{2n}\).
 For all \(w\in\pip\), then
\begin{align*}
 \mpka{2n}{w}&=\rCa{n}{w},&
 \lpsdrka{2n}{w}&=\rga{n}{w},&
&\text{and}&
 \rpsdrka{2n}{w}&=\rda{n}{w}.
\end{align*}
\eprop
\bproof
 For all \(w\in\pip\) we have \(\abs{w-\ko w}=\abs{2\iu\rim w}=2\rim w=\iu\rk{\ko w-w}\).
 Combining \rpropss{P1157}{P1229} with \rcor{C1645} and \rlem{R1704} completes the proof.
\eproof

\appendix
\section{Some particular facts on matrix theory}\label{A-A}

\breml{R1241}
 Let \(M\in\Cpq\).
 If \(R\in\Cqq\) is invertible, then \(\ran{MR}=\ran{M}\).
 If \(L\in\Cpp\) is invertible, then \(\nul{LM}=\nul{M}\).
\erem

 We will write \(\ipE{.}{.}\) for the (left) Euclidean inner product in \(\Cq\), \tie{}, for all \(x,y\in\Cq\), let \(\ipE{x}{y}\defeq y^\ad x\).
 If \(\mM\) is a \tne{} subset of \(\Cq\), then the set \(\mM^\oc\) of all \(x\in\Cq\) which fulfill \(\ipE{x}{y}=0\) for all \(y\in\mM\) is a linear subspace of \(\Cq\) and is called (left) orthogonal complement of \(\mM\).
 
\breml{PM} 
 If \(\mU\) is a linear subspace of \(\Cq\), then there exists exactly one matrix \(\OPu{\mU}\in\Cqq\) such that both \(\OPu{\mU} x\in\mU\) and \(x-\OPu{\mU} x\in\mU^\oc\) are fulfilled for each \(x\in\Cq\).
 This matrix \(\OPu{\mU}\) is called the orthogonal projection matrix onto \(\mU\).
 In particular, \(\OPu{\mU} u=u\) for all \(u\in\mU\).
 A complex \tqqa{matrix} \(P\) is said to be an orthogonal projection matrix, if there exists a linear subspace \(\mU\) of \(\Cq\) such that \(P=\OPu{\mU}\).
\erem

\breml{R.P}
 Let \(\mathcal{U}\) be a linear subspace of \(\Cq\).
 Then \(\OPu{\mathcal{U}}\) is the unique complex \tqqa{matrix} satisfying \(\OPu{\mathcal{U}}^2=\OPu{\mathcal{U}}\) and \(\OPu{\mathcal{U}}^\ad=\OPu{\mathcal{U}}\) as well as \(\ran{\OPu{\mathcal{U}}}=\mathcal{U}\).
 Furthermore, \(\nul{\OPu{\mathcal{U}}}=\mathcal{U}^\oc\).
\erem

\breml{A.R.0<P<1}
 Let \(\mathcal{U}\) be a linear subspace of \(\Cq\).
 Then \(\Oqq\lleq\OPu{\mathcal{U}}\lleq\Iq\) and \(\OPu{\mathcal{U}}+\OPu{\mathcal{U}^\oc}=\Iq\).
\erem

\breml{tsa2}
 If \(A\in\Cpq\), then \(\ran{A^\ad}=\nul{A}^\oc\) and \(\nul{A^\ad}=\ran{A}^\oc\).
\erem

\breml{tsa4}%
 Let \(A\in\Cggq\).
 For every \(B\in\Cqp\), then \(B^\ad AB\in\Cggp\).
\erem

\breml{A.R.r-sqrt}
 If \(A\in\Cggq\), then \(\ran{\sqrt{A}}=\ran{A}\) and \(\nul{\sqrt{A}}=\nul{A}\).
\erem

\breml{ZR2}
 The mapping \(w\colon\Cggq\to\Cggq\) defined by \(w(A)\defeq\sqrt{A}\) is continuous.
\erem

\breml{tsa15}
 Let \(C\in\Cpq\).
 Then it is readily checked that \(C\) is contractive if and only if the matrix \(\Iq-C^\ad C\) is \tnnH{}.
\erem

\breml{KTK}
 Let \(\seq{C_n}{n}{1}{\infi}\) be a sequence of contractive complex \tpqa{matrices}.
 Due to the established theorem of Bolzano--Weierstra\ss{}, there are a contractive complex \tpqa{matrix} \(C\) and a subsequence \(\seq{C_{n_k}}{j}{1}{\infi}\) of \(\seq{C_n}{n}{1}{\infi}\) such that \(\lim_{k\to\infty}C_{n_k}=C\) (see also, \teg{}, \zitaa{MR2380292}{\csatz{11.7}{104}}).
\erem

\breml{A.R.A++*}
 If \(A\in\Cpq\), then \(\rk{A^\mpi}^\mpi=A\) and \(\rk{A^\ad}^\mpi=\rk{A^\mpi}^\ad\).
\erem

\breml{A.R.r-mpi}
 If \(A\in\Cpq\), the \(\ran{A^\mpi}=\ran{A^\ad}\) and \(\nul{A^\mpi}=\nul{A^\ad}\).
\erem

\breml{tsb3}
 If \(A\in\Cpq\), then \(AA^\mpi=\OPu{\ran{A}}\) and \(A^\mpi A=\OPu{\ran{A^\ad}}\).
\erem

\breml{tsa12}
 Let \(A\in\Cpq\).
 Using \(AA^\mpi A=A\), one can easily see that the following statements hold true:
\benui
 \il{tsa12.a} Let \(B\in\Coo{p}{r}\).
 Then \(\ran{B}\subseteq\ran{A}\) if and only if \(AA^\mpi B=B\).
 \il{tsa12.b} Let \(C\in\Coo{r}{q}\).
 Then \(\nul{A}\subseteq\nul{C}\) if and only if \(CA^\mpi A=C\).
\eenui
\erem

\bpropnl{\tcf{}\ \zitaa{MR1105324}{\cthm{10.4.1}{218}}}{MPK}
 Suppose that \((A_n)_{n=1}^\infi\) is a sequence of complex \tpqa{matrices} which converges to a complex \tpqa{matrix} \(A\).
 Then \((A_n^\mpi )_{n=1}^\infi\) is convergent if and only if there is a positive integer \(m\) such that \(\rank A_n=\rank A\) for each integer \(n\) fulfilling \(n\geq m\).
 In this case, \((A_n^\mpi )_{n=1}^\infi\) converges to \(A^\mpi \).
\eprop

\bleml{LemC1}
 Let \(\mG\) be a non-empty subset of \(\C\).
 Furthermore, let \(F\colon\mG\to\Cpq\) be a continuous function such that \(\rank F\rk{z} = \rank F\rk{w}\) for every choice of \(z\) and \(w\) in \(\mG\).
 Then \(G\colon\mG\to\Cqp\) given by \(G\rk{z}\defeq\ek{F\rk{z}}^\mpi\) is continuous.
\elem
\bproof
 Let us consider an arbitrary \(z\in\mG\).
 Let \(\seq{z_n}{n}{1}{\infi}\) be a sequence of points belonging to \(\mG\) such that \(\lim_{n\to\infi}z_n=z\).
 Since \(F\) is continuous, then we perceive \(\lim_{n\to\infi}F\rk{z_n}=F\rk{z}\).
 Thus, in view of \(\rank F\rk{z_n}=\rank F\rk{z}\) for all \(n\in\N\), from \rprop{MPK} we conclude \(\lim_{n\to\infi} G\rk{z_n} =G\rk{z}\).
 Consequently, \(G\) is continuous at \(z\).
\eproof

\breml{A.R.A+>}
 If \(A\in\Cggq\), then \(A^\mpi\in\Cggq\) and \(\sqrt{A^\mpi}=\sqrt{A}^\mpi\).
\erem

\breml{A.R.A+sqrt}
 If \(A\in\Cggq\), then \(\sqrt{A}^\mpi A=\sqrt{A}=A\sqrt{A}^\mpi\).
\erem

\bleml{KB-5}
 Let both \(L\in\Cpp\) and \(R\in\Cqq\) be invertible matrices.
 Let \(M\in\Cpq\) and let \(N\defeq LMR^\inv\).
 Then \(X\defeq MR^\inv\) admits the representation \(X=L^\inv N\) and the Moore--Penrose inverse of \(X\) admits the representations \(X^\mpi =N^\mpi NRM^\mpi \) and \(X^\mpi =N^\mpi LMM^\mpi \).
\elem
\bproof
 Clearly, \(X=L^\inv N\).
 From \rrem{R1241} we can infer then \(\ran{X}=\ran{MR^\inv}=\ran{M}\) and \(\nul{X}=\nul{L^\inv N}=\nul{N}\).
 Using \rrem{tsa2}, in particular \(\ran{X^\ad}=\ran{N^\ad}\) follows.
 Consequently, \rrem{tsb3} yields \(XX^\mpi =MM^\mpi \) and \(X^\mpi X=N^\mpi N\).
 Thus, we receive \(
 X^\mpi
 =X^\mpi XX^\mpi
 =X^\mpi MM^\mpi
 =X^\mpi MR^\inv RM^\mpi
 =X^\mpi XRM^\mpi
 =N^\mpi NRM^\mpi
\).
 Taking into account \(NR=LM\), then as well \(X^\mpi=N^\mpi LMM^\mpi\) follows.
\eproof

\bleml{B5}
 Let \(M\in\Cqq\), let \(A\defeq\rre M\), and let \(B\defeq\rim M\).
 If \(\ran{A}\subseteq\ran{B}\), then 
\(%
 MB^\mpi M^\ad
 =M^\ad B^\mpi M
\).
\elem
\bproof
 Suppose \(\ran{A}\subseteq\ran{B}\).
 In view of \(A^\ad=A\) and \(B^\ad=B\), using \rrem{tsa2}, we obtain \(\nul{B}\subseteq\nul{A}\).
 \rrem{tsa12} then yields \(BB^\mpi A=A\) and \(AB^\mpi B=A\).
 Consequently, in view of \(M=A+\iu B\) and \(M^\ad=A-\iu B\), then we receive 
\[\begin{split}
 &MB^\mpi M^\ad-M^\ad B^\mpi M
 =\rk{A+\iu B}B^\mpi \rk{A-\iu B}-\rk{A-\iu B}B^\mpi \rk{A+\iu B}\\
 &=\rk{AB^\mpi A-\iu AB^\mpi B+\iu BB^\mpi A+BB^\mpi B}-\rk{AB^\mpi A+\iu AB^\mpi B-\iu BB^\mpi A+BB^\mpi B}
 =\Oqq,
\end{split}\]
 which completes the proof.
\eproof

\bleml{C2}
 Let \(A\in\Cpq\) and let \(B\in\Coo{q}{r}\).
\benui
 \il{C2.a} If \(\nul{B^\ad}\subseteq\nul{A}\), then \(\ran{AB}=\ran{A}\).  
 \il{C2.b} If \(\ran{B}\subseteq\ran{A^\ad}\), then \(\nul{AB}=\nul{B}\).
\eenui
\elem
\bproof
 The results in this lemma are considered established.
 For the reader's convenience, we include a proof, though, since we did not succeed in finding an appropriate one.

 \eqref{C2.a} Suppose \(\nul{B^\ad}\subseteq\nul{A}\).
 According to both \rremp{tsa12}{tsa12.b} and \rrem{A.R.A++*}, then \(A= A\rk{B^\ad}^\mpi B^\ad=A\rk{B^\mpi}^\ad B^\ad=A\rk{BB^\mpi}^\ad=ABB^\mpi \).
 Consequently, \(\ran{AB}=\ran{A}\) follows, by virtue of \(\ran{AB}\subseteq\ran{A}=\ran{ABB^\mpi}\subseteq\ran{AB}\).
 
 \eqref{C2.b} Assume \(\ran{B}\subseteq\ran{A^\ad}\).
 Using both \rremp{tsa12}{tsa12.a} and \rrem{A.R.A++*}, we obtain \(B=A^\ad\rk{A^\ad}^\mpi B=A^\ad\rk{A^\mpi}^\ad B=\rk{A^\mpi A}^\ad B=A^\mpi AB\),  implying \(\nul{AB}=\nul{B}\), since \(\nul{AB}\subseteq\nul{A^\mpi AB}=\nul{B}\subseteq\nul{AB}\) holds true.
\eproof

\section{On some connections between contractive matrices and matrices with \tnnH{} imaginary part}\label{S1624}

 This appendix contains several algebraic results, which are used for proofs in \rsec{Cha13}.
 
\breml{R1043}
 The matrix \(\Jdefpq\) given by 
\beql{Jdef}
 \Jdefpq
 \defeq
\Mat{
 \Ip&\Opq\\ 
 \Oqp&-\Iq
} 
\eeq
 is a \taaa{\rk{p+q}}{\rk{p+q}}{signature} matrix, \tie{}, \(\Jdefpq^\ad=\Jdefpq\) and \(\Jdefpq^2=\Iu{p+q}\) hold true.
 Moreover, \(\tmat{A \\ B}^\ad\rk{-\Jdefpq}\tmat{A \\ B}=B^\ad B-A^\ad A\) for all \(A\in\Cpq\) and \(B\in\Cqq\).
 In particular, the case \(B=\Iq\) is of interest.
\erem

 As usual, let \(\CEPq\defeq\setaca{M\in\Cqq}{\ran{M}=\ran{M^\ad}}\).
 
\breml{P1316}
 Let \(M\in\Cqq\).
 From \rrem{tsa2} one can easily see then that \(M\in\CEPq\) if and only if \(\nul{M}=\nul{M^\ad}\).
\erem

\bnotal{N.psdreim}
 Let \(\psdreq\defeq\setaca{M\in\Cqq}{\rre M\in\Cggq}\) and \(\psdimq\defeq\setaca{M\in\Cqq}{\rim M\in\Cggq}\).
\enota

\blemnl{\zitaa{MR3014198}{\clem{A.8(c)}{112} and \clem{A.10(c)}{113}}}{L1321}
 If \(M\in\psdreq\cup\psdimq\) then \(M\in\CEPq\).
\elem

\bleml{L1058}
 Let \(E\in\Cqq\), let \(b\in\R\setminus\set{0}\), let \(B\defeq b^\inv\rim E\), and let \(P\defeq\OPu{\ran{E}}\).
 Suppose \(B\in\Cggq\) and \(\ran{B}=\ran{E}\).
 Then the matrices 
\begin{align}\label{L1058.0}
 H &\defeq
 \Mat{
  \sqrt{B}^\mpi  & -\sqrt{B}^\mpi E\\
  \sqrt{B}^\mpi  & -\sqrt{B}^\mpi E^\ad}&
&\text{and}&
 K &\defeq
 \Mat{-E^\ad\sqrt{B}^\mpi  & E\sqrt{B}^\mpi\\
 -\sqrt{B}^\mpi  & \sqrt{B}^\mpi}
\end{align}
 fulfill the equations \(H K =2\iu b\diag\rk{P,P}\) and \(K H =2\iu b\diag\rk{P,P}\) 
 as well as
\begin{align}
 H^\ad\rk{-\Jdefqq} H &
 =2b\ek*{\diag\rk{P,P}}^\ad\rk{-\Jimq}\diag\rk{P,P}\label{L1058.A}
\intertext{and}
 K^\ad\rk{-\Jimq}K 
 &=2b\ek*{\diag\rk{P,P}}^\ad\rk{-\Jdefqq}\diag\rk{P,P}.\label{L1058.B}
\end{align}
\elem
\bproof
 Let \(W\defeq\sqrt{B}\) and let \(R\defeq\sqrt{B}^\mpi\).
 From \rrem{A.R.A+>} we can infer
\begin{align}\label{L1058.2} 
 R^\ad&=R&
&\text{and}&
 R^2&=B^\mpi.
\end{align}
 Clearly, \(W^\ad=W\) and \(W^\mpi=R\).
 From \rrem{A.R.A+sqrt} we see \(BR=W\).
 \rrem{tsb3} yields \(W^\mpi W=\OPu{\ran{W^\ad}}\). %
 \rrem{A.R.r-sqrt} shows \(\ran{W}=\ran{B}\).
 Taking additionally into account the assumption \(\ran{B}=\ran{E}\), we obtain
\beql{L1058.9}
 RBR
 =W^\mpi W
 =\OPu{\ran{W^\ad}}
 =\OPu{\ran{W}}
 =\OPu{\ran{B}}
 =\OPu{\ran{E}}
 =P.
\eeq
 In view of \(b\in\R\setminus\set{0}\), we have \(B=\rim\rk{b^\inv E}\).
 Regarding the assumption \(B\in\Cggq\) and \rnota{N.psdreim}, we can thus apply \rlem{L1321} to obtain \(b^\inv E\in\CEPq\), \tie{}, \(\ran{b^\inv E}=\ran{\rk{b^\inv E}^\ad}\), which implies \(\ran{E}=\ran{E^\ad}\).
 Using additionally \(\ran{B}=\ran{E}\) and \(b\neq0\), we conclude
\beql{N33}
 \ran{\rre E}
 =\ran{E+E^\ad}
 \subseteq\ran{E}+\ran{E^\ad}
 =\ran{E}+\ran{E}
 =\ran{E}
 =\ran{B}
 =\ran{bB}.
\eeq
 Clearly, \(bB=\rim E\).
 Regarding \eqref{N33}, consequently, \rlem{B5} shows that \(E\rk{bB}^\mpi E^\ad=E^\ad\rk{bB}^\mpi E\).
 Because of \(b\neq0\) and \eqref{L1058.2}, then
\beql{L1058.10}
 ER^2 E^\ad
 =E^\ad R^2E
\eeq
 follows.
 In view of \(\ran{B}=\ran{E}\) and \(B^\ad=B\), from \rrem{tsb3} we obtain
\begin{align}\label{L1058.8} 
 P&=BB^\mpi&
&\text{and}&
 P&=B^\mpi B.
\end{align}
 By virtue of \rrem{R.P}, we furthermore can infer \(P^\ad P=P\).
 In view of \(bB=\rim E\), we have \(E-E^\ad=2\iu bB\).
 Because of \eqref{L1058.0}, \eqref{L1058.2}, \eqref{L1058.10}, and \eqref{L1058.8}, we conclude
\[\begin{split}
 K H 
 &=
 \Mat{
 -E^\ad R   & E R\\
 - R   &  R}
 \Mat{
   R   & - R  E\\
   R   & - R  E^\ad}
 =
 \begin{pmat}[{|}]
 -E^\ad  R^2+E R^2 & E^\ad  R^2E-E R^2E^\ad\cr\- 
 - R^2+ R^2& R^2E- R^2E^\ad\cr
 \end{pmat}\\
 &=\diag\rk*{\rk{E-E^\ad}B^\mpi,B^\mpi\rk{E-E^\ad}}
 =\diag\rk{2\iu bBB^\mpi,2\iu bB^\mpi B}
 =2\iu b\diag\rk{P,P}. 
\end{split}\]
 Using \eqref{L1058.0} and \eqref{L1058.9}, we receive
\[%
 H K 
 =
 \diag\rk*{R \rk{E-E^\ad} R,R \rk{E-E^\ad} R}
 =
 \diag\rk{2\iu bR  B R,2\iu bR  B R}
 =2\iu b\diag\rk{ P,P}.
\]%
 Once more, taking into account \eqref{L1058.0}, \eqref{Jdef}, \eqref{L1058.2}, \eqref{L1058.10}, and \eqref{L1058.8}, we perceive
\[\begin{split}
 H^\ad\rk{-\Jdefqq} H 
 &=
 \Mat{
   R^\ad   & R^\ad\\
   \rk{- R  E}^\ad   & \rk{- R  E^\ad}^\ad}
 \Mat{
  -\Iq & \Oqq \\ 
  \Oqq & \Iq
 }
 \Mat{
   R   & - R  E\\
   R   & - R  E^\ad}\\
 &=
 \begin{pmat}[{|}]
 -R^2+R^2 & R^2E-R^2E^\ad\cr\-
 E^\ad R^2-ER^2& -E^\ad R^2E+ER^2E^\ad\cr
 \end{pmat}\\
 &=\begin{pmat}[{|}]
 \Oqq & B^\mpi\rk{E-E^\ad}\cr\-
 -\rk{E-E^\ad}B^\mpi  &\Oqq\cr
 \end{pmat}
 =
 \Mat{
 \Oqq & 2\iu bP\\ 
 -2\iu bP&\Oqq
 },
\end{split}\]
 which, in view of \(P^\ad P=P\) and \eqref{JQ}, implies \eqref{L1058.A}.
 According to \eqref{L1058.0}, \eqref{JQ}, \eqref{L1058.2}, and \eqref{L1058.9}, analogous straightforward calculations yield \(K^\ad\rk{-\Jimq}K=\diag\rk{-2bP,2bP}\), 
 which, by virtue of \(P^\ad P=P\) and \eqref{Jdef}, implies \eqref{L1058.B}.
\eproof

\bleml{L1104}
 Let \(E\in\Cqq\), let \(b\in\R\setminus\set{0}\), let \(B\defeq b^\inv\rim E\), and let \(\cpP,\cpQ\in\Cqq\) be such that \(B\in\Cggq\) and \(\ran{\cpP}\subseteq\ran{E}=\ran{B}=\ran{\cpP-E^\ad\cpQ}\) 
 hold true.
 Then
\beql{B4-7}
 C
 \defeq\sqrt{B}^\mpi\rk{\cpP -E\cpQ}\rk{\cpP -E^\ad \cpQ}^\mpi\sqrt{B}
\eeq
 fulfills 
\beql{L1104.A}
 \OPu{\ran{E}}-C^\ad C
 =4b\ek*{\rk{\cpP -E^\ad \cpQ}^\mpi\sqrt{B}}^\ad\rim\rk{\cpQ^\ad \cpP }\ek*{\rk{\cpP -E^\ad \cpQ}^\mpi\sqrt{B}}.
\eeq
 Furthermore, if \(b\rim\rk{\cpQ^\ad \cpP }\) is \tnnH{}, then \(C\) is contractive.
\elem
\bproof
 Let \(R\defeq\sqrt{B}^\mpi\).
 Clearly, the matrix \(H\defeq\smat{R&-RE\\R&-RE^\ad}\) fulfills \(H\tmat{\cpP \\\cpQ}=\tmat{R X \\R Y }\), where \(X\defeq \cpP -E\cpQ\) and \(Y\defeq \cpP -E^\ad \cpQ\).
 In view of the assumptions \(B\in\Cggq\) and \(\ran{B}=\ran{E}\), we can apply \rlem{L1058} to obtain \eqref{L1058.A}, where \(P\defeq\OPu{\ran{E}}\).
 From \rrem{R.P} we can conclude \(P^\ad P=P\).
 Regarding the assumption \(\ran{\cpP }\subseteq\ran{E}\), we get \(P\cpP =\cpP \).
 Hence, we have \(\rk{P\cpQ}^\ad P\cpP =\cpQ^\ad P\cpP =\cpQ^\ad \cpP \).
 Taking additionally into account \rremss{R1043}{AB53N} as well as \eqref{L1058.A}, we conclude
\beql{L1104.1}\begin{split}%
 &\rk{ R  Y }^\ad\rk{ R  Y }-\rk{ R  X }^\ad \rk{ R  X }
 =
 \Mat{
  R  X  \\ 
  R  Y  
 }^\ad\rk{-\Jdefqq}
 \Mat{
  R  X  \\ 
  R  Y  
 }
 =\rk*{H\Mat{\cpP \\\cpQ}}^\ad\rk{-\Jdefqq}\rk*{H\Mat{\cpP \\\cpQ}}\\
 &=\Mat{\cpP  \\ \cpQ}^\ad\ek*{H^\ad\rk{-\Jdefqq} H}\Mat{\cpP  \\ \cpQ}
 =2b\Mat{\cpP  \\ \cpQ}^\ad\ek{\diag\rk{P,P}}^\ad\rk{-\Jimq}\diag\rk{P,P}\Mat{\cpP  \\ \cpQ}\\
 &=2b\Mat{P\cpP  \\P\cpQ}^\ad\rk{-\Jimq}\Mat{P\cpP  \\P\cpQ}
 =4b\rim\rk*{\rk{P\cpQ}^\ad P\cpP }
 =4b\rim\rk{\cpQ^\ad \cpP }.
\end{split}\eeq
 Let \(W\defeq\sqrt{B}\).
 Clearly, \(W^\ad=W\) and \(W^\mpi=R\).
 \rrem{A.R.r-sqrt} shows \(\ran{W}=\ran{B}\).
 In view of the assumption \(\ran{B}=\ran{\cpP -E^\ad \cpQ}\), from \rrem{tsb3} we then infer \(WW^\mpi=BB^\mpi=YY^\mpi\).
 \rrem{tsb3} yields furthermore \(W^\mpi W=\OPu{\ran{W^\ad}}\).
 Taking additionally into account \(\ran{B}=\ran{E}\), we consequently get
\beql{L1104.2}
 R Y  Y^\mpi W
 =W^\mpi WW^\mpi W
 =W^\mpi W
 =\OPu{\ran{W^\ad}}
 =\OPu{\ran{W}}
 =\OPu{\ran{B}}
 =\OPu{\ran{E}}
 =P.
\eeq
 According to \eqref{B4-7}, we have \(C=RXY^\mpi W\).
 Thus, in view of \eqref{L1104.2} and \eqref{L1104.1}, we obtain 
\[\begin{split}
 P-C^\ad C
 &=P^\ad P-C^\ad C
 =\rk{ R  Y  Y^\mpi W }^\ad \rk{ R  Y  Y^\mpi W }- \rk{ R  X  Y^\mpi W }^\ad \rk{ R  X  Y^\mpi W }\\
 &=\rk{ Y^\mpi W }^\ad\ek*{\rk{ R  Y }^\ad \rk{ R  Y } -\rk{ R  X }^\ad \rk{ R  X }}\rk{ Y^\mpi W }
 =4b\rk{ Y^\mpi W }^\ad\ek*{\rim\rk{\cpQ^\ad \cpP }}\rk{ Y^\mpi W },
\end{split}\]
 which implies \eqref{L1104.A}.
 Now suppose \(b\rim\rk{\cpQ^\ad \cpP }\in\Cggq\).
 Using \rrem{tsa4}, from \eqref{L1104.A} then \(P-C^\ad C\in\Cggq\) follows.
 Since \rrem{A.R.0<P<1} yields \(\Iq-P\in\Cggq\), we can conclude \(\Iq-C^\ad C\in\Cggq\).
 Because of \rrem{tsa15}, thus \(C\) is contractive.
\eproof

\bleml{L2043}
 Let \(E\in\Cqq\), let \(b\in\R\setminus\set{0}\), let \(B\defeq b^\inv\rim E\), and let \(P\defeq\OPu{\ran{E}}\) and \(Q\defeq\OPu{\nul{E}}\).
 Suppose \(B\in\Cggq\) and \(\ran{B}=\ran{E}\).
 Let \(C\in\Cqq\) and let
\begin{align*}%
 \cpP &\defeq E\sqrt{B}^\mpi-E^\ad\sqrt{B}^\mpi CP,&
&\text{and}&
 \cpQ&\defeq\sqrt{B}^\mpi-\sqrt{B}^\mpi CP+Q. 
\end{align*}
 Then \(\rank\smat{\cpP \\\cpQ}=q\),
 \(P\cpP =\cpP \),
 \(\cpP P=\cpP \),
 \(\cpQ P=\cpQ-Q\), as well as
\begin{align}\label{L2043.A}
 b^\inv\rim\rk{\cpQ ^\ad \cpP }&=P^\ad\ek*{\Iq-\rk{PC}^\ad\rk{PC}}P&
&\text{and}&
 \sqrt{B}^\mpi\rk{\cpP -E\cpQ }\rk{\cpP -E^\ad \cpQ }^\mpi\sqrt{B}
 &=PCP.
\end{align}
 Furthermore, if \(C\) is contractive, then \(b^\inv\rim\rk{\cpQ^\ad\cpP}\) is \tnnH{}.
\elem
\bproof
 In view of \(b\in\R\setminus\set{0}\), we have \(B=\rim\rk{b^\inv E}\).
 Regarding \(B\in\Cggq\) and \rnota{N.psdreim}, we can thus apply \rlem{L1321} to obtain \(b^\inv E\in\CEPq\).
 \rrem{P1316} then shows \(\nul{b^\inv E}=\nul{\rk{b^\inv E}^\ad}\), implying \(\nul{E}=\nul{E^\ad}\).
 According to \rrem{tsa2}, we have \(\nul{E^\ad}=\ran{E}^\oc\).
 Consequently, \(\nul{E}=\ran{E}^\oc\) follows.
 From \rrem{A.R.0<P<1} we can thus conclude \(P+Q=\Iq\).
 Using \rrem{R.P}, we get furthermore \(\ran{Q}=\nul{E}=\nul{E^\ad}\), implying
\begin{align}\label{L2043.1}
 EQ&=\Oqq&
&\text{and}&
 E^\ad Q&=\Oqq.
\end{align}
 Let \(W\defeq\sqrt{B}\) and let \(R\defeq\sqrt{B}^\mpi\).
 As in the proof of \rlem{L1058} we have then \(R^\ad=R\) and \eqref{L1058.9}.
 Clearly, \(W^\ad=W\) and \(W^\mpi=R\).
 Using \rrem{A.R.r-mpi}, we can infer then \(\ran{R}=\ran{W^\ad}=\ran{W}\).
 \rrem{A.R.r-sqrt} shows \(\ran{W}=\ran{B}\).
 Consequently, \(\ran{R}=\ran{B}\) follows.
 \rrem{tsb3} yields \(B^\mpi B=\OPu{\ran{B^\ad}}\).
 Taking additionally into account \(B^\ad=B\), we obtain
\beql{L2043.2}
 R
 =\OPu{\ran{R}}R
 =\OPu{\ran{B}}R
 =\OPu{\ran{B^\ad}}R
 =B^\mpi BR.
\eeq
 Clearly, \(bB=\rim E\).
 Hence, we have \(E-E^\ad=2\iu bB\).
 Furthermore, observe that
\begin{align}\label{L2043.0}
 \cpP &=E R -E^\ad R  CP&
&\text{and}&
 \cpQ &=R - R  CP+Q.
\end{align}
 
 Consider now an arbitrary \(v\in\nul{\tmat{\cpP \\\cpQ }}\).
 Then \(\cpP v=\Ouu{q}{1}\) as well as \(\cpQ v=\Ouu{q}{1}\), and, in view of \eqref{L2043.0}, consequently,
\begin{align}\label{L2043.6}
 E R  v&=E^\ad R  CPv&
&\text{and}&
  R  v+Qv&= R  CPv.
\end{align}   
 Thus, \eqref{L2043.6} and the second equation in \eqref{L2043.1} imply \(E R  v =E^\ad R  CPv =E^\ad \rk{ R  v+Qv} =E^\ad R  v\).
 Therefore, \(\rk{E-E^\ad} R  v=\Oqq\).
 In view of \(E-E^\ad=2\iu bB\) and \(b\neq0\), hence \(BRv=\Oqq\).
 Consequently, \eqref{L1058.9} and \eqref{L2043.2} show that
 \begin{align}\label{L2043.12}%
 Pv&=RBRv=\Oqq&
&\text{and}&
 Rv&=B^\mpi BRv=\Oqq.
\end{align}
 From the second equation in \eqref{L2043.6} and \eqref{L2043.12} we discern \(Qv= R  CPv- R  v=\Oqq\).
 Regarding \(P+Q=\Iq\) and the first equation in \eqref{L2043.12}, then \(v=(P+Q)v=\Ouu{q}{1}\) follows.
 Hence, \(\nul{\tmat{\cpP \\\cpQ }}\subseteq\set{\Ouu{q}{1}}\).
 Thus, \(\rank\tmat{\cpP \\\cpQ }=q\) is proved.
 
 In view of \(b^\inv E\in\CEPq\), we have \(\ran{b^\inv E}=\ran{\rk{b^\inv E}^\ad}\), implying \(\ran{E}=\ran{E^\ad}\).
 By virtue of \(P=\OPu{\ran{E}}\), then \(PE=E\) and \(PE^\ad=E^\ad\) follow.
 Regarding \eqref{L2043.0}, consequently, we get \(P\cpP =\cpP \). 
 Because of \(\ran{R}=\ran{B}\) and the assumption \(\ran{B}=\ran{E}\), we have 
\( 
 R
 =\OPu{\ran{R}}R
 =\OPu{\ran{B}}R
 =PR
\).
 \rrem{R.P} yields \(P^2=P\) and \(P^\ad=P\).
 Regarding \(R^\ad=R\), we get then \(RP=\rk{PR}^\ad=R\).
 Taking additionally into account \(P^2=P\) and \eqref{L2043.0}, consequently \(\cpP P=\cpP \) follows. 
 According to \rrem{tsa2}, we have \(\ran{E^\ad}=\nul{E}^\oc\).
 Using \rrem{R.P} and \(\ran{E}=\ran{E^\ad}\), we can infer then \(\ran{P}=\ran{E}=\ran{E^\ad}=\nul{E}^\oc=\nul{Q}\), implying \(QP=\Oqq\).
 Taking additionally into account \eqref{L2043.0} as well as \(RP=R\) and \(P^2=P\), we acquire
\[
 \cpQ P
 =RP-RCP^2+QP
 =R-RCP
 =\cpQ -Q.
\]
  
 Let
\begin{align}\label{L2043.7}
 H&\defeq\Mat{R&-RE\\R&-RE^\ad}&
&\text{and}&
 K&\defeq\Mat{-E^\ad R&ER\\-R&R}.
\end{align}
 In view of \eqref{L2043.0}, we have then
\beql{L2043.84}
 K\Mat{CP\\\Iq}
 =\Mat{-E^\ad RCP+ER\\-RCP+R}
 =\Mat{\cpP \\\cpQ -Q}.
\eeq
 Using \eqref{L2043.1}, we can infer \(Q^\ad E=\Oqq\) and \(Q^\ad E^\ad=\Oqq\).
 Because of \eqref{L2043.0}, consequently \(Q^\ad \cpP =\Oqq\) follows.
 Hence, we have \({\rk{\cpQ -Q}}^\ad \cpP =\cpQ ^\ad \cpP \).
 We can apply \rlem{L1058} to obtain \(K^\ad\rk{-\Jimq}K=2b\ek{\diag\rk{P,P}}^\ad\rk{-\Jdefqq}\diag\rk{P,P}\).
 Taking additionally into account \eqref{L2043.84} and \rremss{AB53N}{R1043}, we receive
\[\begin{split}
 \frac{1}{b}\rim\rk{\cpQ ^\ad \cpP }
 &=\frac{1}{b}\rim\rk*{\rk{\cpQ -Q}^\ad \cpP }
 =\frac{1}{2b}\Mat{\cpP \\\cpQ -Q}^\ad\rk{-\Jimq}\Mat{\cpP \\\cpQ -Q}
 =\frac{1}{2b}\Mat{CP\\\Iq}^\ad K^\ad\rk{-\Jimq}K\Mat{CP\\\Iq}\\
 &=\Mat{CP\\\Iq}^\ad\ek*{\diag\rk{P,P}}^\ad\rk{-\Jdefqq}\diag\rk{P,P}\Mat{CP\\\Iq}
 =\Mat{PCP\\P}^\ad\rk{-\Jdefqq}\Mat{PCP\\P}\\
 &=P^\ad P-\rk{PCP}^\ad\rk{PCP},
\end{split}\]
 which implies the first equation in \eqref{L2043.A}.
 From \eqref{L2043.7} and \eqref{L2043.1} we realize
\beql{L2043.91} 
 H\Mat{\cpP \\\cpQ -Q}
 =\Mat{R\cpP -RE\rk{\cpQ -Q}\\R\cpP -RE^\ad\rk{\cpQ -Q}}
 =\Mat{R X \\R Y },
\eeq
 where \(X\defeq \cpP -E\cpQ \) and \(Y\defeq \cpP -E^\ad \cpQ \).
 We can apply \rlem{L1058} to obtain \(HK=2\iu b\diag\rk{P,P}\).
 Using additionally \eqref{L2043.84}, we acquire 
\beql{L2043.88}
 H\Mat{\cpP  \\  \cpQ -Q }
 =HK\Mat{CP \\ \Iq}
 =2\iu b\diag\rk{P,P}\Mat{CP \\ \Iq}
 =\Mat{2\iu bPCP \\ 2\iu bP }.
\eeq
 Comparing \eqref{L2043.91} and \eqref{L2043.88}, we obtain 
\beql{L2043.92}
 RX 
 =2\iu bPCP.
\eeq
 \rrem{A.R.A+sqrt} yields \(BR=W\).
 Applying additionally \eqref{L2043.0} and the second equation in \eqref{L2043.1}, we infer
\beql{L2043.93}\begin{split}
 Y
 &=\cpP -E^\ad \cpQ 
 =\rk{E R -E^\ad R  CP}-E^\ad\rk{R-R  CP+Q}\\
 &=ER-E^\ad R -E^\ad Q
 =\rk{E-E^\ad} R 
 =2\iu bBR 
 =2\iu bW.
\end{split}\eeq
 \rrem{tsb3} yields \(W^\mpi W=\OPu{\ran{W^\ad}}\).
 Taking additionally into account \(W^\ad=W\) and \(\ran{W}=\ran{B}\) as well as \(\ran{B}=\ran{E}\), we consequently get \(W^\mpi W=\OPu{\ran{W^\ad}}=\OPu{\ran{W}}=\OPu{\ran{B}}=\OPu{\ran{E}}=P\).
 In view of \eqref{L2043.92} and \eqref{L2043.93} as well as \(b\neq0\) and \(P^2=P\), then 
\[
 RXY^\mpi W 
 =2\iu bPCP\rk{2\iu bW}^\mpi W
 =PCPW^\mpi W 
 =PCPP
 =PCP
\]
 follows, implying the second equation in \eqref{L2043.A}.
 
 Now suppose that \(C\) is contractive.
 Clearly, \(P\) is contractive.
 Thus, \(PC\) is contractive as well.
 Consequently, \rrem{tsa15} provides that \(\Iq-\rk{PC}^\ad\rk{PC}\) is \tnnH{}.
 In view of \rrem{tsa4}, then \(P^\ad \ek{\Iq-\rk{PC}^\ad \rk{PC}}P\) is \tnnH{} as well.
 From the first equation in \eqref{L2043.A}, therefore \(b^\inv\rim\rk{\cpQ ^\ad \cpP }\in\Cggq\) follows.
\eproof

\section{Some facts on the integration theory of \tnnH{} measures}\label{A-B}

 In this section, we present basic facts regarding the integration theory with respect to \tnnH{} measures.
 Throughout the section, let \(\K\in\set{\R,\C}\).
 We write \(\BsA{\K}\) denoting the \(\sigma\)\nobreakdash-algebra of all Borel subsets of \(\K\).
 Let \(\Omega\) be a \tne{} set and let \(\fA\) be a \(\sigma\)\nobreakdash-algebra on \(\Omega\).
 Consider a measure \(\nu\) on the measurable space \(\OA\).
 We use \(\LOK{\nu}\) to denote the set of all \(\fA\)\nobreakdash-\(\BsA{\K}\)\nobreakdash-measurable functions \(f\colon\Omega\to\K\) such that \(\int_{\Omega}\abs{f}\dif\nu<\infty\).
 We will write \(\BsA{\pq}\) for the \(\sigma\)\nobreakdash-algebra of all Borel subsets of \(\Cpq\).
 An \(\fA\)\nobreakdash-\(\BsA{\pq}\)\nobreakdash-measurable function \(F\colon\Omega\to\Cpq\) is said to be integrable with respect to \(\nu\) if \(F=\mat{f_{jk}}_{\substack{j=1,\dotsc,p\\k=1,\dotsc,q}}\) belongs to \(\ek{\LOC{\nu}}^\xx{p}{q}\), \tie{}, all entries \(f_{jk}\) belong to the class \(\LOC{\nu}\).
 In this case, let \(\int_\Omega F\dif\nu\defeq\mat{\int_\Omega f_{jk}\dif\nu}_{\substack{j=1,\dotsc,p\\k=1,\dotsc,q}}\).
 A matrix-valued function \(\mu\) the domain of which is \(\fA\) and the values of which belong to the set \(\Cggq\) of all \tnnH{} complex \tqqa{matrices} is called \tnnH{} \tqqa{measure} on \(\OA\) if it is \(\sigma\)\nobreakdash-additive, \tie{}, if \(\mu\) fulfills \(\mu\rk{\bigcup_{n=1}^\infi A_n} =\sum_{n=1}^{\infi}\mu\rk{A_n}\) for each sequence \(\seq{A_n}{n}{1}{\infi}\) of pairwise disjoint sets belonging to \(\fA\).
 By \(\Mggqa{\Omega,\fA}\) we denote the set of all \tnnH{} \tqqa{measures} on \(\OA\), \tie{}, the set of all \(\sigma\)\nobreakdash-additive mappings \(\mu\colon\fA \to\Cggq\).
 Let \(\mu=\matauuuo{\mu_{jk}}{j}{k}{1}{q}\in\Mggqa{\Omega,\fA}\).
 For each \(j\in\mn{1}{q}\) and for each \(k\in\mn{1}{q}\), the mapping \(\mu_{jk}\) describes a complex measure on \(\OA\) and the variation \(\nu_{jk}\) of \(\mu_{jk}\) is a finite measure on \(\OA\).
 Especially, \(\mu_{11},\mu_{22}, \dotsc,\mu_{qq}\) and the so-called trace measure \(\tau\defeq\sum_{j=1}^q\mu_{jj}\) of \(\mu\) are finite measures on \(\OA\). 
 For each function \(f\) belonging to \(\LOK{\mu}\defeq\bigcap_{j,k=1}^{q}\LOK{\nu_{jk}}\), we use the notation \(\int_\Omega f\dif\mu\defeq\matauuuo{\int_\Omega f\dif\mu_{jk}}{j}{k}{1}{q}\).
 For this integral, we write \(\int_\Omega f\rk{\omega}\mu\rk{\dif\omega}\) as well.

\breml{M.8-1}
 Let \(\OA\) be a measurable space, let \(\mu=\matauuuo{\mu_{jk}}{j}{k}{1}{q}\in\Mggqa{\Omega,\fA}\), and let \(f \colon\Omega\to\K\) be an \(\fA\)\nobreakdash-\(\BsA{\K}\)\nobreakdash-measurable mapping.
 For each \(U\in\Cqp\), then \(U^\ad\mu U\) belongs to \(\Mggoa{p}{\Omega,\fA}\).
 Using standard arguments of measure and integration theory, easily one can see moreover that the following statements are equivalent:
\baeqi{0}
 \il{M.8-1.i} \(f\in\LOK{\mu}\).
 \il{M.8-1.ii} \(f\in\bigcap_{j=1}^q\LOK{\mu_{jj}}\).
 \il{M.8-1.iii} \(f\in\LOK{\tau}\), where \(\tau\) is the trace measure of \(\mu\).
 \il{M.8-1.iv} \(f\in\LOK{U^\ad\mu U}\) for each \(U\in\Cqp\).
\eaeqi
 If~\ref{M.8-1.i} holds true, then \(U^\ad\rk{\int_\Omega f\dif\mu}U=\int_\Omega f\dif\rk{U^\ad\mu U}\) for all \(U\in\Cqp\).
\erem

\bibliography{211arxiv}
\bibliographystyle{bababbrv}

\vfill\noindent
\begin{minipage}{0.5\textwidth}
 Universit\"at Leipzig\\
 Fakult\"at f\"ur Mathematik und Informatik\\
 PF~10~09~20\\
 D-04009~Leipzig\\
 Germany
\end{minipage}
\begin{minipage}{0.49\textwidth}
 \begin{flushright}
  \texttt{fritzsche@math.uni-leipzig.de}\\
  \texttt{kirstein@math.uni-leipzig.de}\\
  \texttt{kley@math.uni-leipzig.de}\\
  \texttt{maedler@math.uni-leipzig.de} 
 \end{flushright}
\end{minipage}

\end{document}